# The Edge Geometry of Regular N-gons  (Part 2 for 26 ≤ N ≤ 50)
## G. H. Hughes

We will not repeat the preliminary material from Part 1 except for this introduction. In [H5] (*First Families of Regular Polygons and their Mutations*) we defined the First Family of a regular N-gon to be the canonical regular 'tiles' that survive in the complement of the singularity set W under the outer-billiards map $\tau$. Since all the dynamics of $\tau$ are based on this complement, these tiles play an important role in determining the topology of W. Because this topology is very complex we suggested that it might be fruitful to begin with a study the 'edge-geometry' of these regular polygons to yield some insight into the global topology.

Every regular N-gon has a locally invariant region which will include at least the S[1] and S[2] 'tiles 'in the First Family of N. Typically these regions will be bounded by the 'major' resonances of N  so they would be expected to include least 1/3 or 1/4 of all the S[k] in the First Family of N. These tiles would be expected to share a similar geometry and dynamics as shown by N = 60 below – where the inner invariant region extends out to S[25].

In [H2] (*Outer-billiards, digital filters and kicked Hamiltonians*) we describe two additional mappings that can be used to generate portions of this edge geometry of $\tau$. These mappings are the digital-filter map *Df* of Chou and Lin [Cl] and the dual-center map *Dc* of E. Goetz [Go]. In the Appendix of [H5] we show how these two mappings can be used to trace the evolution of W.

In [H5] we show that the evolution of the $\tau$-web can be reduced to a simple 'shear and rotation' and in [H2] we show how the *Df* map can also be reduced to a shear and rotation – but on a toral space. The *Dc* map is *defined* to be a shear and rotation in the complex plane  Therefore it is not surprising that these two mappings have singularity sets that are locally conjugate to $\tau$. This supports our premise that this geometry is inherent in the N-gon. The Df and Dc maps offer a significant computational advantage, but their dynamics are quite different

The $\tau$- singularity set W (a.k.a. the 'web') can be generated by iterating the extended edges of the N-gon as shown below. If these extended edges are truncated they form classical star polygons and the cyan First Family S[k] tiles are defined to be 'conforming' to these nested star polygons. For N = 14 the maximal S[k] is S[5] – also known as D. Here it is congruent to N but for N-odd it is a 2N-gon with edge length identical to N. This D tile is always globally maximum among convex tiles that can evolve in W (with any extended edge length), and rings of these D tiles guarantee that this resulting 'generalized star polygon' is invariant, so it is sufficient to iterate just the star polygon edges to define what we call the web W. If full extended edges are iterated, concentric rings of these D tiles will guarantee global stability (and introduce no additional scaling or geometry). See [VS].

**Figure 1**  The web development for the regular tetradecagon known as N = 14.

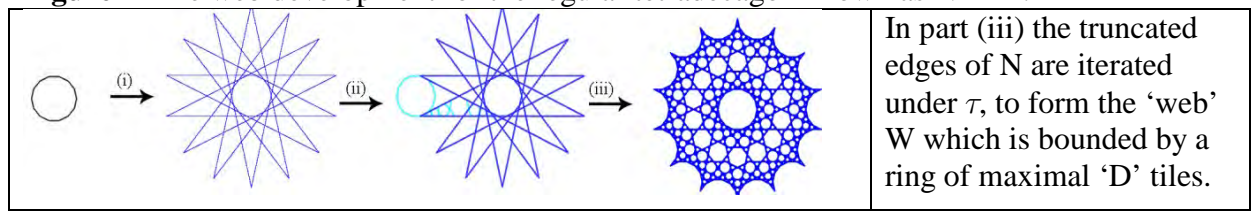

In part (iii) the truncated edges of N are iterated under $\tau$, to form the 'web' W which is bounded by a ring of maximal 'D' tiles.

These 'generalized star polygons' in (iii) share the same scaling and dihedral symmetry as N so they are inherent in N and their geometry is determined by the matching cyclotomic field $\mathbb{Q}_N$. Each S[k] in the First Family of N defines a star[k] point and a scale[k]. The number of 'primitive' scales (gcd(k,N) = 1) matches the 'algebraic complexity' of N, namely $\phi(N)/2$ where $\phi$ is the Euler totient function. This is the rank of the maximal real subfield $\mathbb{Q}_N^+$ of $\mathbb{Q}_N$.

Based on a 1949 result of C.L.Siegel communicated to S. Chowla [Ch], the primitive scales form a unit basis for $\mathbb{Q}_N^+$, so this is what we call the 'scaling field' of N. The traditional generator of $\mathbb{Q}_N^+$ is $\lambda_N = 2\cos(2\pi/N)$ which is $\zeta + \zeta^{-1}$ where $\zeta = \exp(2\pi i/N)$. We will typically use primitive scales as alternate generators of $\mathbb{Q}_N^+$ because the resulting scaling polynomials will be more meaningful than the generic polynomials in $\lambda_N$.

Our ultimate goal is to understand the topology of W, but currently the only cases where the topology of W is understood are the 'linear' or 'quadratic' cases of N = 3,4,5,6,8,10 and 12. The N = 14 case described here has $\phi(N)/2 = 3$ so it is classified as 'cubic' and there are two non-trivial primitive scales along with scale[1] = 1, so W is probably multi-fractal.

As noted by J. Moser in 1978, this $\tau$-web can be regarded as a discontinuous version of the phase space for a Hamiltonian system, so all three maps will have geometry related to the classical 1969 Standard Map of Chirikov. This connection is outlined in [A3] of [H5]

At this time the $\tau$-representations are the most meaningful since there is a well-developed theory of dynamics for both regular and non-regular N-gons. But for computational purposes the Df and Dc maps are simpler and more efficient and we will sometimes use these alternative maps to generate W.

Our primary concern here is the geometry local to N and this will always include the S[1] and S[2] tiles of N, but as noted above, this geometry is typically shared by adjacent tiles in an invariant region local to N. In this way the edge geometry interacts with the overall geometry of N. In Appendix II of [H5] we describe the global evolution of W based on concentric rings of D tiles and show that the initial Ring 0 guarantees that the canonical 'inner star' region is invariant. But invariance exists at all scales and typically this inner region is sub-divided into secondary regions as shown below for N = 60. For N twice-even these secondary rings are fairly well behaved, but for N odd or twice-odd they can be difficult to predict.

**Figure 2** N = 60 has 7 invariant regions which are (almost) symmetric about the S[28] 'M' tile

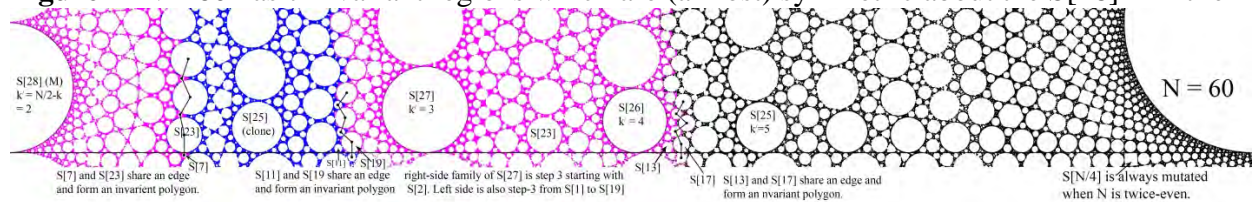

These regions are defined by invariant polygons in the same fashion as the D rings. There are a total of 6 such polygons for N = 60 and the three right-side polygons are shown here. They are formed from the centers of adjacent tiles which share an *edge*. This is a 'stronger' link than the vertex sharing of D tiles and it seems to be a necessary condition for invariance of these regions.

# Catalog of singularity set edge geometry for regular N-gons with $26 \leq N \leq 50$

**Table 1** The 8 dynamical classes of edge geometry based on the Rule of 4 for N even (top) and the Rule of 8 for N odd (bottom).

| 8k family (8,16,24,32,40,48) | 8k + 2 family (10,18,26,34,42,50) | 8k + 4 family (4,12,20,28,36,44) | 8k + 6 family (6,14,22,30,38,46) |
|---|---|---|---|
| N = 16 | N = 18 | N = 20 | N = 22 |
| **8k + 1 family** (9,17,25,33,41,49) | **8k + 3 family** (3,11,19,27,35,43) | **8k+ 5 family** (5,13,21,29,37,45) | **8k+7 family** (7,15,23,31,29,47) |
| N = 17 | N = 19 | N = 21 | N = 23 |

These even and odd cases are linked because any odd N-gon can be regarded as embedded in the 2N-case. Therefore the 8k+1 and 8k+5 families can be embedded in 8k+2 families while the 8k+3 and 8k+7 families can be embedded in 8k+6. This embedding of an odd N-gon in 2N is algebraically sound (with the proper transformation), but the dynamics are different and the web local to N may be very different from the web local to D or 2N.

The **8k+2 Conjecture** says that these N-gons will have an edge geometry driven by sequences of self-similar D[k] and M[k] tiles with known geometric and temporal scaling, The **8k+4 Conjecture** says that the geometry is driven by the mutation of S[2]. The **8k+1 Conjecture** says that these families will have a volunteer DS[2] to go along with the predicted DS[5]. The **8k+7 Conjecture** says that the predicted DS[3] will generate dual DS[1]s with step-2 webs at S[N-3]. The **Twice-even S[1] Conjecture** says that since S[1] has a step-2 web it can support 'step-2'tiles called Skx which are D tiles relative to S[k] (like N-odd).These Skx include S[2] and S3x is an S[2] tile of S[3] for $N \geq 12$.

Every N-gon has a local web which is invariant and this web would be expected to contain at least 1/4 of the S[k], so there is a link between edge geometry and the large scale geometry. Both are driven by the cyclotomic field and the corresponding scaling field $S_N$ with complexity $\varphi(N)/2$. Hopefully the examples below may shed some light on the issue of 'nature' (algebraic complexity) vs. 'nurture' (web and edge complexity under $\tau$ ). As N increases there appears to be a surprising amount of diversity within the 'algebraic families' shown below.

**Table 2** Algebraic Complexity of regular N-gons for $N \leq 50$

| $\phi(N)/2$ | 1 | 2 | 3 | 4 | 5 | 6 | 8 | 9 | 10 | 11 | 12 | 14 | 15 | 18 | 20 | 21 | 23 |
|---|---|---|---|---|---|---|---|---|---|---|---|---|---|---|---|---|---|
| N | 3 | 5 | 7 | 15 | 11 | 13 | 17 | 19 | 25 | 23 | 35 | 29 | 31 | 37 | 41 | 43 | 47 |
|  | 4 | 8 | 9 | 16 | 22 | 21 | 32 | 27 | 33 | 46 | 39 |  |  |  |  | 49 |  |
|  | 6 | 10 | 14 | 20 |  | 26 | 34 | 38 | 44 |  | 45 |  |  |  |  |  |  |
|  |  | 12 | 18 | 24 |  | 28 | 40 |  | 50 |  |  |  |  |  |  |  |  |
|  |  |  |  | 30 |  | 36 | 48 |  |  |  |  |  |  |  |  |  |  |
|  |  |  |  |  |  | 42 |  |  |  |  |  |  |  |  |  |  |  |

●**N = 26**  N = 26 has algebraic complexity 6 and is the first 8k+2 case which is neither quadratic nor cubic in complexity so we used it in Part 1 to help explain the 8k+2 Conjecture. That explanation will be repeated here followed by a full analysis of the geometry. The Edge Conjecture predicts that S[1] will be DS[11] and there will also be DS[7]s and DS[3]s as shown in Fig. 26.4 and Fig. 26.1 below where blue center lines track the local web symmetry. This symmetry also applies to the First Family of N as shown by the S[2] line of symmetry below.

**Figure 26.1** The geometry of N = 26 showing 'clusters of D[2] tiles anchored by S[3]s (not shown). The blue lines are lines of symmetry

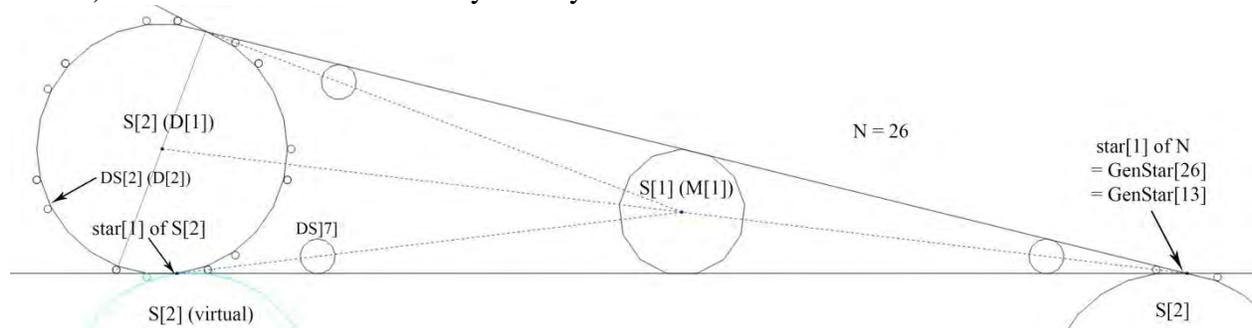

Because the combined web of S[1] and S[2] is step-4 there will be k (from 8k+2) of these S[3] 'clusters' on either side of this line of symmetry. Each cluster has 2 D[2]s and there a lone D[2] on the line of symmetry for a total of 4k+1 D[2]s. This is the 'n = N/2' in the difference equation below. Since the period of any S[k] is N/gcd(k,N), S[2] has period n also. Therefore the period of the D[2]s is $n^2$. Every member of the 8k+2 family will have D[1] at period n and D[2] at period $n^2$. These are the initial conditions of the p[k] equation. (Note that the S[2] local to N is in a different group of S[2]s so its small D[2] is not counted here. However when D[1] is replaced with D[2] to count the next generation D[3]s, these 'outliers' must be counted.)

**Figure 26.2** At D[1] there will be 13 magenta D[3] around each D[2] - each with their outliers. So there will be 14 outlier D[3]s which are just visible here as magenta dots.

The periods of the first three D[k] are:
D[1] at 13, D[2] at $13^2$ and D[3] at $13^3 + 14 \cdot 13$
so the proposed difference equation for the periods $D_k$ of the D[k] is
$$D_k = nD_{k-1} + (n+1)D_{k-2}$$
where n = N/2 and $D_1 = n$ and $D_2 = n^2$
Mathematica's solution is:

**RSolve**[{D[k]= D[k-1]n + D[k-2](n+1), D[1]=n, D[2]=n^2}, D[k], k]

$$D[k] \to -\frac{n\left((-1)^k - (1+n)^k\right)}{2+n}$$

as in the 8k+2 Conjecture

Therefore: $D[n,k] = -\dfrac{n\left((-1)^k - (1+n)^k\right)}{2+n}$ gives the $\tau$-period of each D[k] for N = 2n

**Table**[D[13,k], {k,1,8}] = {13, 169, 2379, 33293, 466115, 6525597, 91358371, 1279017181}

As noted in the 8k+2 Conjecture the periods of the M[k] have the same basic equation but different initial conditions of M[1] = n and $M_2$ = period of M[2] = n(3(n-1)/2 + 2  to give

$$M_k = \frac{n\left((-1)^{1+k} + (-1)^k n + 3(1+n)^k\right)}{2(2+n)}$$  **Table**[M[13,k], {k,1,5}] ={{13, 260, 3562, 49946, 699166}}

It is easy to find the centers of these tiles and verify these periods. Just the D[k] centers will suffice because these are the reference for the other centers. These centers must lie on the line joining the center of S[1] with the local limit point which is star[1] of S[2]. This point has known coordinates of $\{x_0, y_0\}$ = {StarS2[[1]][[1]],-1}and the center line has slope that we call cslope. The center of D[k] will have coordinates $\{x_k, y_k\}$ where $x_k = x_0 + hS[2]*GenScale^{k-1}$/cslope and $y_k = hS[2]*GenScale^{k-1}$-1; where for N twice-odd, GenScale = GenerationScale[N/2].

Then cM[k] = TranslationTransform[cD[k]][cS[1]*GenScale$^k$]; and this works for any DSk. For example with N = 26, cDS7[k] = TranslationTransform[cD[k]][cS[7]*GenScale$^k$] so the first few periods of DS7[k] are : 52 (13*4) , 858, 11882, 166478 . It is not surprising that these will satisfy the same difference equation as the D[k] except with initial conditions 4n and n(4n+n+1) (where n = N/2). The case for DS3 is very similar with initial conditions 6n and n(6n+n+1) and the M[k] at DS[11] will have initial conditions 2n and n(2 + (3/2)(n-1)) as noted earlier. All of these ratios will approach N/2+1, and the geometric scaling will be GenScale(n).

This explains our conjecture that the entire 'generations' between D[k[and M[k] will share the same fractal dimension – even though these regions are not globally invariant. But all N-gons appear to have a globally invariant region local to N and here for N = 26 the region extends out to S[6]. These regions would be expected contain locally invariant regions on all scales. For N even the edges of the S[2] tiles should be locally invariant because the 'even' and 'odd' S[2] map to each other. Our periods here are based on just the 'odd' S[2] and that is sufficient.

**Figure 26.3** The region around S[1] and S[2] for N = 26 is globally invariant out to S[6], but invariance exists on all scales and the blue and magenta edges of the S[2] appear to be invariant.

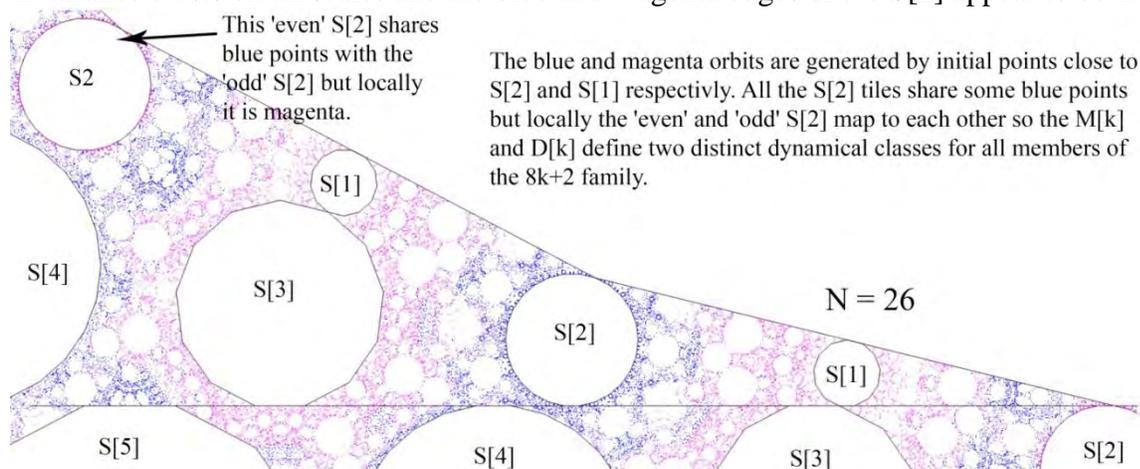

This plot makes it clear that the M[k], D[k] convergence at star[1] of N is magenta and not blue, but there is clearly a magenta version of the 8k+2 Conjecture that describes this equivalent convergence, so we can concentrate here on star[1] of S[2] and the region between S[1] and S[2].

The 8k+2 Conjecture makes some very strong predictions about the dynamics of this region and since N = 26 is relatively accessible we will attempt to explore the resulting geometry in as much detail as practical. In particular we will probe the first 5 generations below.

For N = 18 it appears that each generation defined by M[k] and D]k] is self-similar to the 2$^{nd}$ generation but that is clearly not true in general. Because of the cw-ccw dichotomy of each D[k] there will always be some variation between even and odd D[k] generations and the major issue is whether there is self-similarity within these two classes. What we will see here, and most likely everywhere in the 8k+2 family beyond N = 18, is a similarity within even and odd generations, but not perfect self-similarity of any two generations. This prediction is now part of the revised 8k+2 Conjecture.

**Figure 26.4** The early web of N = 26 showing blue lines of symmetry. We expect that these lines of symmetry will scale with generations in the same fashion as the DS[k].

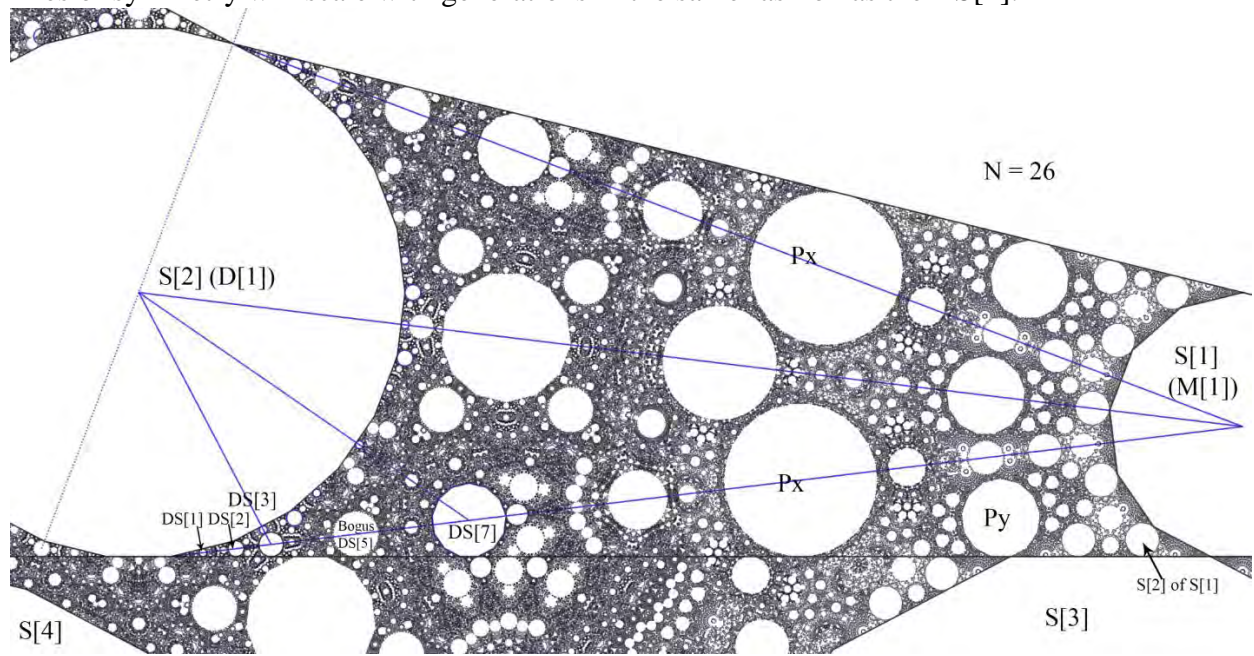

The Rule of 4 predicts the DS[3]s and DS[7]s along with S[1] at S[11]. There is also a close match with the DS[5] of S[2]. Since N is even these DS[k] are simply the S[k] in the First Family of S[2]. It is easy to derive the parameters of the Px and Py volunteers .It appears that the Py may be preserved in future generations but the Px tiles will not survive in their current form.

As noted earlier, because the combined S[1]-S[2] web is step-4, there will always be k clusters of DS[3]s on each side of the line of symmetry and one shared DS[2]. Here k is 3 so there will be 13 DS[2]s just barely visible above. These DS[2]s can serve as the next-generation D[2]s and foster DS[3]s on their edges. Based on Fig.26.2 earlier, the count for the DS[3]s should be 13 for each of the 13 D[2]s but in addition each 13 D[2]s has the potential for a left and right-side D[3] as part of their families. However because the D[2]s come in clusters of two, these outliers are shared and the total is just 14. Therefore the total count for D[3]s is $13^3 + 13*14 = 2379$.

Below is the web local to one of these D[2]. Since these D[2]s are vertex tiles of S[2] it is natural that their web rotations are now cw as shown below. But since the matching M[k] exist the new joint web will be step-4 just like D[1]. Therefore the local geometry is very similar with 13 D[3]s with plus the outlier on the left.

**Figure 26.3** The 3rd generation presided over by D[2] is similar to the $2^{nd}$ but the web now evolves cw which is opposite of D[1]

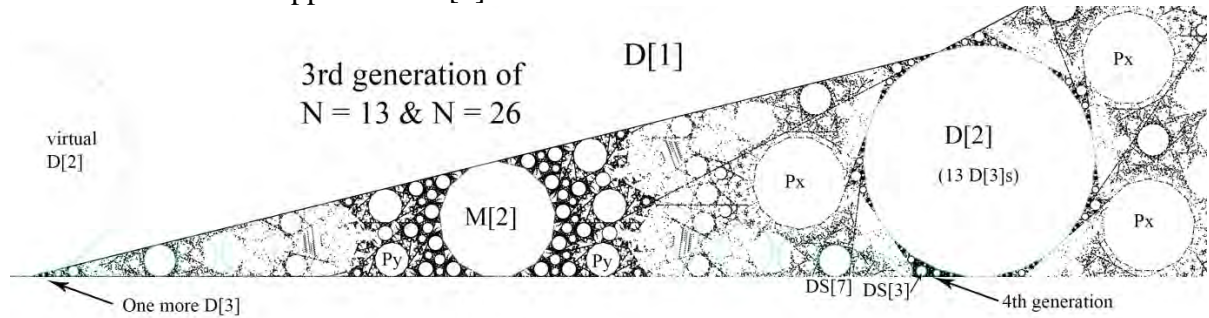

It appears that in 'most' cases there is no global self-similarity of these generations. Here with N = 26 it is clear that the $3^{rd}$ generation is quite different from the earlier generations and below we will see that the $4^{th}$ generation is also distinct – but naturally has more in common with the $2^{nd}$ then the $3^{rd}$.

**Figure 26.5** The $4^{th}$ generation presided over by D[3] and M[3].

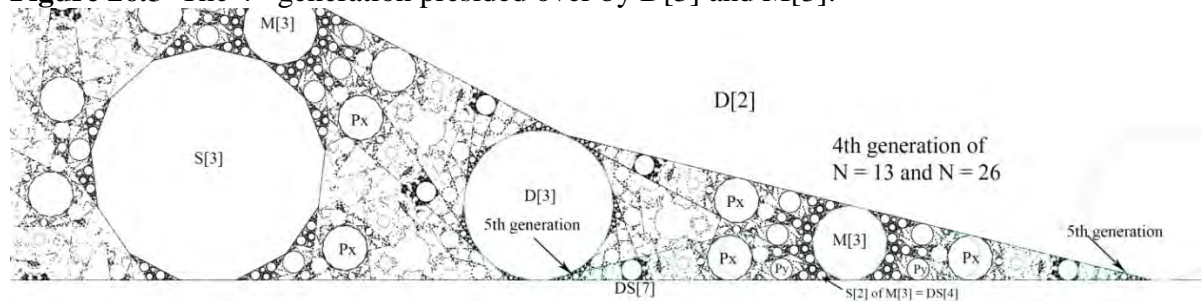

The Px tiles shown here are not just scaled versions of the Px from the $2^{nd}$ or $3^{rd}$ generations. Of course their heights are in the scaling field but the various Px polynomials have no obvious relationship with each other, while the Py local to M[3] are simply scaled versions of each other.

Below is the right side of S[1] where the parameters of the original Py can be derived. Each star[k] point of S[1] defines a potential S[k], but the only one that exists here is the S[2] of S[1]. This tile is a clone of S[4] of D[1]. The virtual S[5] of M[1] is the same as S[10] of D[1] and it is the only menber of of the First Families of both D[1] and M[1] – but it is our convention to combine these into a single shared 'First Family' of M[1] or D[1]. If all the First Family tiles of M[1] are included they can be regarded a part of the 'extended' family of M[1] or D[1].

**Figure 26.6** A vector plot of the geometry on the right-side of S[1]

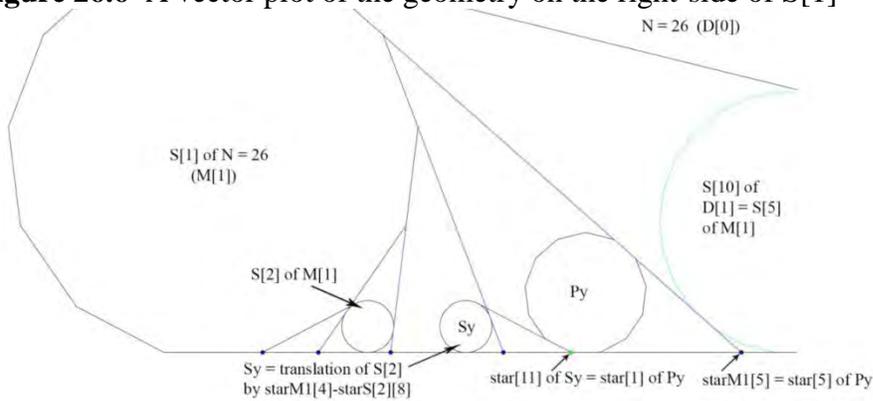

As shown here the missing green point that definesf Py can be found using the displaced S[2]. This geometry appears to survive in future generations. Of course the M[k] scale along with the D[k], but they have high-low ratios instead of low-high, because the bulk of the M[k] are satellites of the D[k-1]. Here M[6] has period 6884150*13 and it has ratio 14.00008 with respect to M[5]

For N = 26 , S[11] is the N = 13 proxy but it is in another invariant region with very different dynamics, so it is not clear what we can conclude about N = 13 based on the geometry here. The Twice-Odd Lemma implies that S[11] will be a surrogate N = 13 but this applies only to the 1$^{st}$ generation and makes no conjecture about M[1] above being a surrogate 2$^{nd}$ generation N = 13 because D[1] is not a clone of N. It is true that hM[1]/hD[1] is identical to hM[0]/hN but the local geometry of M[0] (N = 13) may be very different from M[1] above.

**Note**: In the analysis of the first 50 N-gons to follow, the 8k+2 family appears to be the only one where we can make predictions about the evolution of generations. The web development of most all N-gons is so complex that it is impossible to make long-term predictions and we are left with brute force web iterations to resolve detail beyond the first few generations. For an N-gon like N = 48, to resolve detail for the 3$^{rd}$ generation D[2] tile on the edges of S[2] will require at least 100,000 local web points, but the Dc map 'return rate' is on the order of 1 in 500,000 iterations because hD[2] = (GenScale[N]/scale[2])$^2$ ≈ .0000075. This means that the image will require about 50 billion iterations. A reasonably fast I7 computer takes about 90 minutes for 1 billion iterations so we were left with 3 days of round-the clock iterations for a single image, but our new I9 HP workstation running compiled code in Mathematica 14 can do this in about 1 day, Of course this is just a just a temporary reprieve since 1/GenScale[N] scales by more than 1.8.

In [H5] we noted that for N =11, 'generations' in the scaling field $S_{11}$ scale by GenScale[11] = Tan[π/11]·Tan[π/22] ≈ .042217, so the 25th generation would be at the Plank scale of 1.6·10$^{-35}$. In the words of Richard Schwartz at Brown, "A case like N =11 "may be beyond the reach of current technology".

- **N = 27**

N = 27 is the 3rd non-trivial member of the 8k+3 family along with N = 11 and N = 19. N = 27 has complexity 9 and this matches the complexity of N = 19, so these two are in the same 8k+3 family and also have the same algebraic complexity. This 8k+3 family has been a difficult one to characterize because the first predicted DS[k] is DS[7]. However in the 8k+3 family there does seem to be a propensity of conforming volunteers between predicted DS[k]. Both N = 19 and N = 27 have a $D_1$ volunteer which could play the role of a surrogate M tile between S[1] and S[2]. When it exists, this M tile is the only tile shared by S[1] and S[2]. N = 27 does indeed have a $D_2$ volunteer between DS[15] and DS[7], and here may be a volunteer $D_3$ in a manner similar to N 19. As expected the volunteer $S_k$ tiles of S[1] occur primarily at the mod-4 star points and this includes the shared $D_3$ at star[N-2] of S[1] = star[2] of S[2].

**Figure 27.1** – The web of N = 27 showing the early web in magenta.

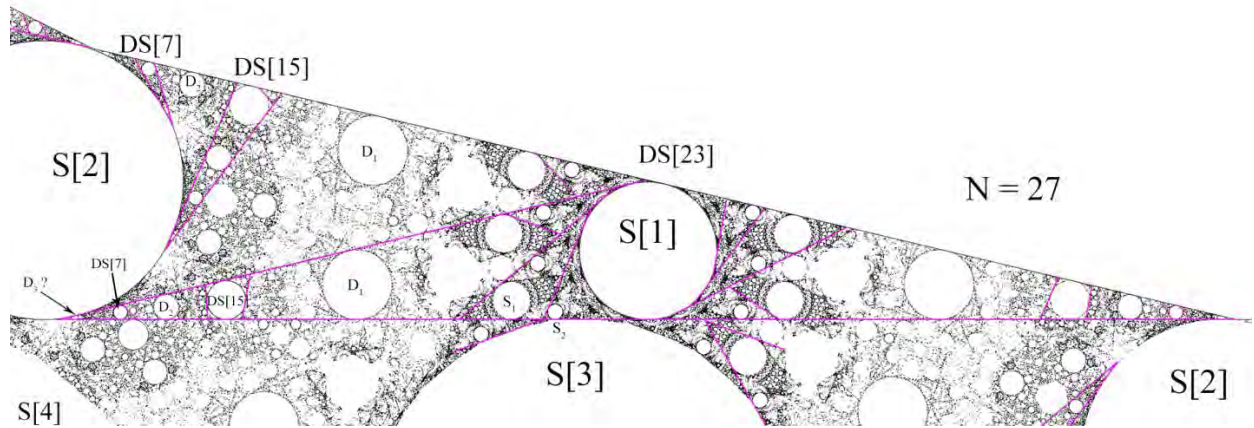

**Figure 27.2**  The web local to S[2]

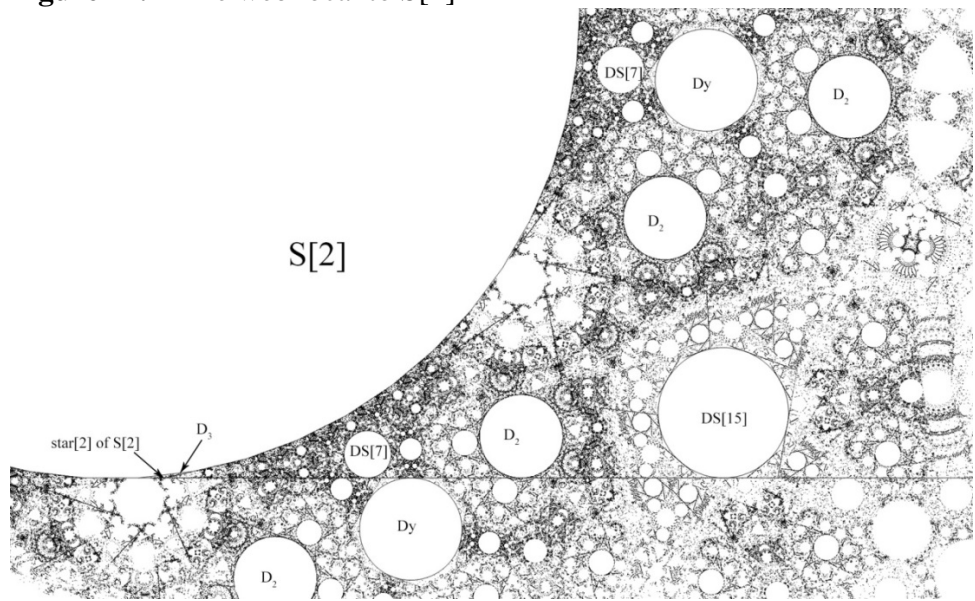

**Figure 27.3** Detail of the web local to S[2]

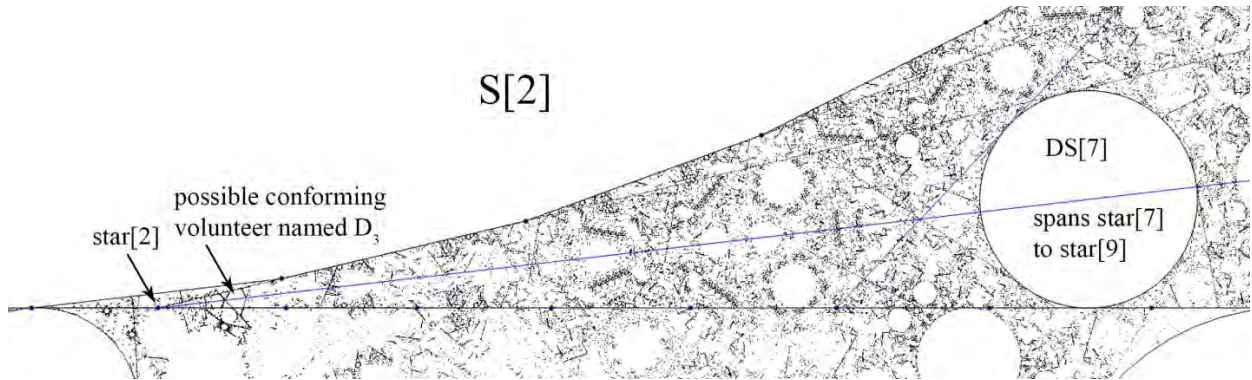

It seems that $D_3$ is mutated in the form of a hexagon and this is consistent with it being the weave of two triangles. Typically the star[2] point of S[2] should support convergent sequences of conforming tiles. The $D_4$ tile below might be in the First Family of $D_3$

**Figure 27.4** Detail of the star[2] region of S[2]

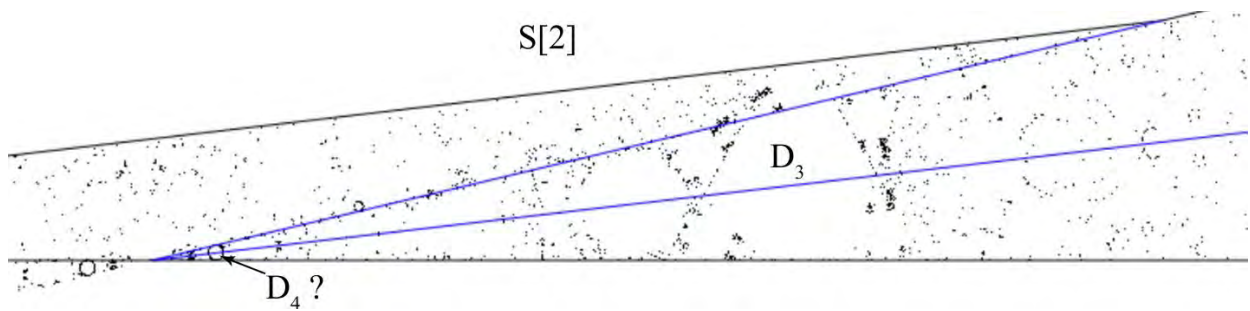

When N is odd the web geometry dictates that the 'GenStar' of S[1] is the penultimate star[N-2] and the 'effective' star points which survive the early web are always mod 4. Therefore star[1] of S[1] is effective in the 8k+3 and 8k+7 families. This does not appear to enhance the geometry local to S[1] because there is an large gap from star[1] to star[4]. There is more evidence of local structure when star[3] of S[1] is effective as in the 8k+1 and 8k+5 families. When N is even, the local web of S[1] is step-2 and the gaps are smaller. Since S[1] now shares its S[N/2-1] point with star[1] of S[2], star[1] of S[1] will always be effective.

**Figure 27.5** The web local to S[1]

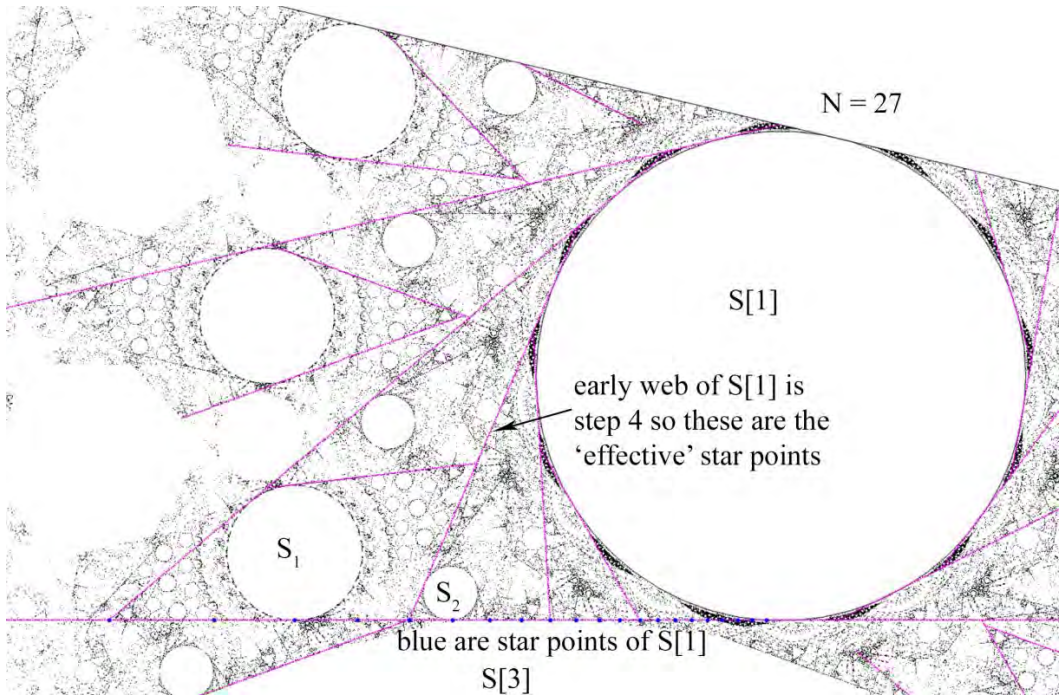

**Figure 27.6** Detail of the web local to S[1]

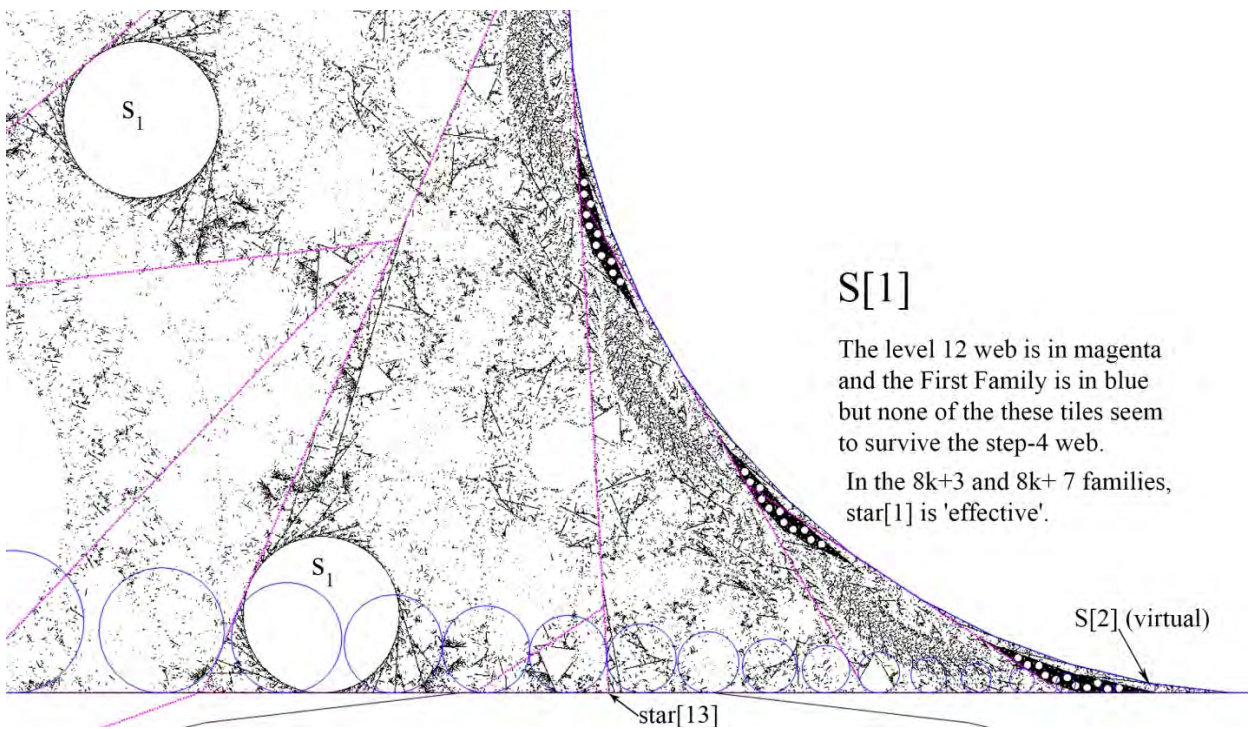

- N = 28

N = 28 has complexity 6 and is the 3$^{rd}$ member of the 8k+4 family. The Mutation Conjecture predicts that for N twice-even, the S[k] tiles will be mutated iff gcd(N/2-k,N) > 2, so the 8k+ 4 family is the only family where the S[2] tile is mutated. Because S[2] will have gcd(N/2-k,N) = 4, the web cycle will be also be step-4 and S[2] will be the weave of two regular N/4-gons.

**Figure 28.1**  The level 10 web in magenta and limiting web in black

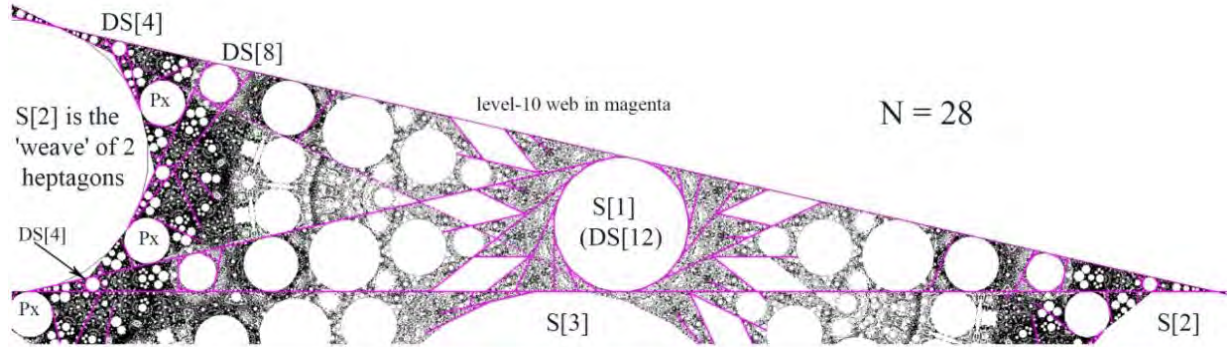

Below is a schematic of the mutation. The Mutation Conjecture states that the equilateral edges of the resulting weave will span 4 star points of the underlying S[2] because gcd(k′,N) = 4. The most common cases are when one of the star points is a star[1] of the underlying S[k] and this is true here, with the 'base' spanning right-side star[1] to left-side star[3]. The general rule as outlined in Section 2 of Part 1 is that for N-even the minimum index should be N/2-1-jk because this mimics the choices the 'lazy' web has. Here the minimum is 17-2j = 1.

**Figure 28.2**  The mutation of the S[2] tile. These 'star point' mutations will be 'canonical' because their side length s will be the difference of star points. We conjecture that all tiles that arise in the web must be canonical with s/s$_N$ in the scaling field S$_N$.

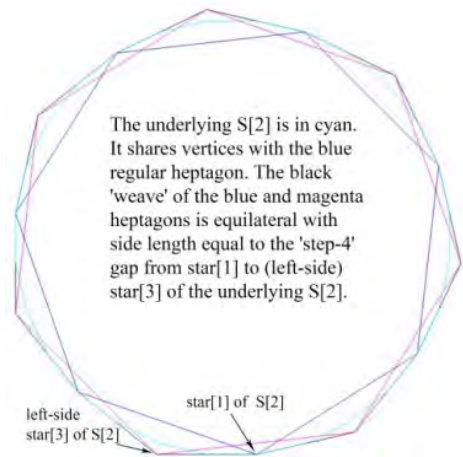

This mutation is consistent with the Rule of 4 which predicts that DS[k] will exist at least mod-4 counting down from S[1] at DS[N-4]. Here the predicted DS[k] are DS[8] and DS[4] and these are members of the 'normal' First Family of the underlying S[2].

For the 8k+4 family this means there will always be DS[4]s which are consistent with the mutation of S[2]. These DS]4\s will 'almost' share a vertex with the extended edges of S[2]. As shown below, this small gap creates a very 'awkward' local geometry around the DS[4]s and prevents them from having a predictable local family structure.

The 8k+4 Conjecture predicts that there will be two possible consequences of this break-down. N = 12 and N = 28 illustrate one case where there will be a volunteer Px tile whose reflection about star[1[ of S[2] will be a 'parent' of DS[4]. The other case is illustrated by N = 20 and N = 36 is where this 'parent' Px may be real or simply the virtual D tile of DS[4] which always exists for any DS[k] of S[2]. (Initially we conjectured that this case would always yield a real Px as found in N = 20 and N = 36, but looking forward to N = 52 shows that this is false.) In both cases this Px has the potential to generate a $3^{rd}$ generation on the edges of S[2].

**Figure 28.3** Detail local to S[2]

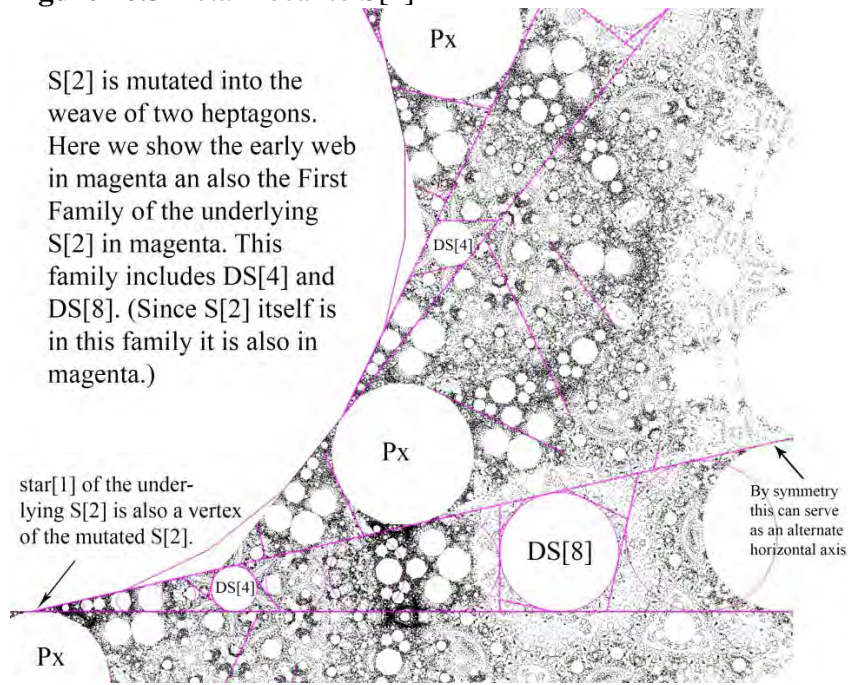

Here the blue vertices support volunteer Px tiles. The 8k+4 Conjecture says that the 8k+4 family will split into two mod-16 subfamilies, where N = 12 is the initial tile in one branch, followed by N = 28, N = 44 and N = 60 of the Introduction. The S[2] tiles in this branch will have Px tiles at the blue vertices and these will be 'parent' tiles of the DS[4]s. Here we will show that DS[4] is the S[11] tile of Px. The second branch of the 8k+4 Conjecture is anchored by N = 20 and for these cases Px is a D tile relative to DS[4].

In terms of edge geometry N = 28 can be regarded as a 'heptagon' version of N = 12. Our convention for numbering the vertices of S[2] will be ccw starting with $v_1$ at star[1] which is in both the mutated and un-mutated S[2]. These odd vertices will always be vertices of the underlying S[2] while the even vertices are extended outwards. This implies that the canonical DS[4]s will be roughly 'vertex-tiles' of the large heptagons while the volunteer matching Px will be actual vertex tiles of the smaller blue heptagons. Therefore the local geometry of these Px tiles will determine what we call the 3$^{rd}$ generation at star[1] of S[2]. For N =12 ,the (rotated) S[3] tile of N was the Px, but that was a special case. Figure 28.5 below shows that the (rotated) S[4] tile of N appears to share a vertex with Px, but this is only an approximate alignment. Px is not in the First Family of S[4].

Any S[k] tile of N can generate N because they share at least two star points, so DS[4] can be used to generate Px as follows:

(i) The 8k+ 4 Conjecture predicts that DS[4] will be DS[3 + 8] of Px so it must share star[11] of Px. The local index is N/2-k, so this star[11] point of Px is star[3] of DS[4] (which is known.)

(ii) Since DS[4] is assumed to be conforming to Px, it also shares its GenStar point with star[1] of Px. Therefore star[N/2-1] of DS[4] is star[1] of Px.

(iii) Therefore hPx = d/(Tan[11π/28]-Tan[π/28] where d is the horizontal displacement of star[3] and star[13] of DS[4].

(iv) Once hPx is known any star point can be used to find its midpoint and center. For example using (right-side) star[1] of Px, MidPx = star[1] – {hPx·Tan[π/28],0} (Of course this yields a copy of Px which is embedded in S[2]. But this is a natural embedding and very useful.)

For N = 28 this Px plays the role of a 'vertex based'step-4 tile relative to both S[2] and the mutated 14-gon. This raises the issue of when mutated tiles have 'natural' families. For N = 12, it is clear that the mutated S[2] has a self-similar web which is locally equivalent to the N = 12 web. This is clearly a special case, but there may be subtle relationships in general.

**AlgebraicNumberPolynomial[ToNumberField[hPx/hN,GenScale],x]** gives

$$-\frac{433}{512} + \frac{41903x}{512} - \frac{39573x^2}{256} + \frac{16255x^3}{256} - \frac{3157x^4}{512} + \frac{35x^5}{512}$$

**Figure 28.5** The First Family of Px in magenta showing S[2],S[4] and S[11] as survivors

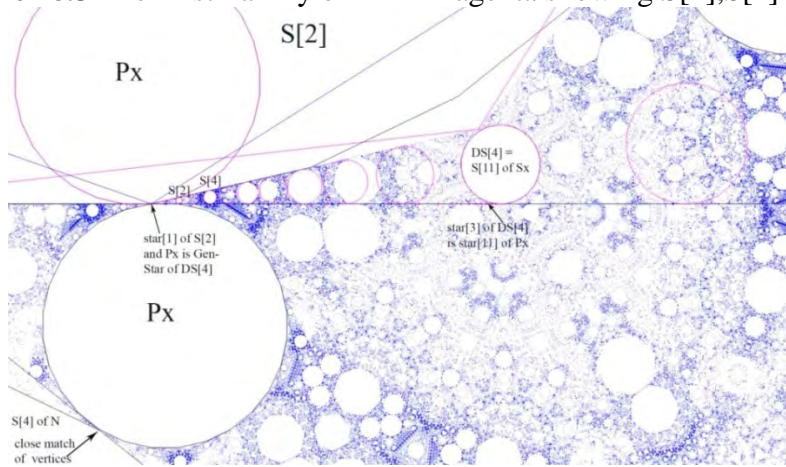

- **N = 29**

N = 29 and the matching N = 58 are the only N-gons with complexity 14. N = 29 is in the 8k+5 family so it will have a DS[1] which will be congruent to an S[2] of S[2], which can serve as a D[2] in a potential chain of D[k]. The existence of this chain in the 8k=2 family is strongly dependent on the existence of matching M[k] which will share a local web with D[k]. Here for N-odd any M[k] are 'out of bounds' because they arise from star[1] and the S[2] web is based on star[2]. This means that any M[k] which exist will be 'volunteers' and not part of the normal 'step-2' DS[k] family of S[2].

The only known 8k= 5 cases where an M[2] exists are N = 5, 13 and 21, so N = 29 is the first 'non-exceptional' case for this family and the question is whether a lone D[2] can generate an extended family without a matching M[2]. N = 16 accomplished this but that was an exception in the 8k family. For N odd it will be much more difficult. For N = 5, S[2] was the D tile in the 8k+2 family where the S[2] –S[1] pairing is preserved in future generations. In this sense N = 5 was totally dependent on N = 10 for support. This is clearly an exception and even N = 13 has no M[2] and no obvious extended family structure, even though N = 26 is just 'next door'.

**Figure 29.1** The major lines symmetry showing DS[1], DS[11] and DS[19]

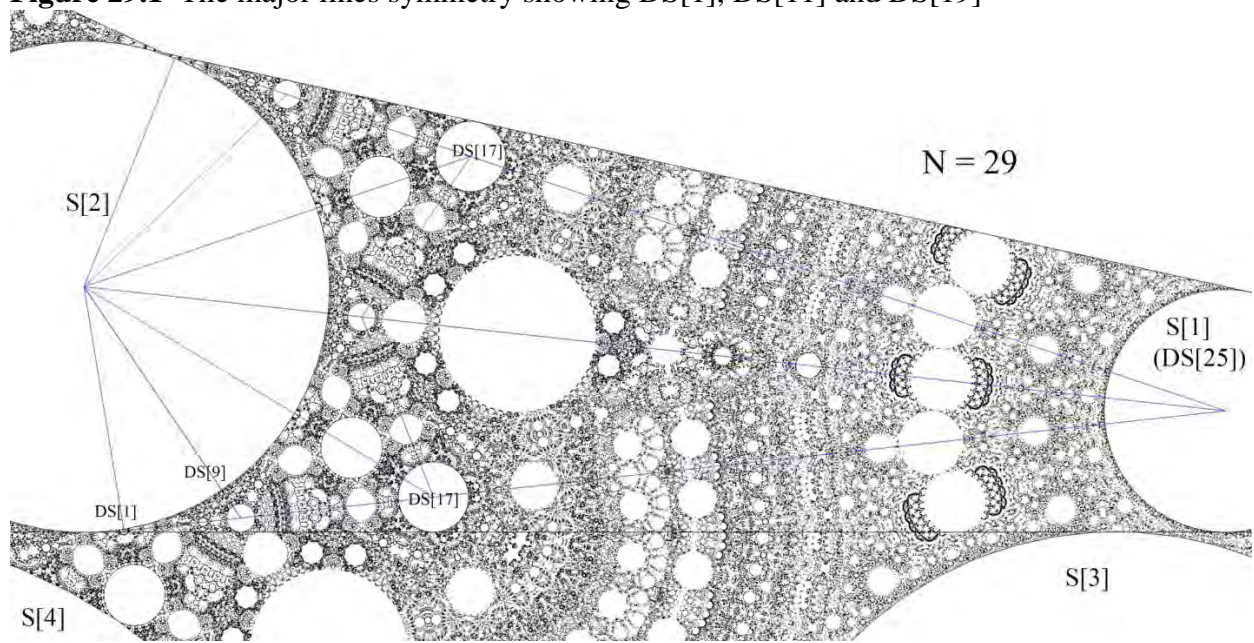

**Figure 29.2** Detail of S[2] showing the early magenta web

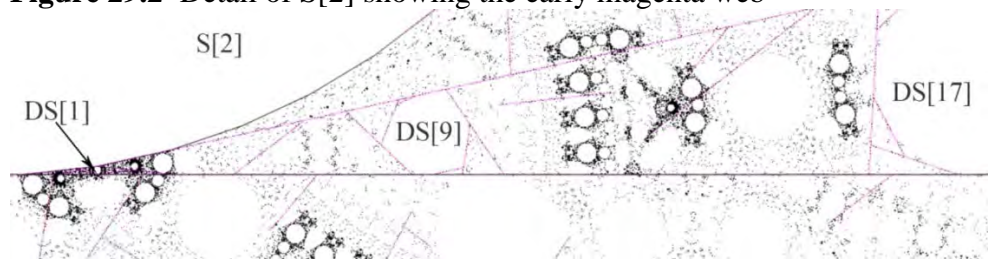

**Figure 29.3** The DS[1] region showing the canonical step-2 cw web of DS[1]

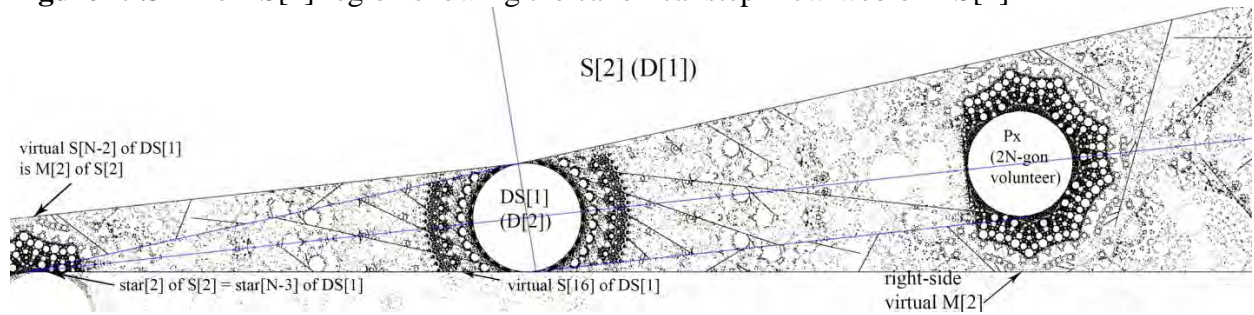

DS[1] is a 2N-gon its First Family has right and left-side S[N-2] tiles that are M[2]'s relative to S[2]. Typically these are virtual in the 8K+4 family but they can still help to explain the local geometry. Unfortunately there are no known survivors of the First Family of DS[1] and we expect this to be generic as N increases. But there is hope for a family of the virtual M[2].

**Figure 29.3** Detail of DS[1] showing the virtual First Family tiles in Cyan

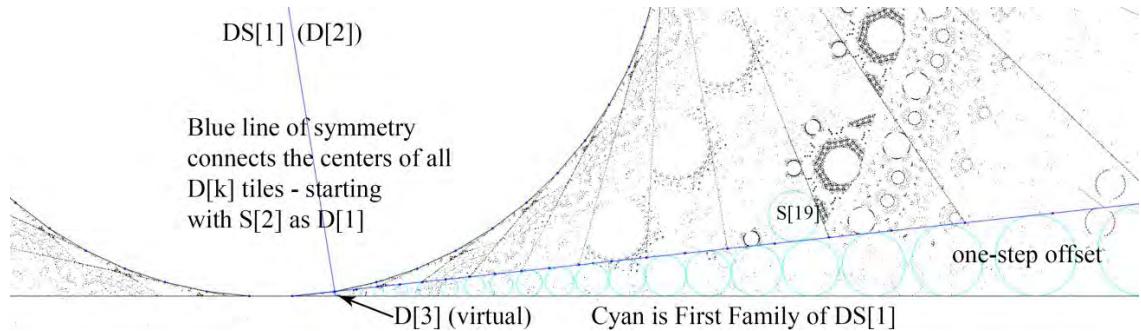

**Figure 29.4** The volunteer Px tile shares an S[2] tile with the right-side virtual M[2]

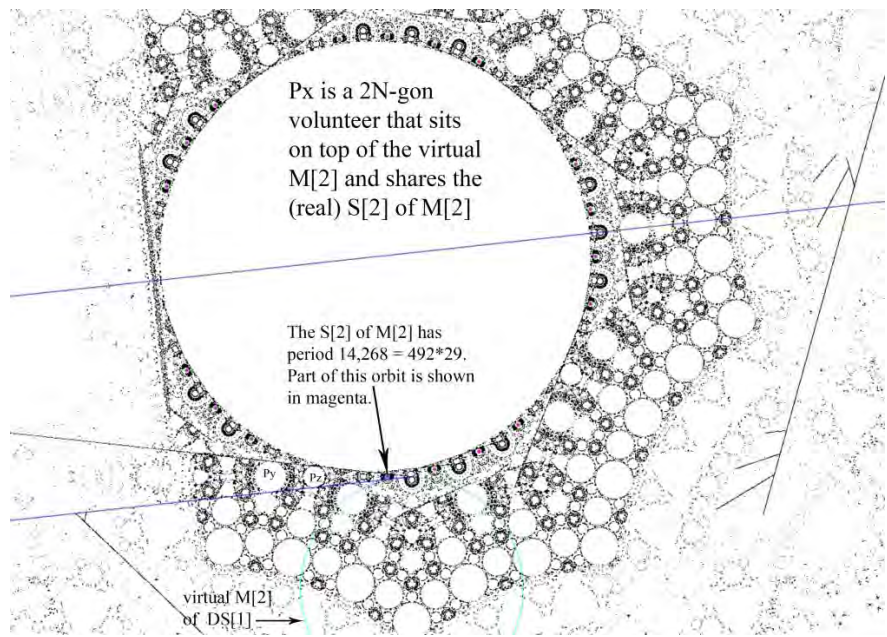

Based on this orbit of the center of S[2] it is clear that Px reflects the geometry local to star[2] of S[2] and the matching star[27] of DS[1]. It seems that the step-2 web of DS[1] forms local regions around DS[1] that look promising for S[k] tiles to form, but the interior points of these regions are pushed outwards into rings which form prior to any local development. This is similar to the way that M[k] tiles evolve.

**Figure 29.5** Detail of the marriage of Px and the virtual M[2]

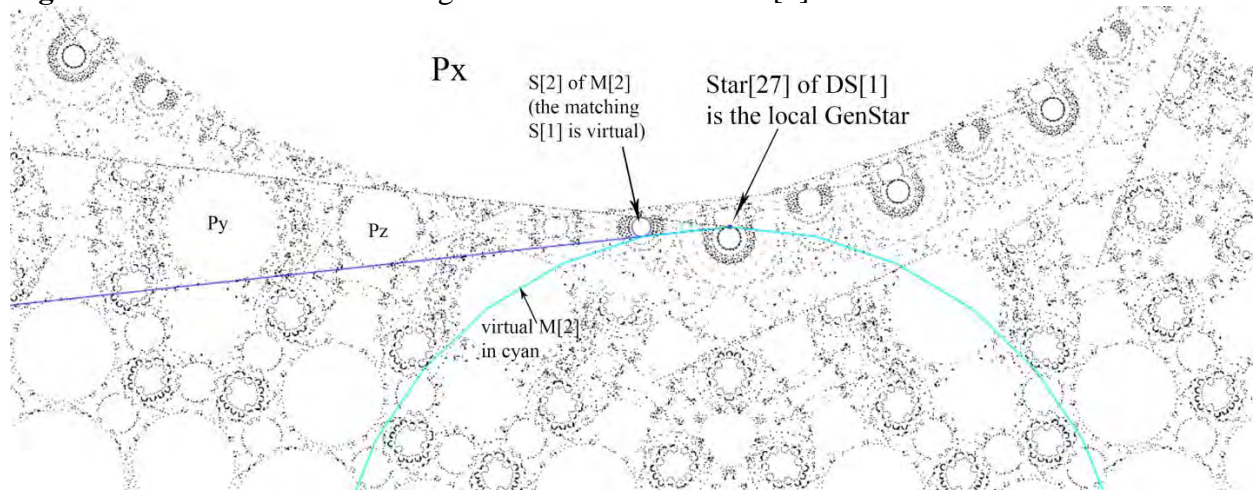

The Py and Pz tiles shown here are clearly part of the Px 'family' so they may have a different origin than S[2]. Their periods are 2900 and 7830 respectively . It is possible that S[2] could foster a 4$^{th}$ generation at the top of M[2] but clearly the matching S[1] is virtual. Note that S[2] has a 'halo' which could be self-similar to DS[1] .Approaching the local GenStar at star[27] of DS[1] the periods grow rapidly but there is not enough detail to guess about the geometry.

**Figure 29.3** The region local to S[1] has no obvious survivors of the First Family of S[1]

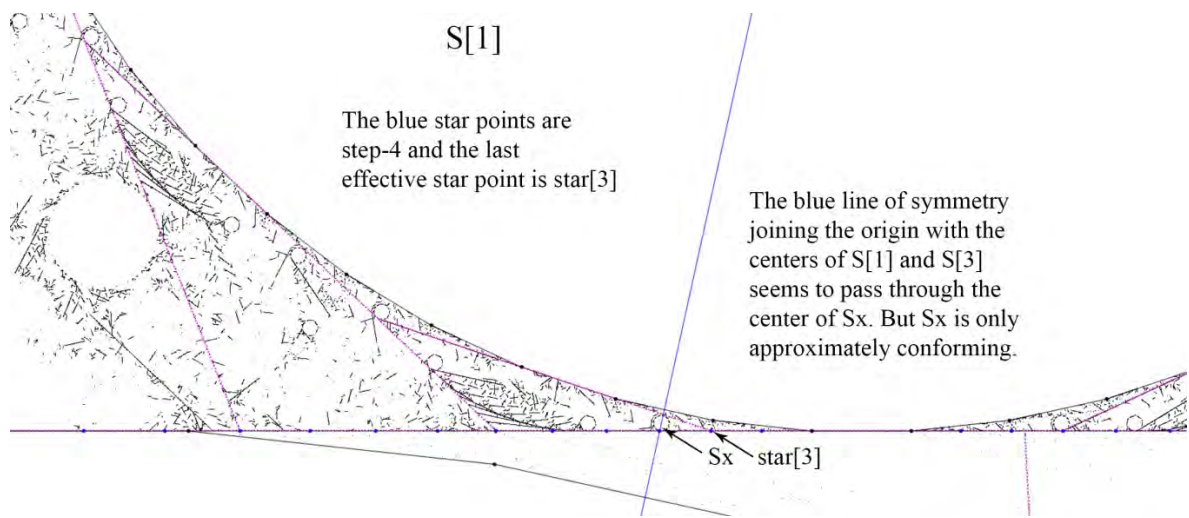

● N = 30

N = 30 has quartic complexity and is a member of the 8k+6 dynamical family along with N = 14 and N = 22. The quartic family consists of N = 15, 16, 20, 24 and 30 and there may be an algebraic connection between these members which is more subtle that the quadratic or cubic cases. The closest connection is between N = 30 and the embedded N = 15 because they share equivalent cyclotomic fields and webs.

**Figure 30.1** The shared First Family of N = 30 and N = 15 with mutations in blue

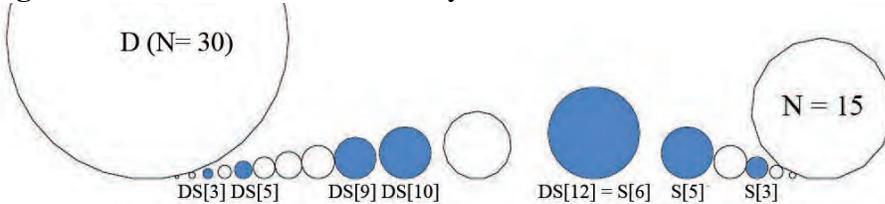

By convention we use the right-side copy of N = 30 as our default N. The Edge Conjecture predicts a step-4 count-down from S[1] at DS[N/2-2] so there should be a DS[9], DS[5] and DS[1]. Since these DS[k] inherit the same k′ = N/2-k as the First Family of N, there will be mutations in DS[3], DS[5] and DS[9] to match the blue mutations above. All these DS[k] have k odd so they will be N/2-gons and the mutation condition is gcd(N/2,k′) > 1, so the volunteer DS[4] is not mutated. DS[3] and DS[5] will have linked mutations and Mx appears to share the DS[9] mutation in a non-trivial fashion.

The prominent volunteer DS[7]s appear at step-8 intervals but not at the normal position on the horizontal 'base' line. This may be due to the influence of S[4]. The volunteer Dx and Mx tiles circle the N/2-gon DS[7]s at step-2 with DS[9] as 'guest'. The two Mx tiles are also rotations about DS[9]. The Dx may form as part of the convex hull of DS[3] and DS[5] as the web evolves. These is more detail in Figure 30.3. Note the DS[7], DS[5], DS[3] odd 'towers' which 'morph' into DS[7]-Dx towers.

**Figure 30.2** Blue lines of symmetry

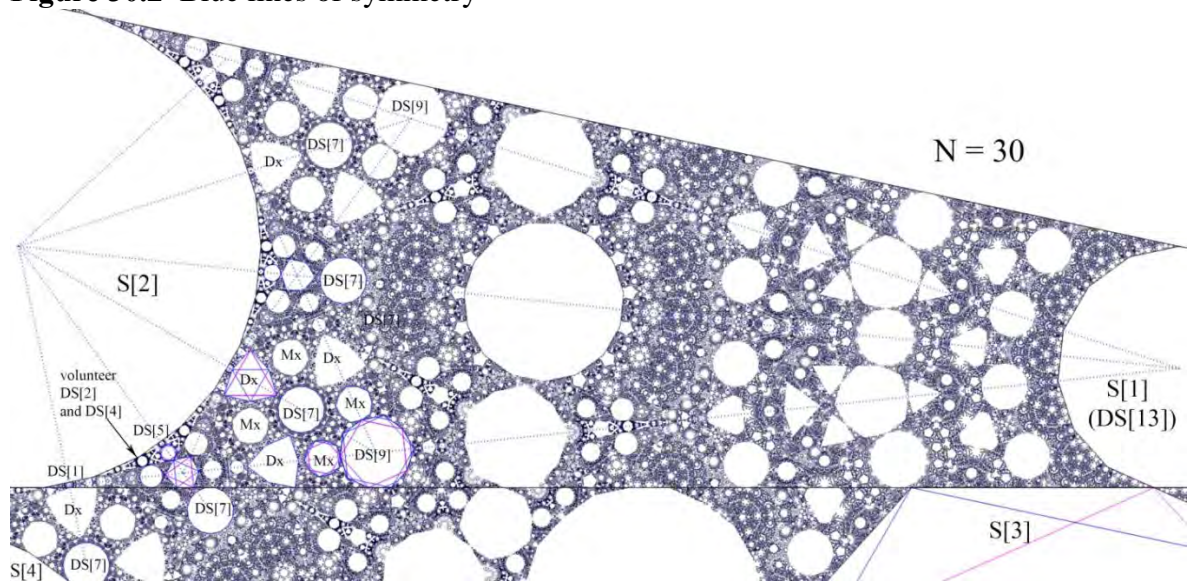

**Figure 30.3** The abundance of volunteers may reflect the abundance of mutations. The step-3 semi-invariant region on the edges of S[2] may be related to the strong step-3 resonance. DS[5] and the DS[4]s are outside the step-3 region but together they define a secondary step-5 region to match the step-5 resonance. N = 54 has a strong step-3 resonance with no sign of invariance so maybe the step-3 and step-5 modes are both contributors.

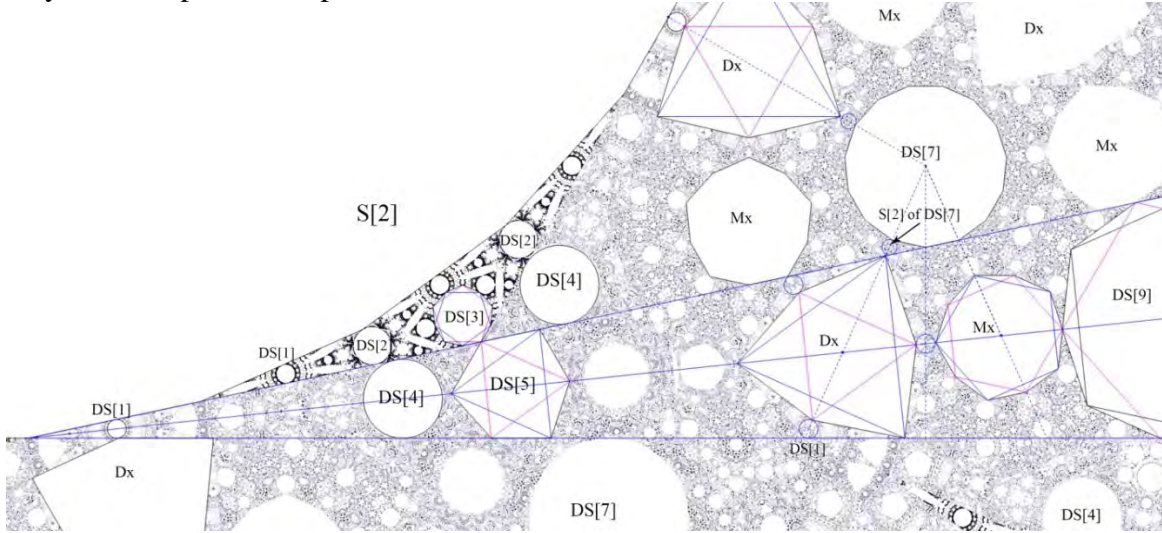

Multiple factors worth trying are N = 150, or 126 with 21 or 198 with 33. It is possible that N = 30 may be quite unique and the quartic nature may be a factor. In investigations by Lowenstein and Vivaldi such as [LV1] the authors found that quartic geometry was topologically not 'well-behaved'.

**Figure 30.4** Detail of the S[2] edge showing the shared-vertex mutations in DS[3] and DS[5]

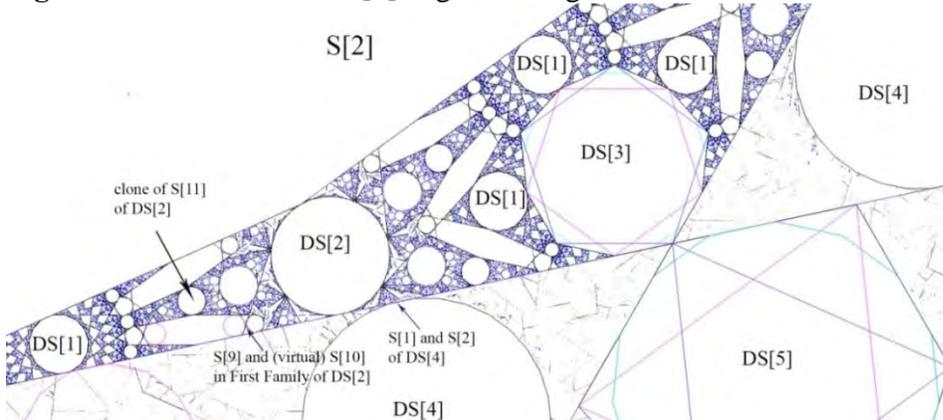

DS[5] and DS[9] have $15/\gcd(15,k') = 3$ and 5 respectively so they will be weave of triangles and pentagons. These are consistent with the mutations of the First Family tiles. DS[3] will have a 'pentagon' mutation to match DS[9] because $15/\gcd(15,12) = 5$, but the volunteer DS[3] shown here is rotated by $2\pi/30$ relative to the First Family DS[3], so it is also rotated relative to DS[9]. The remaining geometry of the invariant region is complex and unique. It does share some features with the DS[5]-DS[7] mutation for N = 15, but the big difference is the volunteer DS[2]s which do not exist for N = 15. Along with the matching DS[1]s, they impart a geometry that has more coherence than anything observed with the S[2] of N = 15. These DS[3] clusters are reminiscent of the 8k+2 family.

● N = 31

N = 31 is in the 8k+7 family and has complexity 15 along with the matching N = 62. Because N is odd the Edge Conjecture predicts that there will be a DS[3], DS[11] and DS[19] to go with S[1] at DS[27]. All of these DS[k] have k odd so they are 2N-gon 'D' tiles relative the normal odd S[k] of S[2] which are N-gons. This means that DS[3] will have a normal 2N-gon (virtual) First Family which includes an S[N-3] tile which is identical to a DS[1] of S[2]. This DS[1] is an S[2] of S[2] so it can possibly act as a D[2] and foster a 'next' generation.

The 8k+7 Conjecture predicts that these DS[1] 'volunteers' will exist but they will have step-2 webs, while S[2] has a step-8 shared web with S[1]. We will investigate here whether it is possible for the DS[1] to generate tiles from a 'normal' $3^{rd}$ generation to possibly create a D[k] sequence converging to star[2] of S[2].

**Figure 31.1** The web showing blue lines of symmetry and magenta early web

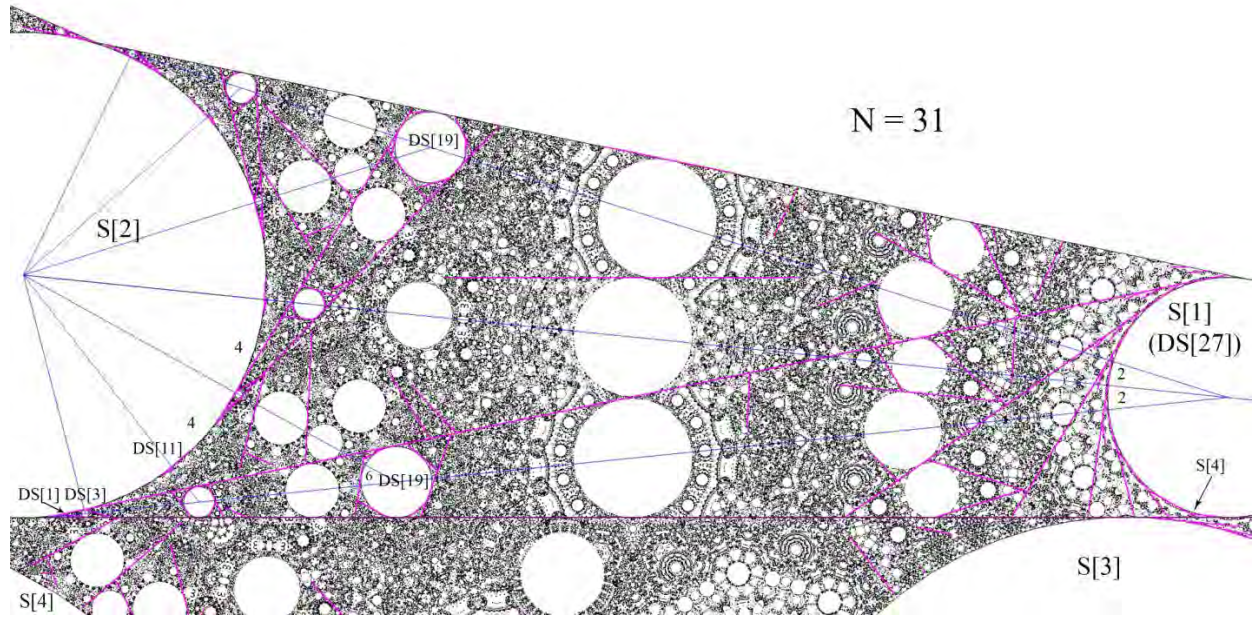

This early magenta web splits mod-8 because the 'effective' star points of S[2] are mod-8 to match the linked web evolution of S[1] and S[2]. This magenta web also makes it clear that the early step-4 web evolution of S[1] is relatively autonomous of S[2], but this web is always linked with S[3]. For N odd S[1] and S[3] will always share their star[4] points. Here the star[4] point of S[1] is 'effective' and the S[4] tile in the First Family of S[1] survives the web. However this is not universal for the 8k+7 family

**Figure 31.2** The DS[3] region at the 'foot' of S[2] showing the volunteer DS[1]

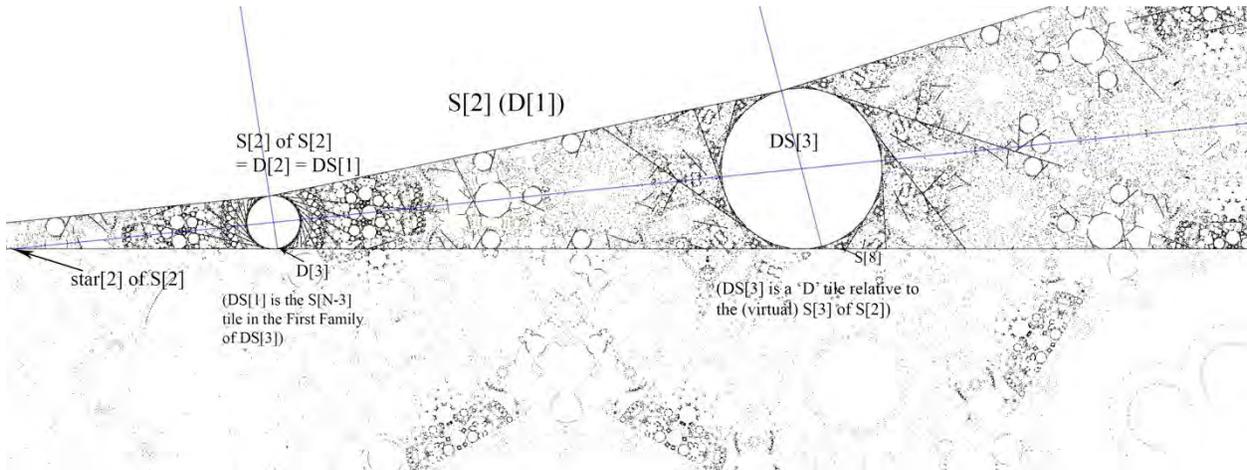

In the 8k+7 family DS[1[ and DS[3] will always have local webs which are step-2 and step-6 respectively. These webs typically give false hope for survival of First Family members and self-similarity because the initial webs show promise and the limiting webs are difficult to resolve. Our best guess is that DS[1] has no surviving First Family members but DS[1] itself always appears to be a survivor from the First Family of S[3]. It is the S[N-3] of S[3] so it is a 2N-gon like all DS[k]. The only other known survivor from the First Family of DS[3] is an S[8] as shown above. It is embedded in a prominent local web and it should exist at 6*2Pi rotations about the center of DS[3].

**Figure 31.3** Detail of DS[3]

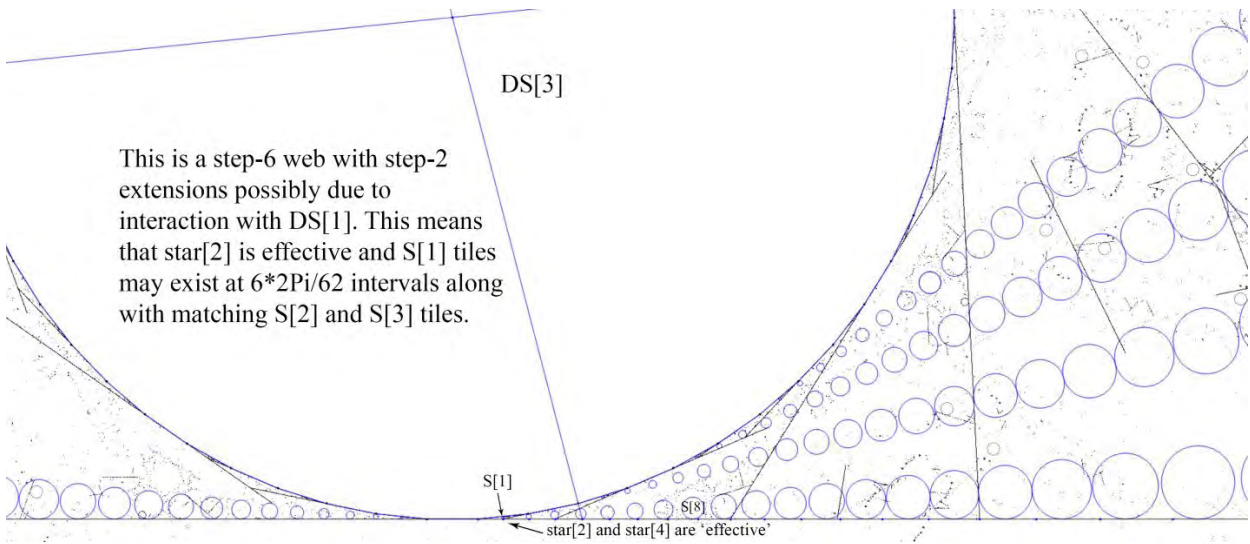

Edge extensions like those shown here are not unusual. It is possible that they may be due to the influence of DS[1] because this is really a 'shared web' in a manner similar to S[1] and S[2]. Since star[2] is now 'effective' this may facilitate the formation of S[1], S[2] and S[3] tiles on the edges of DS[3].

**Figure 31.4** The local web of DS[1] is clearly step-2 so there is a one-step offset between the left and right-side families

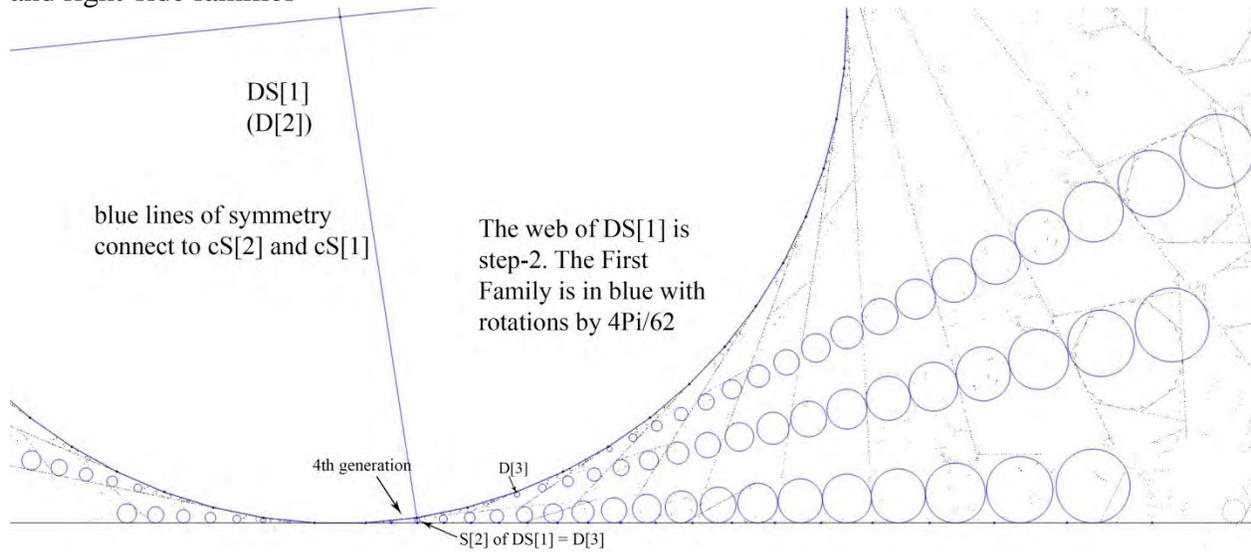

Even though this looks promising it is likely that there are no survivors.

**Figure 31.5** The web local to the N-gon S[1] is always step-4 for N odd and here S[1] supports an S[4], but this is not always true in the 8k+7 family.

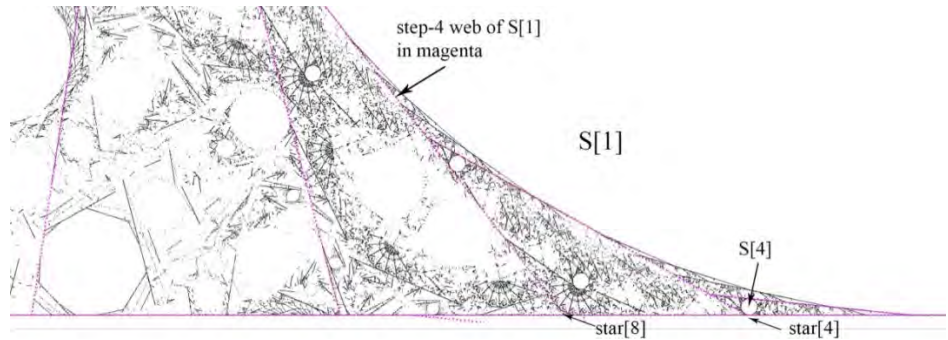

**Figure 31.6** For N-odd, S[1] and S[3] share their star[2] points while S[2] and S[4] share star[3].

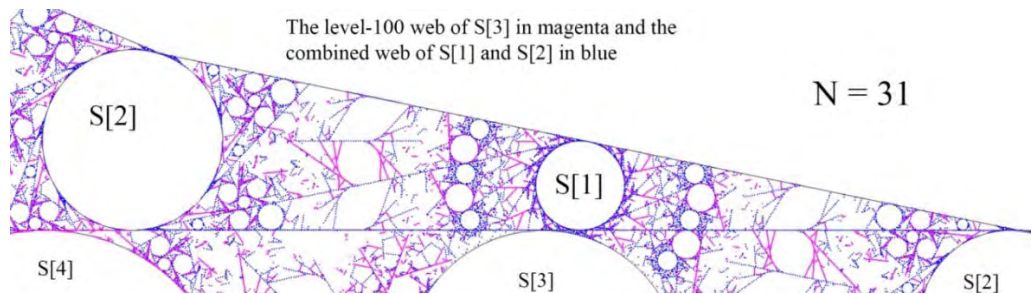

● N = 32

N = 32 is the 4th member of the 8k family. It has complexity 8 along with N = 17. As expected the First Family will have mutations in S[4] and S[8], but not M at S[14] or S[2] because the mutation criteria for the twice-even family is gcd[N/2-k,N] >2.

The twice-even family is poorly understood but we know that the local web at M may be very different from the web of N itself, even though they are both N-gons. S[14] will evolve in the $\tau$-web at step $k' = N/2-k = 2$, so the GFFT says that its right-side family will also be step-2. In the truncated version shown below this yields just S[2] and S[4] on the right and S[2],S[4] and S[6] in the left. We call this the 'short-family' of S[14].

In [H2] we discuss the Digital Filter (Df) map and show how the web local of a tile like M at S[14] can be obtained using a step-7 web. This 'maximal' step-7 Df web can be useful because it reproduces the local dynamics of the S[14] tile in a very efficient fashion. Since the 'twist' is $\rho = 7/32$, the Df map $\theta$ value is $2\pi(7/32) = 7\pi/16$ and the corresponding Df web is shown below.

**Fig 32.1** The Digital Filter map with a step-7 increment yields the following local web of S[14]

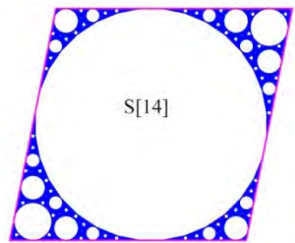

This Df web is locally identical to the $\tau$-web and since the Df map is very efficient, it is easy to probe the edge geometry of S[14] and this is done below. Since S[14] is step-2, the local family does not include an S[1] and there is no sign of First Family members for this S[2] of S[14]. Since S[14] has step-2 rotational symmetry there will be two distinct edge geometries.

**Fig.32,2.** The web local to the S[2] tile of S[14]

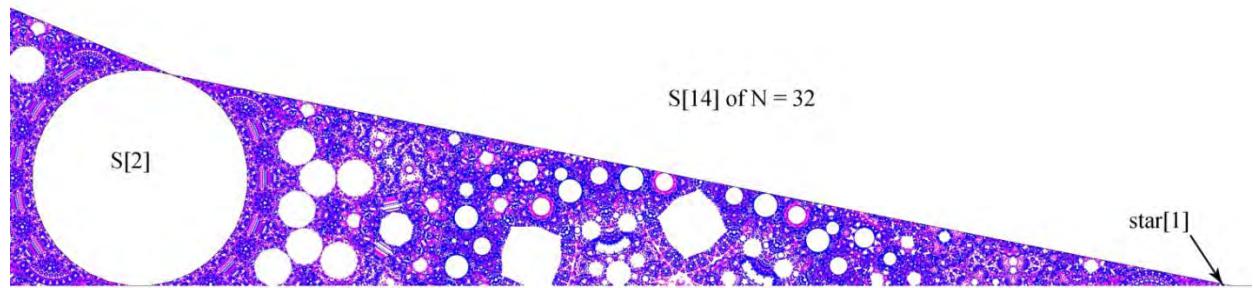

Returning to the traditional 8k family on the edges of N = 32, the Edge Conjecture predicts that there will be a DS[2], DS[6], DS[10] and S[1] at DS[14]. Here there is also a volunteer DS[4] which is displaced slightly horizontally and also mutated. The Twice-even S[1] Conjecture says that since S[1] is step-2, it may support 'step-2' families where each Skx tile is a D tile relative to S[k] in a manner similar to N-odd. For N = 24 earlier, the (rotated) S[3] tile had a real S[3] which was an S7x of S[1]. As shown below the S[3] here has S[2] tiles which are S3x of S[1].

**Figure 32.3** The blue lines of symmetry and an enlargement of the S[2] region

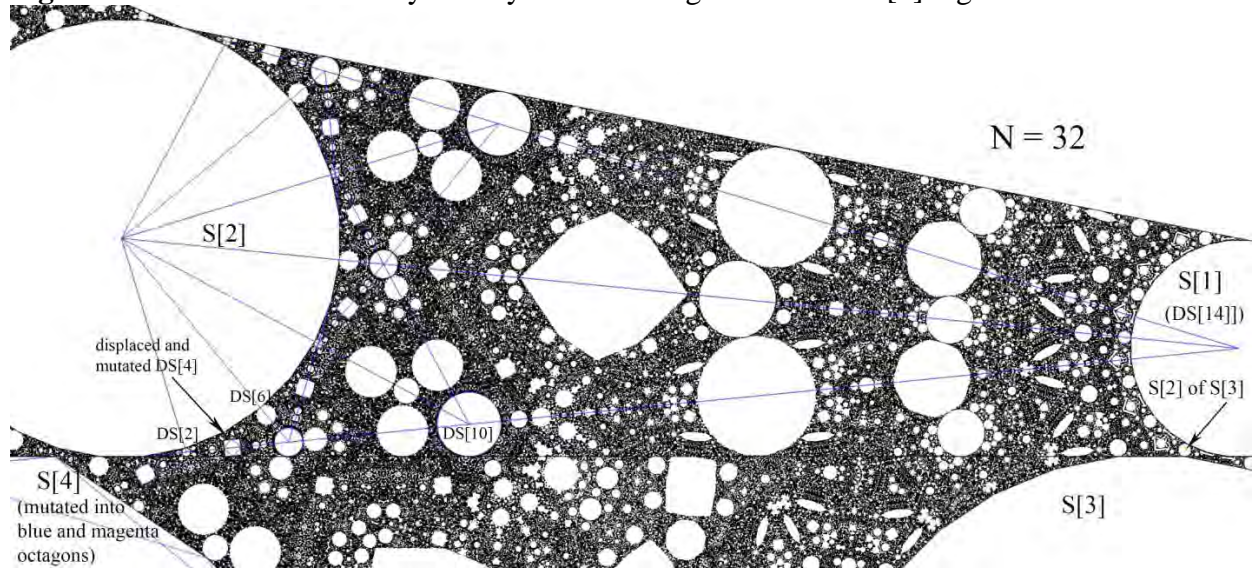

The S[4] of N is is mutated since $k' = N/2-k = 12$ and $N/\gcd[12,N) = 8$, Therefore S[4] will be the 'weave' of two regular octagons with base of length 4 spanning star[1] to right-side star[3] of the underlying S[4]. The volunteer DS[4] of S[2] should have an equivalent mutation since it has the same $k'$ relative to S[2]. However the underlying DS[4] is slightly displaced and this yields a 'mutated 'mutation which is the weave of squares with base of 8 star points based on star[1].

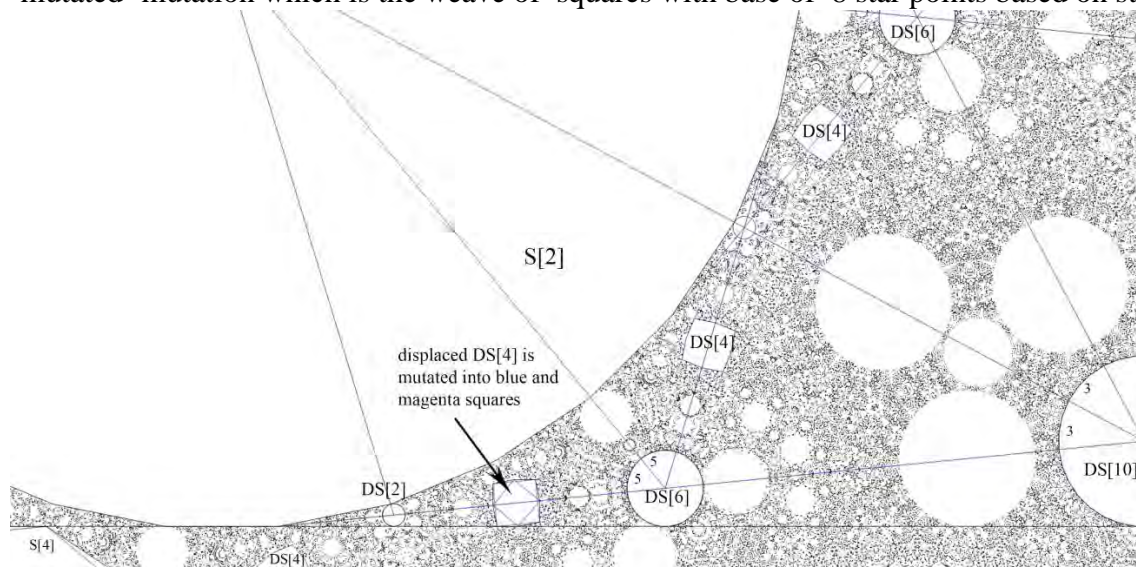

In these symmetry diagrams the vertices skipped between blue lines will increase step-4 from S[1]. Here S[1] is step-2 so DS[10] will be step-6 and DS[6] will have step-10 symmetry. Therefore the rotated DS[4] shown above is rotated $20\pi/32$ relative to cDS6. This is not step-4 relative to S[2], but the step-8 version is always valid. There is another rotated DS[4] just barely visible above and this is a rotation about the center of DS[2]. This DS[2] is step-14, but Mod[32,14] is 4 so it has a limiting step-4 web, like the volunteer DS[2] in the 8k+1 family. In general these step-4 webs are very dense and appear to have no survivors from the First Family of DS[2]. However for N-even DS[2] has a virtual D clone at star[1] of S[2] and this could support web survivors in its First Family. For N = 40 this D has an S[N/2-3] survivor which has a promising web.

**Figure 32.3** The DS[2] region

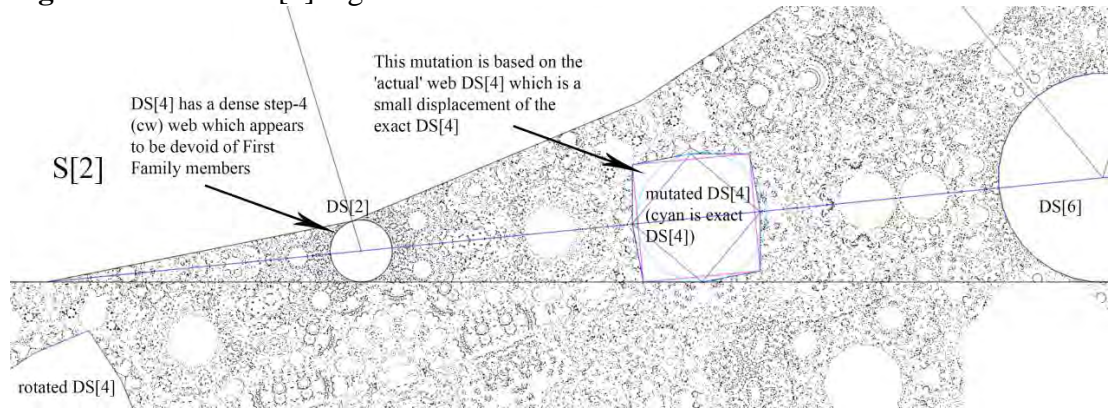

The local web of S[1] is easier to track because it is largely independent of S[2]. Typically these step-2 webs are also not fruitful in terms of First Family survivors but the Twice-even S[1] Conjecture points out that such webs can support 'step-2' families in the form of Skx tiles which are D tiles relative to the matching S[k]. As in the N-odd case this arises because the Skx are conforming to star[2] of S[1] instead of star[1] for normal S[k].. As shown in Figure 32.3 above, the S3x tile of S[1] is identical to the S[2] tile of S[3]. This appears to be the case for all members of the twice-even family – even whenS3x is virtual.

**Figure 32.5** The web local to S[1] is step-2

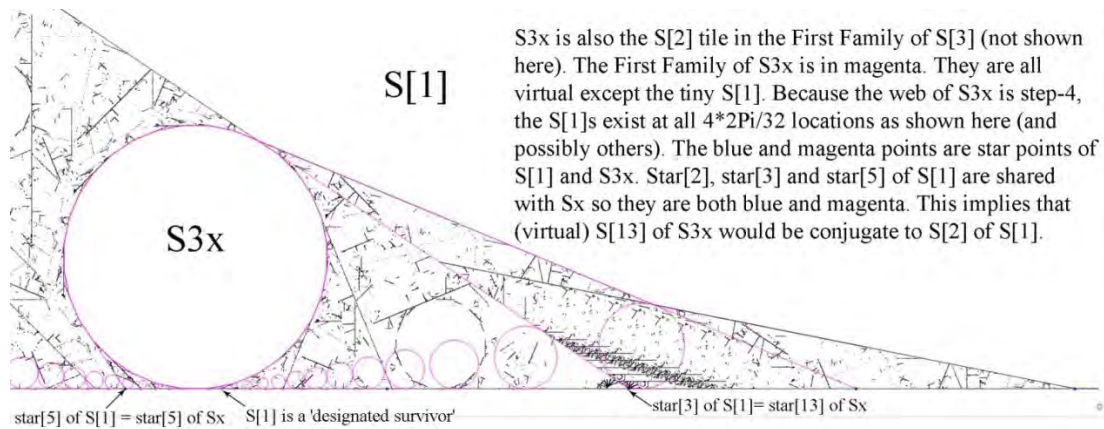

● **N = 33**

N = 33 has algebraic complexity 10 along with N = 25, 42, 50 and 66. It is the 4th member of the 8k+1 family so it has a DS[5] and a volunteer DS[2] as predicted by the 8k+1 Conjecture. N = 33 forms a pairing with N = 25 as they share algebraic complexity and edge family. This is the first member of the 8k+1 family to have a DS[21] and is no surprise that it is mutated as shown in the detail below. For N-odd S[2] is a 2N-gon and the DS[k] would be expected to have $k' = 2N/2\text{-}k$ symmetry, so $k' = 33\text{-}21 = 12$. Therefore $66/\gcd(12,66) = 66/6 = 11$. Therefore DS[21] will be mutated into the weave of two regular 11-gons to form an equilateral 22-gon.

**Figure 33.1** The combined S[1]-S[2] web showing DS[5], DS[13] and DS[21] along with DS[2\

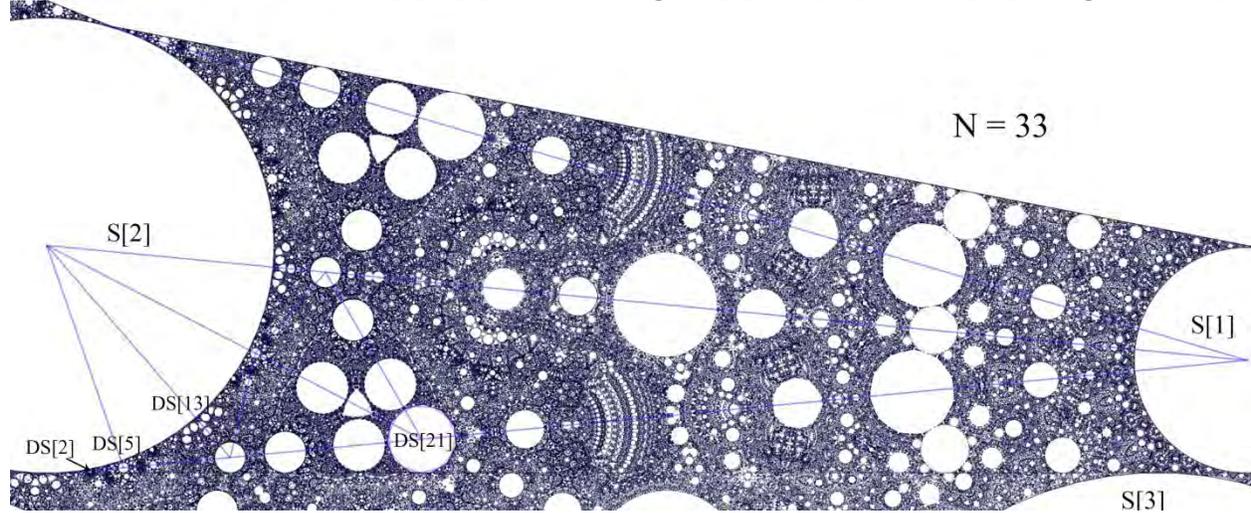

**Figure 33.2** The web local to S[2]

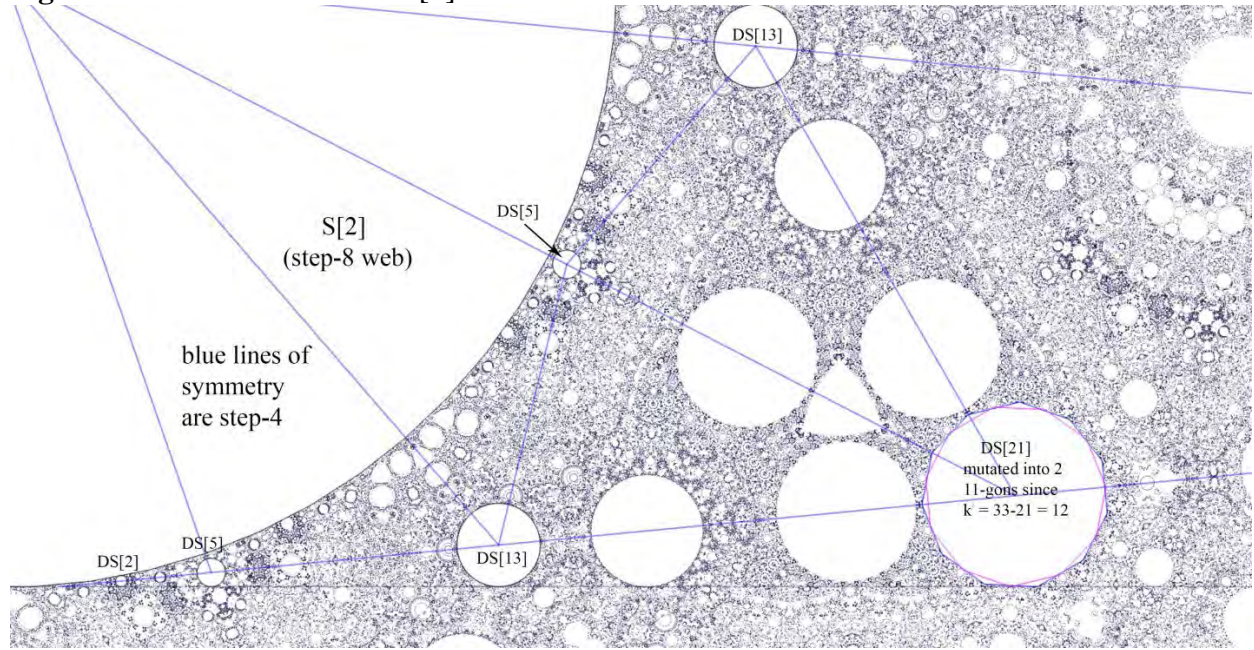

For DS[21] in the graphics above we omitted the resulting 'weave' of the blue and magenta 11-gons because, as shown below, the black weave is just a small offset of the blue 11-gon where every other blue edge has a slight magenta-induced bulge. This 'stocky hendecagon' results because of the 5 to 1 ratio of the blue and magenta star points. This is consistent with the local web having step-6 symmetry as shown by the blue lines of symmetry above.

**Figure 33.3**   The mutation of DS[21]

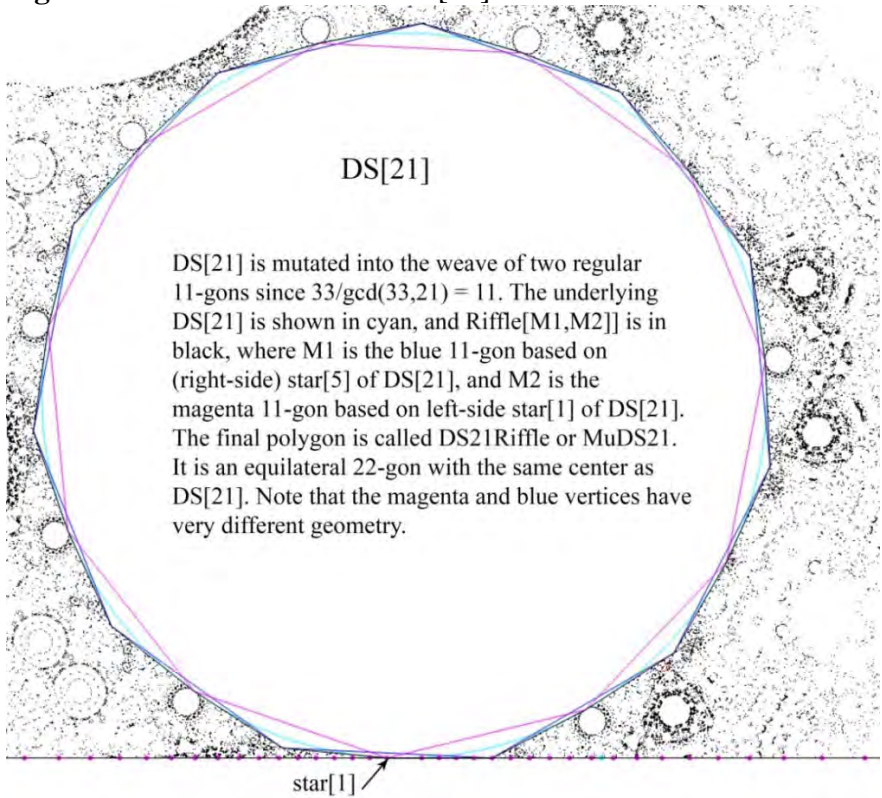

**Figure 33.5** – The region local to DS[2] and the 'wild' DS[5]

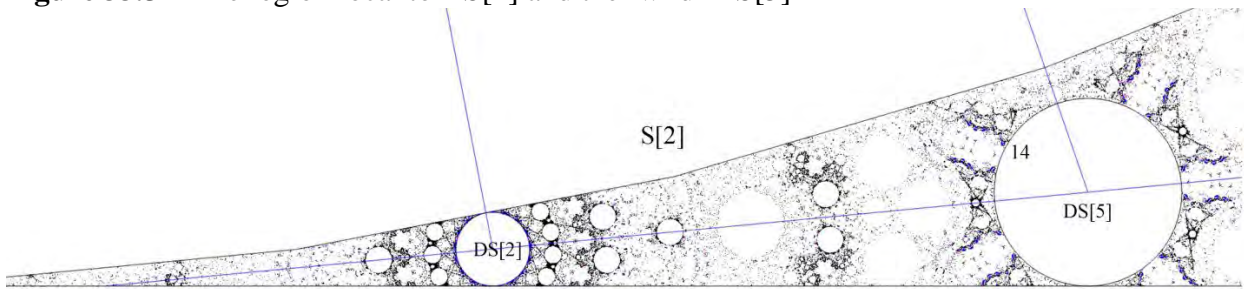

S[1] has step-4 symmetry just like DS[2] and the DS[k[ web symmetry steps increase mod 8 to yield 12, 20 and 28 steps for DS[21],DS[13] and DS[5]. That explains the 14 vertex steps shown here for DS[5]. Mod[66,28 ] = 10 so DS[5] has step-10 symmetry and (except for resonances) this will persist in the 8k+1 family so we can compare DS[5] across the family and of course

step-10 may have a very fifferent geometry depending on N. Note that DS[21] above is mutated but it still retains it's expected step-6 symmetry.

The volunteer DS[2] will always be step-4 for N odd since Mod[2N,N-2] = 4. This implies that the minimal effective star point is star[3] of DS[2] because the mod-4 countdown begins with star[N-2] of DS[2] which must match star[2] of S[2]. Since Mod[N-2,4] = 3 star[3] of DS[2] will always be effective and this creates the impression that DS[2] will support First Family tiles in the 'protected' region between star[3] and star[2]. However there is no evidence that such tiles actually exist in the 8k+1 family (or the neighboring 8k+7 family). For N = 17 the DS[2] had an S[6] that was a close match for the actual tile, but no true candidates have shown up yet.

**Figure 33.6** – Some First Family S[k] of DS[2] in magenta.

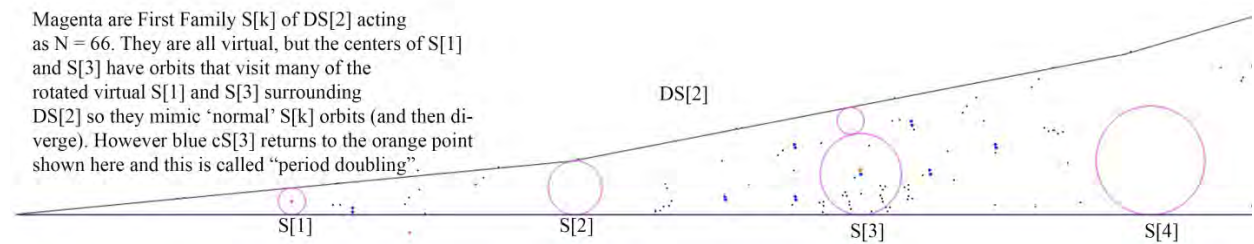

On close inspection,Fig.33.5 above shows portions of the orbit of magenta cS[1] and blue cS[3].Both of these orbits include points local to DS[2], but they also congregate in the region around DS[5] in a manner similar to N = 25. Here we do not know whether cS[1] has a periodic orbit, but the blue cS[3] has what we call a "period doubling" orbit, where blue cS[2] maps to the orange point after 9,074,553 iterations and then back to the blue in a symmetric fashion, so the period is 18,149,106 which is 33*549982.

Among regular polygons only 2N-gons have sufficient symmetry to have period doubling under $\tau$ so hidden inside S[3] is most likely a tiny regular 2N-gon with center halfway between these two points, but it is possible for a non-regular tile to also have period doubling. Here it is likely that all the even S[k] of DS[2] have doubling because DS[2] itself has doubling with 'half-life' 561 = 17*33.

The tiny black dot close to cS[3] in Fig 33.6 above is part of the web, so it should be a boundary point of some tile and it is possible that this point defines the parameters of the tile in question. As we remarked with N = 25, our interest in these tiles comes from the fact that they are like "generalized" First Family members which initially have orbits similar to S[k] satellites of DS[2], but then wander off and possibly return. As N grows it becomes very hard to resolve retailed webs, so these orbits could possibly be useful to illuminate the local geometry like candles in the dark.

Theory predicts that except for points of measure zero, any initial point will be in the interior to some tile but if that tile is small enough and the period is large enough, such a point could mimic a true web point under iteration and potentially yield information about the local web. Here with the 8k+1 family there are some natural choices of candles in the form of cS[k] of DS[2].

- **N = 34**

N = 34 has complexity 8 along with the matching N = 17.  N= 34 is in the 8k+2 family so the step-4 predicted DS[k] include DS[3], DS[7] and DS[11] along with S[1] as DS[15]. The 8k+2 Conjecture predicts that there will be sequences of D[k] and M[k] tiles converging to star[1] of S[2] with temporal scaling 18. One of the important issues to test in the 8k+2 Conjecture is whether the generations of M[k] and D[k] converging to star[1] of S[2] will be self-similar

**Figure 34.1** – The 2$^{nd}$ generation shown below will be compared with the 4$^{th}$ generation.

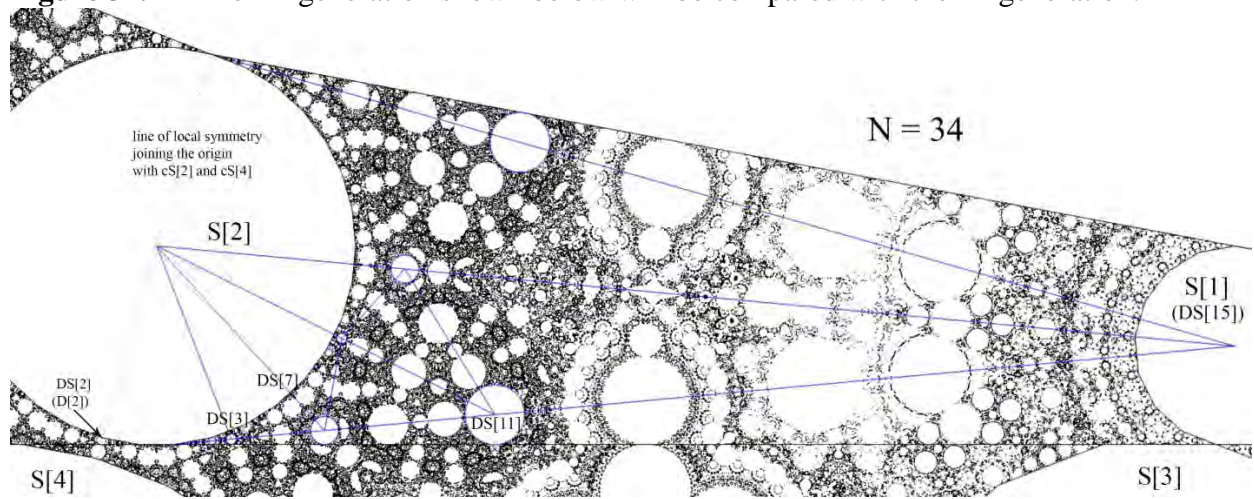

**Figure 34.2** – Enlargement of the region local to S[2] showing the primary DS[3] 'cluster' containing 2 D[2] tiles and 3 M[2] tiles

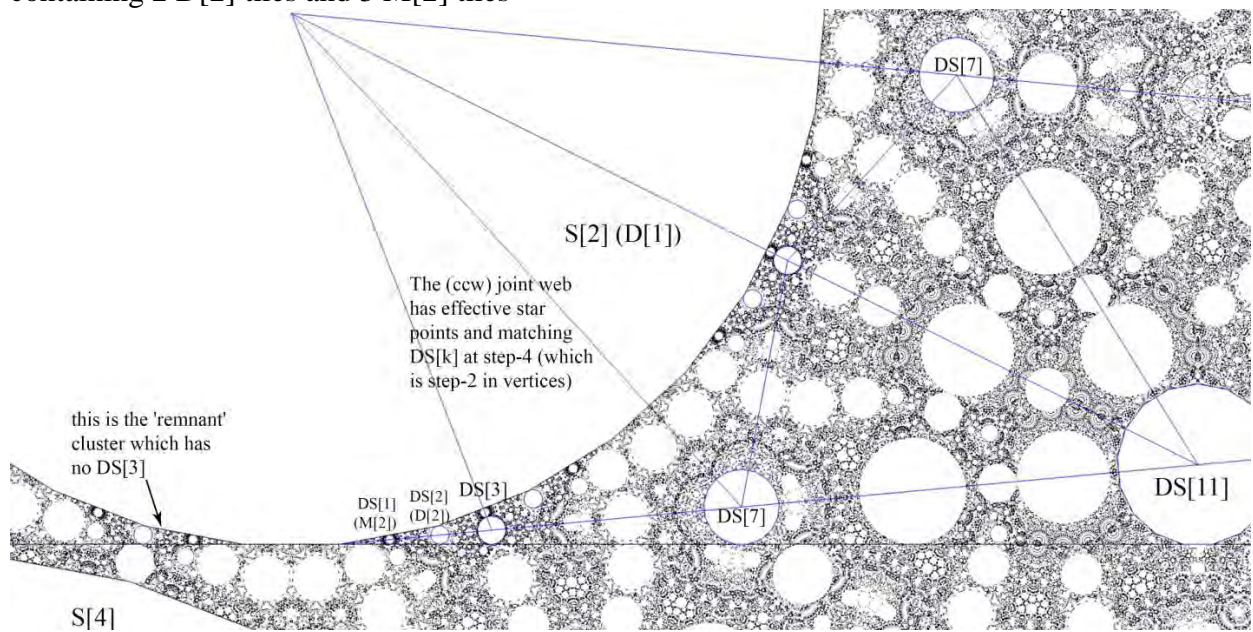

Generating detailed generations is very computationally intensive but we were able to verify that for N = 26 the first few 'even' and 'odd' generations are broadly similar but definitely not self-similar. Here we will probe the first four generations for N = 34 and see that the second generation shown above is only similar to the 4th generation in the region local to S[1]. This is consistent with N = 26 where the local volunteer tiles such as Py of S[1] were preserved between generations 2 and 4, but the region local to S[2] has volunteer Px tiles which are quite different between generations. The 8k + 2 Conjecture predicts that this will continue to be the case as N grows, so it is possible that every generation that arises will be distinct.

In [H5] and Part 1 we derived difference equations for the D[k] and M[k] that should apply to all members of the 8k+2 family. The D[k] equation is very simple and is based on the fact that 'most' of the predicted D[k] on the edges of D[k-1] come in pairs anchored by S[3]'s and the one 'outlier is an isolated D[k] with no matching S[3] as shown above. This is repeated recursively. Initially there will always be n = N/2 D[2] circling D[1] as shown above. Since the period of D[1] is always n, the next generation will have $n^2$ D[3]s along with their outliers which will be (n +1)n, so for N = 34 there will be $17^2 + 18*17$ D[3]s and the difference equation is:
$P_k = nP_{k-1} + (n+1)P_{k-2}$ where n = N/2 and $P_1 = n$ and $P_2 = n^2$. Solving this equation yields the period $P_k$ of any D[k] for N = 2n, namely $D[n,k] = -\frac{n\left((-1)^k - (1+n)^k\right)}{2+n}$

**Table**[D[17,k], {k,1,8}] ={17, 289, 5219, 93925, 1690667, 30431989, 547775819, 9859964725}

The equation for the M[k] is the same except for the initial conditions: which are pp[1] =period of M[1] = n and pp[2] =period of M[2] = n(3(n-1)/2 + 2).
The solution is $M[n,k] = \frac{n\left((-1)^{1+k} + (-1)^k n + 3(1+n)^k\right)}{2(2+n)}$

**Table**[M[17,k] = {17, 442, 7820, 140896, 2535992, 45647992, 821663720, 14789947096}

As noted with N = 26 above, it is easy to find the centers of these tiles and we have used these centers to verify most of the results shown here. Naturally the M[k] outnumber the D[k] but their ratios must have the same limit because of the shared difference equations. We also noted in N = 26 that the DS[k] appear to obey these same basic equations with just different initial conditions, so they should have the same temporal (and geometric) scaling as the D[k] and M[k] for all members of the 8k+2 family.

Here we give the steps necessary to generate the expected parameters of the 4th generation of N = 34. As with N = 26 above, we begin by generating the centers of the D[k] tiles. By convention we track the convergence at star[1] of S[2]. The First Family Theorem will yield the parameters of tiles such as S[1] and S[2]. In particular MidPoint[S2] = {cS[2][[1]],-1} so by definition StarS2=Table[MidPointS2+{hS[2]*Tan[k*Pi/34],0},{k, 1, 34/2-1}]; Then set {$x_0,y_0$} = {StarS2[[1]][[1]].-1} and cslope = slope of line from cS[1] to StarS2[[1]]. Then the center of D[k] will have coordinates {$x_k,y_k$} where $x_k = x_0 + hS[2]*GenScale^{k-1}$/cslope and $y_k = hS[2]*GenScale^{k-1}$-1; where for N twice-odd, GenScale = GenerationScale[N/2].

The S[2] tile itself will be the initial D[1] in this sequence and then D[2] will be the S[2] tile in the First Family of S[2] and this could be continued recursively since the First Families will scale by GenScale and have 'origin' cD[k] as defined above. Note that these First Families will also include the related DS[k]. However any of these tiles can also be constructed from their centers.

Here are two ways to generate the D[2] tile shown below. Since it is an S[2] of S[2], it is part of the First Family of S[2]: FFS2=TranslationTransform[cS[2]/@(FirstFamily*hS2); This takes every member of the First Family of N and scales it to match S[2]. But our default convention is a 'left-side' family, so to get right-side: FFS2R=ReflectionTransform[{1,0},cS2]/@FFS2; Then D[2] = FFS2R[[3]] (because the count starts with N itself). Now all the DS[k] are available, but an alternative is D[2] = RotateVertex[cD[2]+{rD[2],0}, 34, cD[2]]. This rotates the '3:00' vertex for N twice-even. (The radius rD[2] follows directly from the (known) height.)

**Figure 34.3** The 3$^{rd}$ generation presided over by D[2] and M[2]

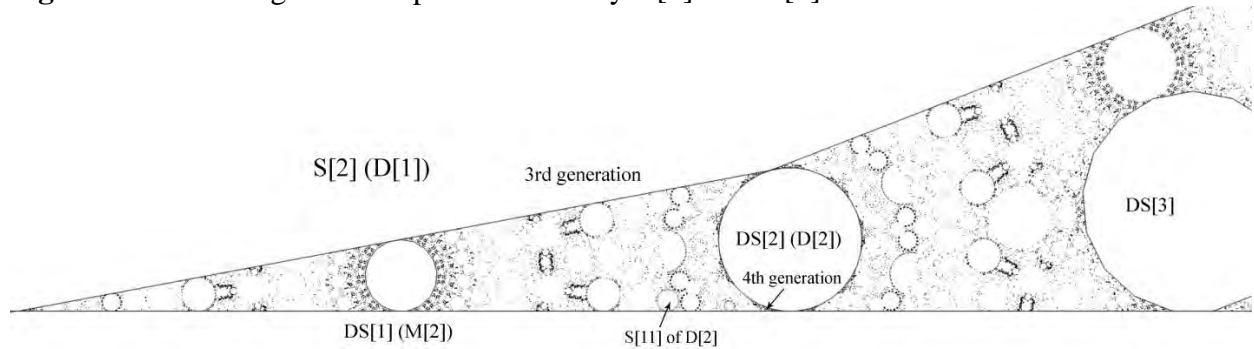

Since cD[2] is known the scaled First Family of D[2] can be imported to find the suspected step-4 survivors as shown in magenta: FFD2=TranslationTransform[cD[2]]/@(FirstFamily*hD[2]); M[2] is in this this family at DS[15] and this tile is also the S[1] of S[2]. It is relatively easy to generate the web at this scale because the period of D[2] is only 289 and this implies frequent returns under iteration, but this gets much harder with the 4$^{th}$ and 5$^{th}$ generations with periods 5219 and 93925. At a return rate of 5000:1 it possible to generate a 4th generation plot like the one below with 200,000 points within a few days, but the 5$^{th}$ generation at 100,000:1 might take more than 200 billion iterations over weeks.

**Figure 34.3** The 4$^{rd}$ generation presided over by D[3] and M[3]

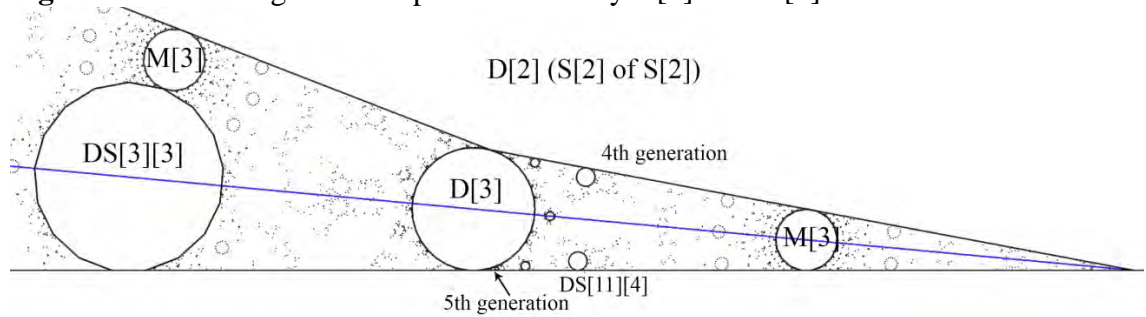

It is clear from this plot that the 4$^{th}$ generation is not self-similar to the 2$^{nd}$ generation above. They both have the expected DS[7] in a the step-4 rotated position (8π/34) on the blue center line, but here there is a large void adjacent to DS[7] that will no doubt contain a sizeable volunteer Px like those seen earlier with N = 26. This is very different from the 2$^{nd}$ generation but in the vicinity of S[1] there is a degree of self-similarity and this is consistent with N = 26 earlier.

● N = 35

N = 35 has algebraic complexity 12 along with N = 39 and 45. It is the 4th non-trivial member of the 8k+3 family preceded by N = 27, 19 and 11. It appears that this family often has volunteers between pairs of DS[k]. In the web plot below the volunteers are labeled $D_k$ but note there is no conforming volunteer between DS[7] and DS[15], and it is not clear if this is an isolated exception. The sample space is small but there are no exceptions with N = 43.

**Figure 35.1** - The $\tau$-web of N = 35 showing volunteers $D_1, D_2$ and $D_3$.

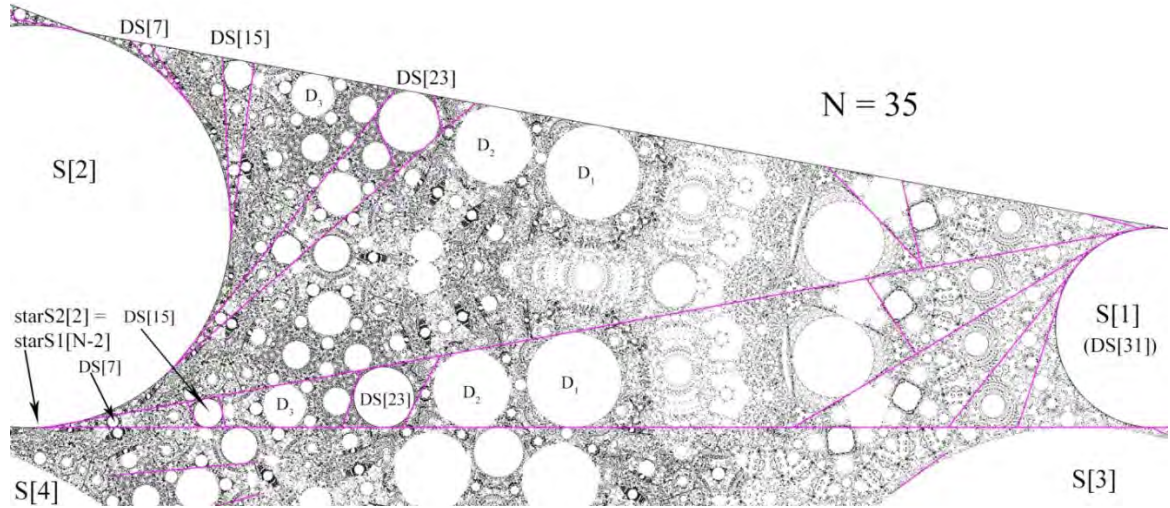

In the enlargement below it is clear that DS[7] and DS[15] (and the volunteer $D_3$) are mutated. For DS[7], our predicted $k' = N-7 = 28$ and $70/\gcd(28,70) = 70/14 = 5$. For DS[15], N-15 = 20

The underlying tiles for the DS[k] are known so they are shown in cyan. We do not know the parameters of the volunteers. Since star[2] of S[2] (a.k.a. StarS2[2]) is the penultimate star point of S[1] and an effective star point, we expect that there will be a limiting tile structure there but this region is difficult to resolve. There may be a small $D_4$ here to the left of virtual S[2] of S[2].

**Figure 35.2** - The region local to S[2]

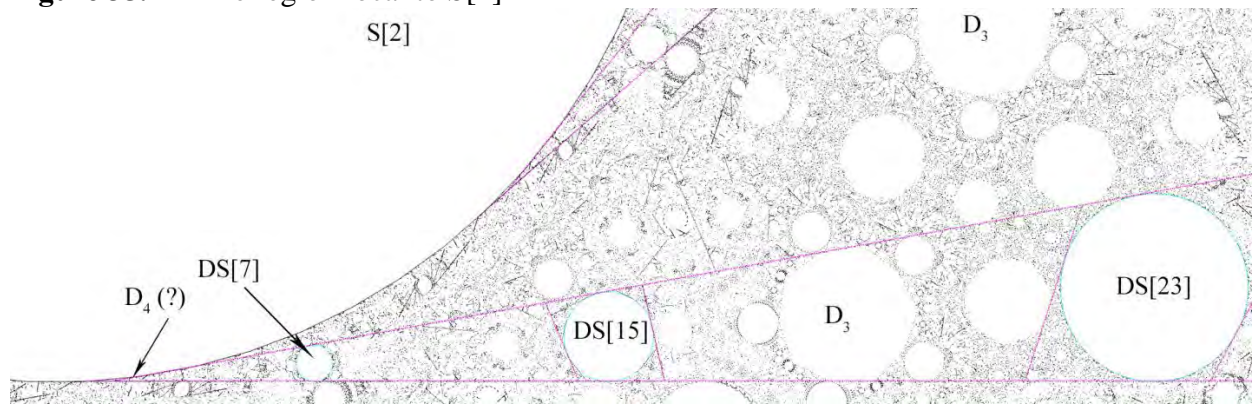

Below is the web local to S[1] showing the step-4 effective star points which include star[1]. This typically yields a step-4 tile structure that does not include an S[1] or S[2] tile of S[1].

**Figure 35.2** - The web local to S[1]

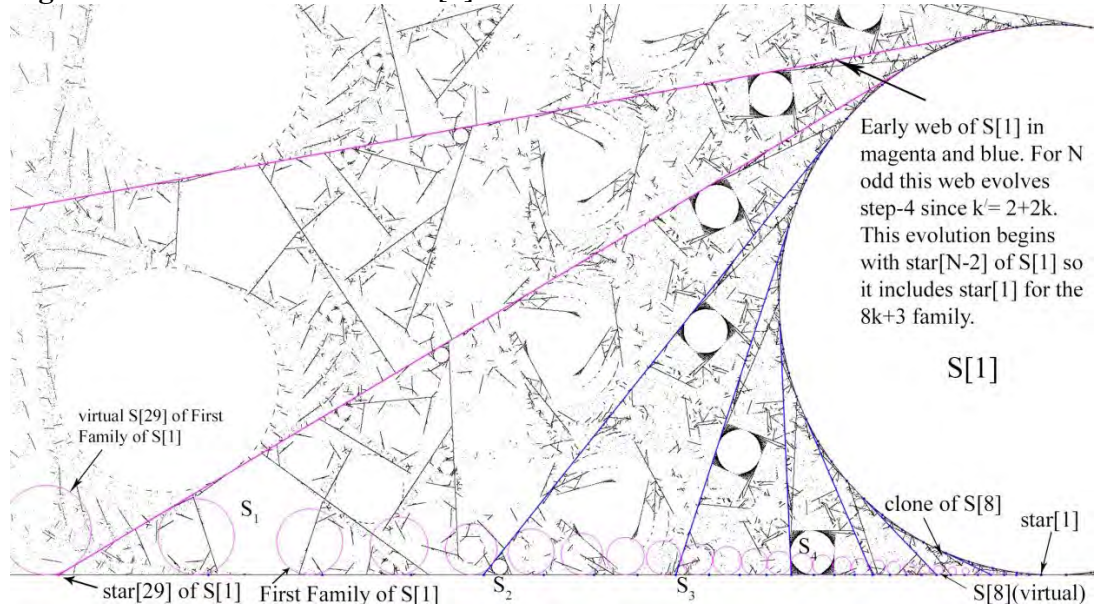

The 'volunteer' edge extensions in black are non-canonical. In general they are 2-steps removed from the ideal edges extensions in blue. Tiles Sx and Sy are volunteers and it appears that there are similar volunteers at each step-4 star point. In particular it seems that Sx has an extended family structure so we will study it below. Sy above is a mutated heptagon which should not be surprising since N is divisible by 7. Sx is surrounded by three mutated heptagons which are not visible here.

**Figure 35.3** –Detail of Sx

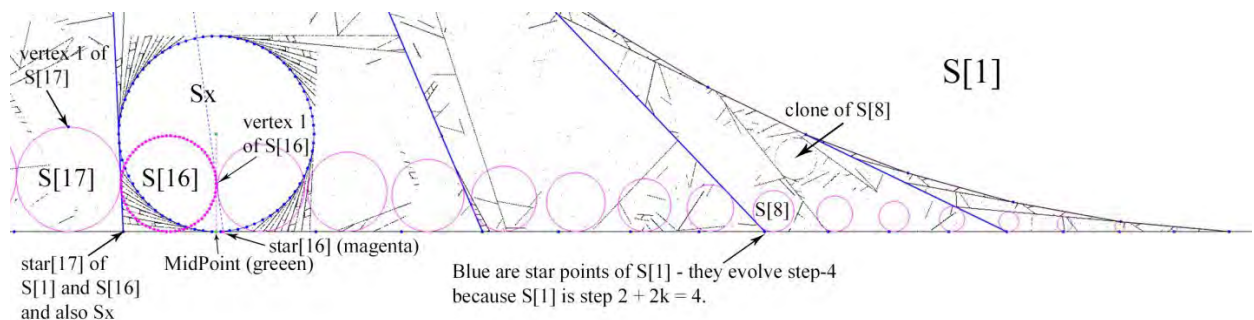

Because Sx shares star[17] with S[1] and the S[16] tile of S[1], $hSx = d/\mathrm{Tan}[17\pi/70]$ where d is the displacement of the Midpoint from star[17]. Normally another star point would be needed to find the Midpoint, but here star[1] of S[16] appears to have the same horizontal coordinate as the Midpoint of Sx.

●**N = 36**

N = 36 has complexity 6 and is the 4th member of the 8k+4 family. Therefore S[2] will be mutated in two N/4-gons since k′ = N/2-k = 16 and gcd(36,16) = 4. N = 36 is just the 2nd member of the mod-16 subfamily anchored by N = 20. It shares some geometry with N = 20, namely the volunteer Px tiles which are D tiles relative to DS[4]. This is a very small sample so we looked forward to N = 52 and discovered that this volunteer Px no longer existed. We had hopes that these Px might foster 3rd generation tiles at star[2] of S[2]. There is some evidence of this here and with N = 20, but there is little evidence of a viable 3rd generation for N = 52. The problem is that DS[4] itself has an 'awkward' small offset from S[2] that yields unpredictable local geometry and the extended First Family of DS[4] is also not fruitful. This means that we have no predications about the 3rd generation on the edges of S[2]. In general very little is known about the extended family structure of mutated tiles.

**Figure 36.1** The blue lines of symmetry and a close-up view below

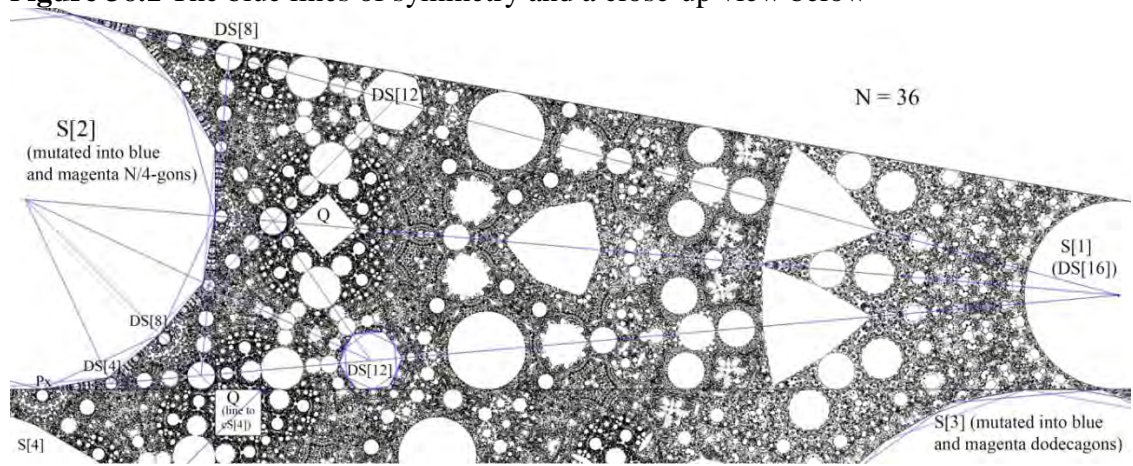

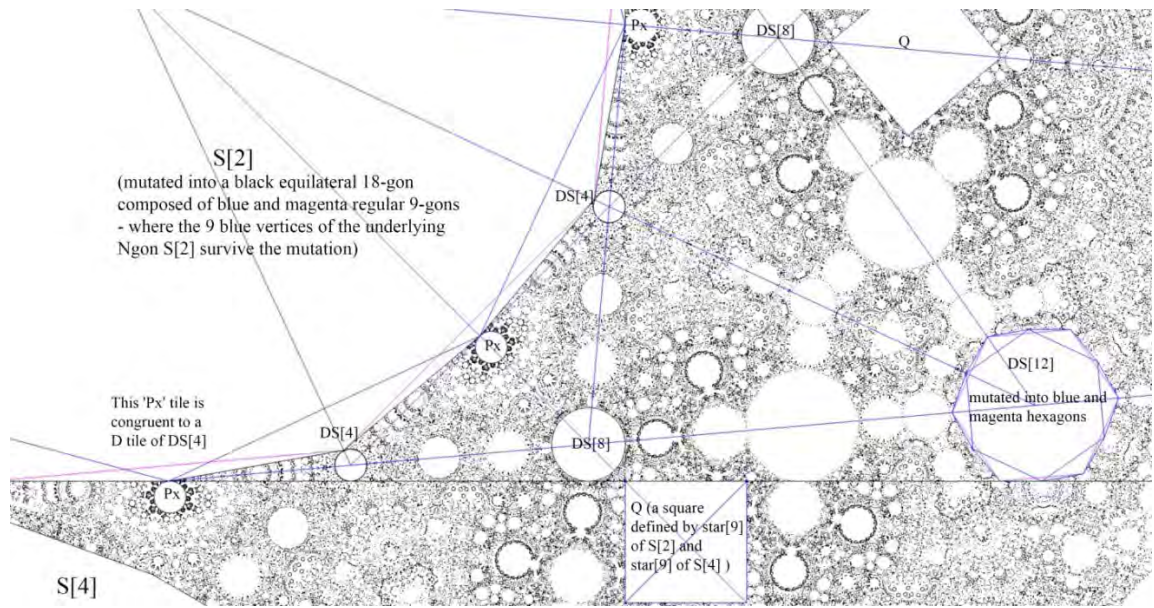

As shown above the volunteer Q tile is a square with top vertices given by star[9] of S[2] (which is also star[9] of DS[8]) and star[9] of the rotated S[4].(S[2] and the rotated S[4] always share their star[3] points). This defines 2hQ and the characteristic polynomial is more closely aligned with S[2] than with S[4] or any other known tile

**AlgebraicNumberPolynomial[ToNumberField[hQ,GenScale],x]** gives

$$-\frac{63}{256} + \frac{8515x}{256} - \frac{14541x^2}{128} + \frac{8781x^3}{128} - \frac{2615x^4}{256} + \frac{19x^5}{256}$$ compared with

$$-\frac{29}{128} + \frac{5957x}{128} - \frac{555x^2}{4} + \frac{1289x^3}{16} - \frac{1515x^4}{128} + \frac{11x^5}{128}$$ for hQ/hS[2]

DS[12] is mutated because $k' = N/2-k = 6$ and $36/\gcd(36,6) = 6$, so S[12] is the 'weave' or Riffle of two regular hexagons with radii determined by the initial star point. In this case the star point gap is also 6 and the initial star points of the underlying DS[12] are right-side star[1] and left-side star[5]. The S[12] of N should have the same mutation with the sides reversed. It is surprising that the S[12] tile of the D tile of DS[4] appears to share this same mutation. See Figure 36.4 below.

S[3] has $k' = N/2-k = 15$ so $36/\gcd(15,36) = 12$ and it has a minimal mutation consisting of blue and magenta dodecagons with 'base' length just 3, as shown above. For N = 20 the S[2] tile of S[3] was in the Step-2 family of S[1], but here there are no shared tiles and the only First Family survivor of S[1] appears to be a small S[1].

**Figure 36.2** As with N = 20, the magenta First Family of DS[4] defines the Px tiles which are D tiles (or reflected D tiles) in this First Family, so they are congruent to DS[4] but with very different dynamics because the blue vertices are vertices of the underlying S[2] while the magenta vertices are extended.

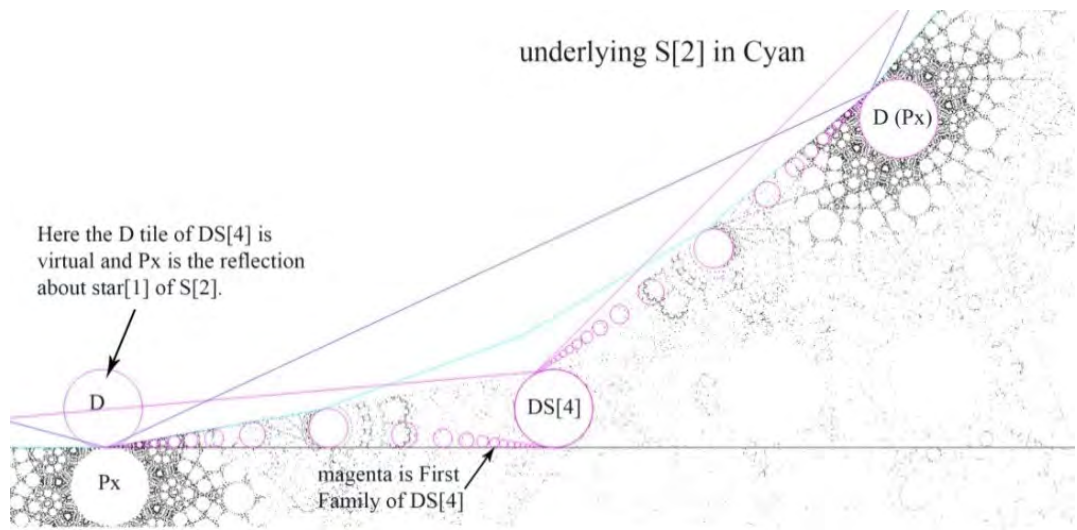

**Figure 36.3** Detail of the star[1] region of S[2] where we define the 3$^{rd}$ generation of N to be the surviving tiles in the First Family of the D tile of DS[4]. This geometry is surprisingly similar to the 2$^{nd}$ generation even though the D tile is scaled by DS[4] and not S[2].

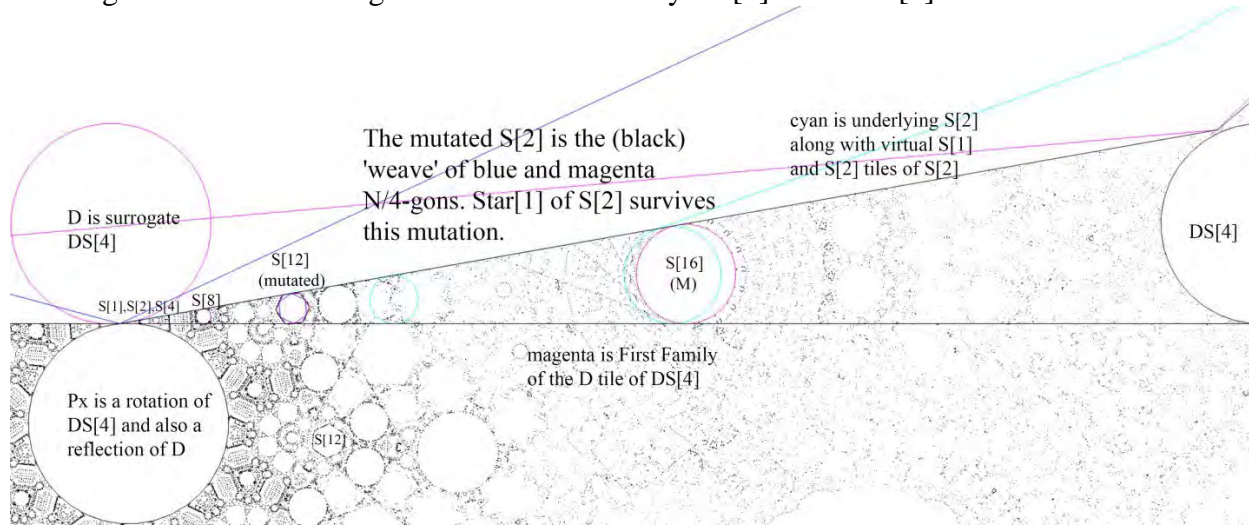

We found with N = 7 and N = 14 that the scale[2] families can and will interact with the scale[3] (GenScale) families and here we should expect this type of interaction with the 6 primitive scales of N = 36. The twice-even family has very simple scaling since GenScale is hS[1] = Tan$^2$[Pi/N] and also hS[k] = hS[1]/scale[k]. Therefore hDS[4] = hS[4]·hS[2] = GenScale$^2$/(scale[2]·scale[4]) and this is the scaling for the '3$^{rd}$ generation' here, but the traditional 3$^{rd}$ generation would be the First Family S[2] of S[2] with height GenScale$^2$/scale[2]$^2$ so these First Family members of D are about half the size of the traditional First Family of D[2] since scales[4]/scale[2]≈ .484 .

**Figure 36.4** Detail of the First Family of D

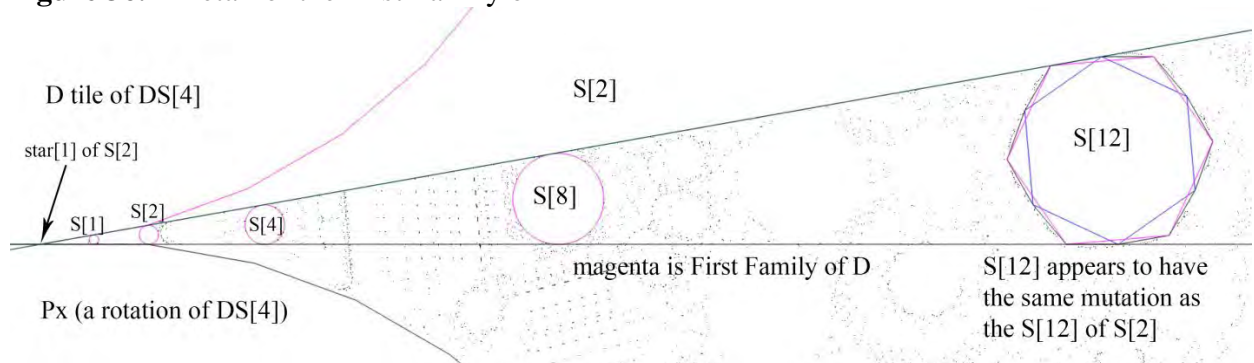

It looks like the 4$^{th}$ generation will revert to the traditional S[2] of D so that generation will be scaled by h[S[2] (of D) = GenScale$^3$/(scale[2]$^2$·scale[4]) ≈ 7.58 x 10$^{-6}$. Resolving detail at this resolution will take many billions of iterations. **Note**: Because this result is so unusual we tested N = 52 and found that the reflected D tile that we called Px was missing, but the D tile itself can always serve as a (virtual) Px 'parent' of DS[4] so the 8k+4 conjecture must be revised to make this clear. The real issue was whether the 3$^{rd}$ generation defined by D was related to the 2$^{nd}$ generation and we found no such relationship and the only verified survivor spanned S[11] and S[12] of the First Family of D. Therefore we have no clue about the general nature of this so-called 3$^{rd}$ generation in the 8k+4 family.

●**N = 37**

N = 37 has complexity18 and is a member of the 8k+5 family, so there will be a DS[1] at the end of the Rule of 8 chain which begins with S[1] at S[N-4] = S[33]. The S[1] and S[2] tiles here are 74-gons which have their counterparts at D, which is 8k + 2. Since D has a relatively 'well-behaved' edge structure, the issue is how this influences the family structure of S[1] and S[2] which are displaced members of this family. In particular we are interested in the local geometry of DS[1] . In the modified 'star[2] family' of S[2] this DS[1] is the D tile relative to the virtual S[1] of S[2], so it is a valid D[2] and the only tile that can exist in both the 'normal' and 'star[2]' families of S[2].

**Figure 37.1**  The blue lines of symmetry

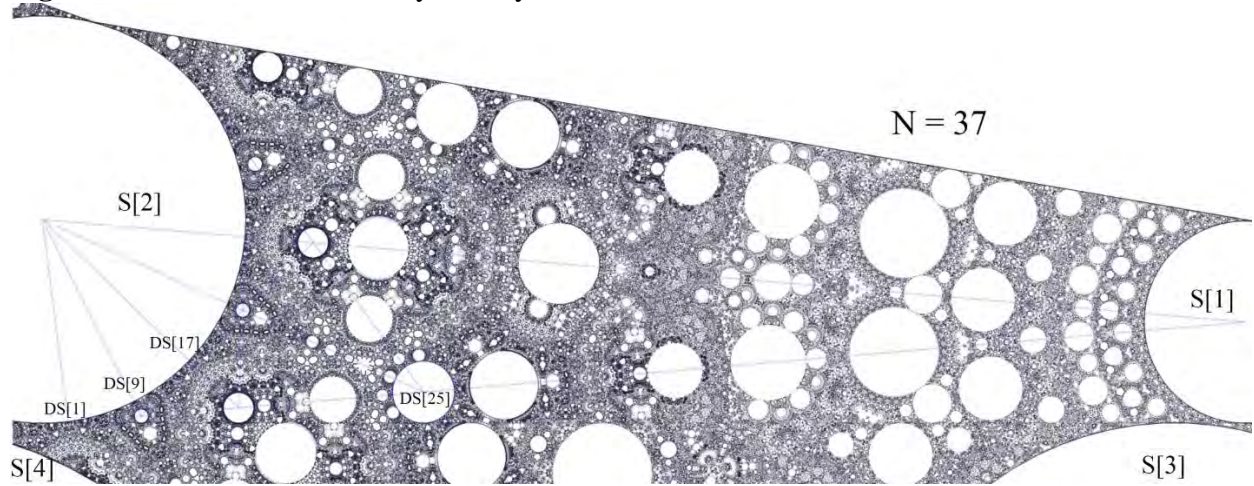

**Figure 37.2**  Symmetry local to S[2]

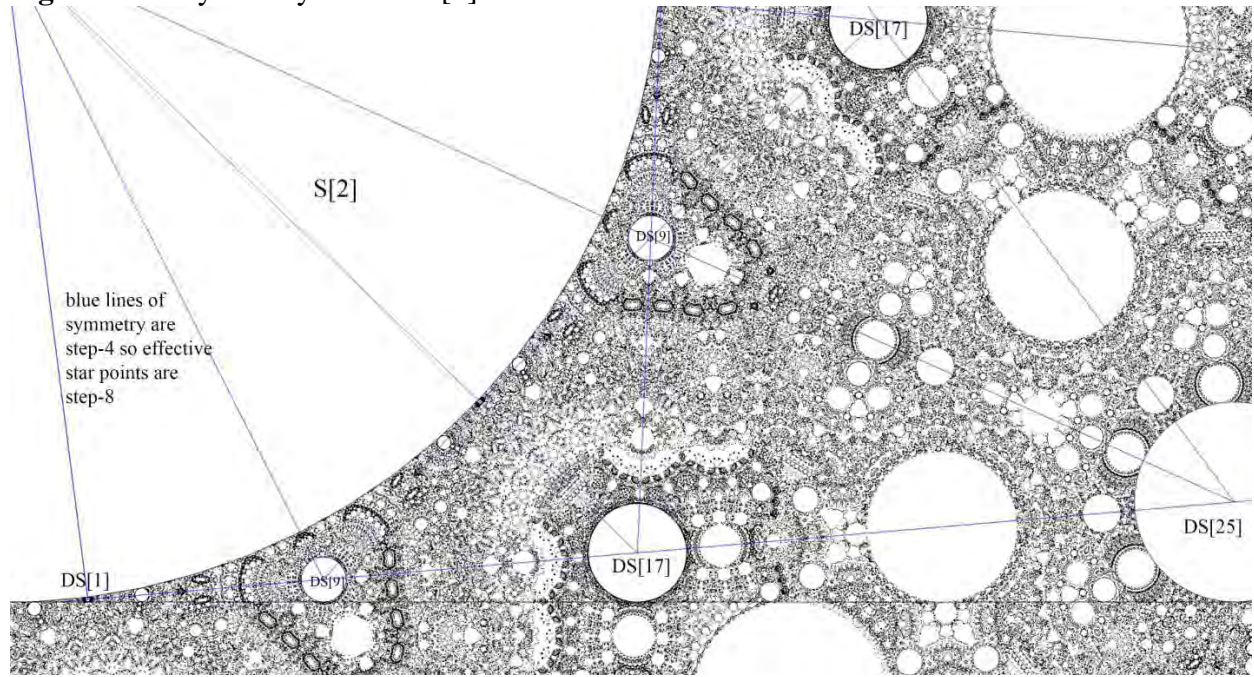

DS9] and DS17] have multiple layers tiles as shown by the dark regions on their edges. Below is an enlargement of the edges of DS[17] where an S[1] or S[2] may exist at the step-2 blue vertices

**Figure 37.3** The edges of DS[17] which is in the N = 74 family

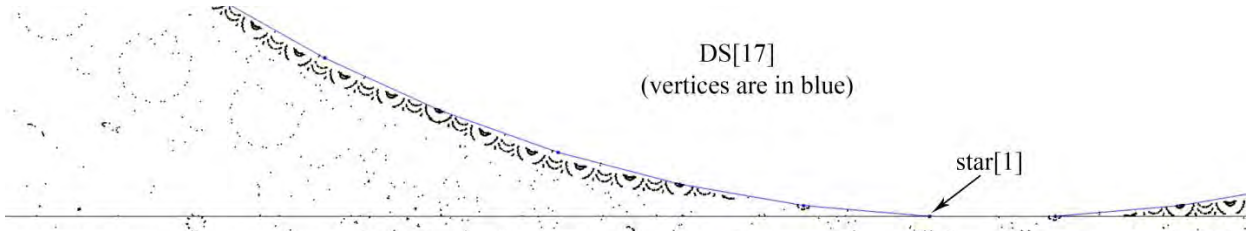

**Figure 37.3** The DS[1] tile of S[2] is also the S[2] of S[2] so I could act as a D[2]. But DS[1] always has a step-2 web in the 8k+5 family and these webs are not well understood.

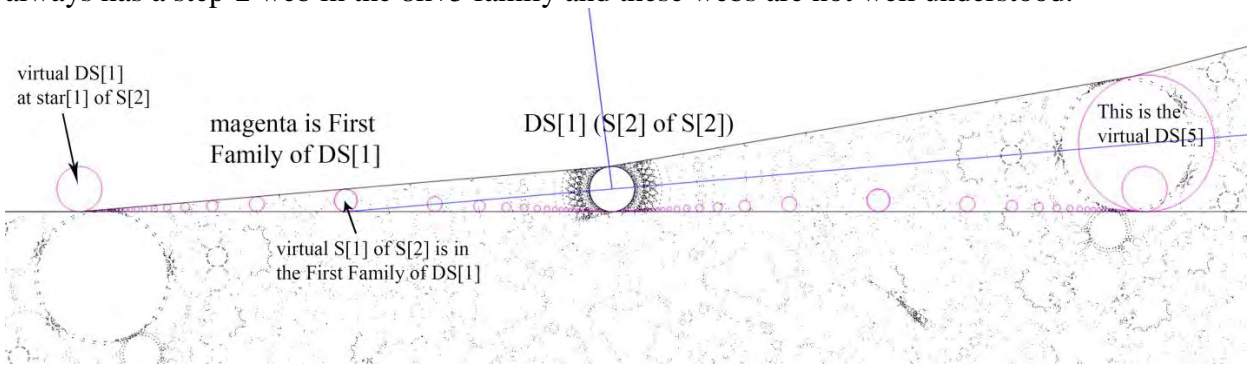

**Figure 37.4** The First Family shows some possible survivors in green. These images took billions of iterations but there is still not sufficient detail. The small S[2] shown here has period 716098 but its neighbors have no obvious period.

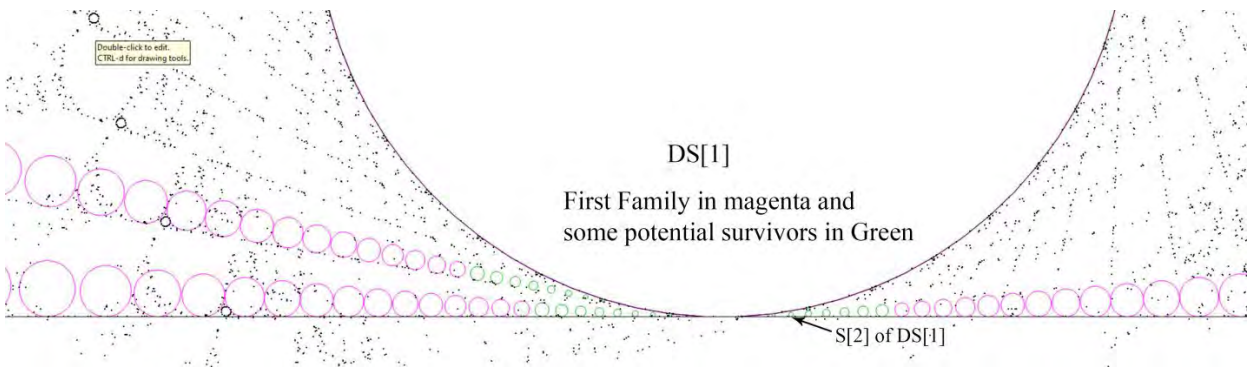

●N = 38

N = 38 has complexity 9 along with the matching N = 19. It is in the 8k+6 family along with N = 14, 22 and 30, so it has a DS[1] which has potential for extended family structure. This DS[1] is the S[1] tile of S[2] so it is an N/2-gon just like S[1] – and they both have step-1 webs.

The 8k+6 family interfaces with two embedded odd families: 8k +3 and 8k+7 and these are very different families, so the natural question is whether the 8k+6 family has an underlying mod-16 structure to distinguish these cases. We will explore this possibility here by comparing N = 38 with the embedded N = 19 and N = 46 with embedded N = 23. These two embedded N/2-gons have very different dynamics and the question is how this affects the geometry local to DS[1]. N = 19 has a minimal second generation structure with just a DS[7] while N = 23 is in the 8k+7 family with a DS[3] and volunteer DS[1]s with potential extended tile structure at S[2].

**Figure 38.1** The S[1]-S[2] combined web

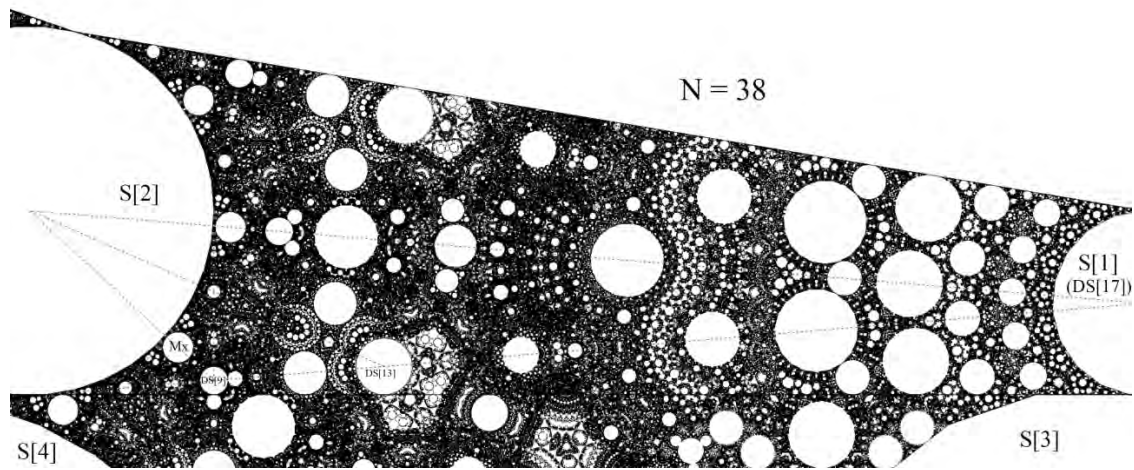

**Figure 38.2** The web local to S[2] showing a volunteer Mx which is similar to the Dx volunteer of N = 30.

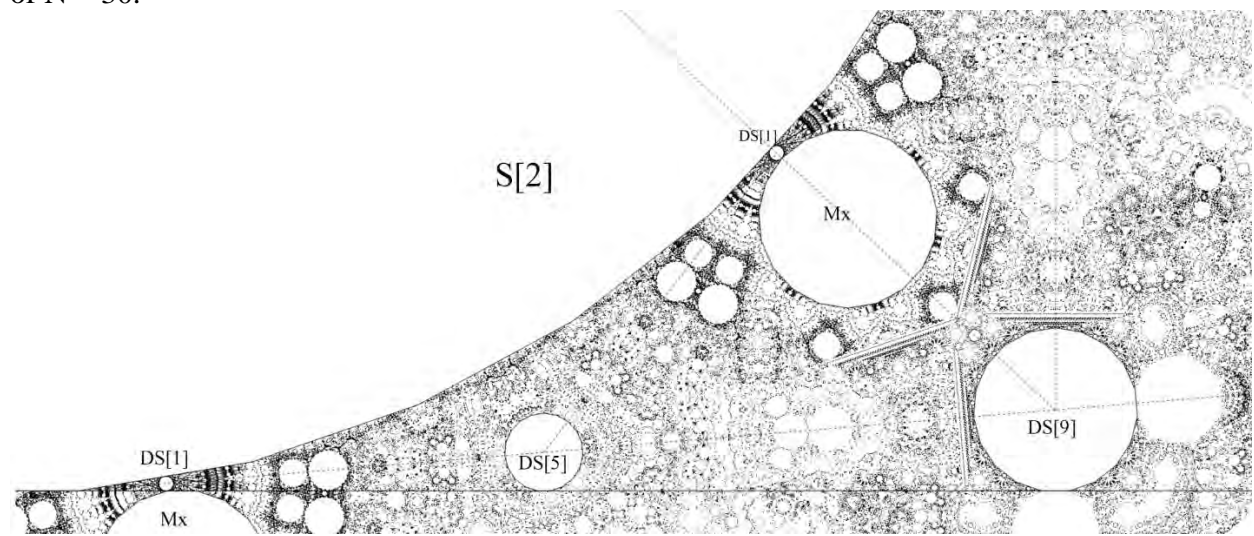

The Mx at the foot of S[2] shares the star[2] point of S[2] with the right-side star[1] of DS[1] but DS[1] is not in The First Family of Mx and their relationship seems to be coincidental except that the height of Mx appears to be hDS[9] + hDS[1]. This is the first known example of such a relationship . If it is valid there should be other examples which will appear.

**Figure 38.3** The DS[5] region

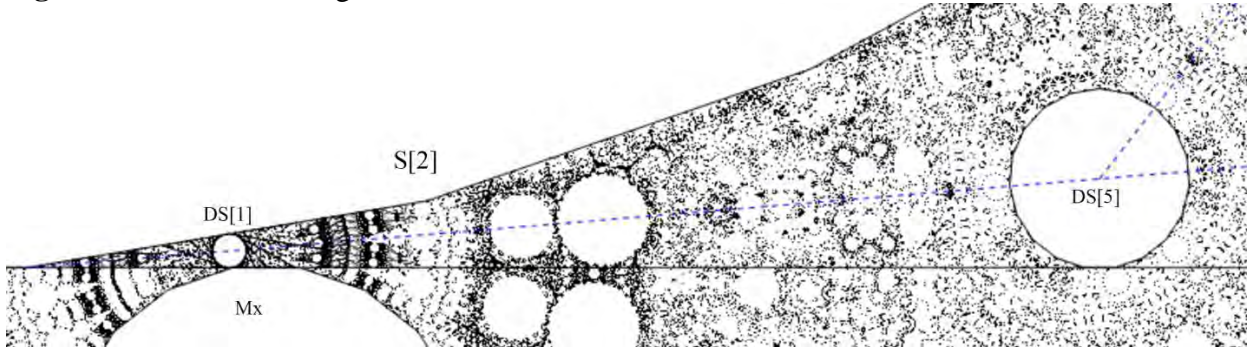

**Figure 38.4**   DS[1] is the S[1] tile of S[2] so it is an N/2-gon with a step-1 web in a manner similar to the S[1] of N.

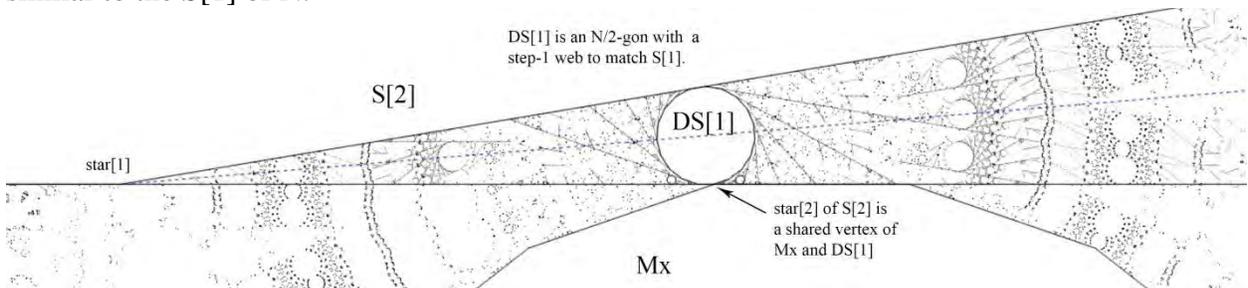

These step-1 webs initially look promising, but we have found no evidence of First Family survivors here or in subsquent 8k+6 cases. See N = 46 for a more detailed analysis of the issues. There we note that the virtual D tile on the First Family of DS[1] shares a vertex with star[1] of S[2] and this is a good place to look for structure since it is the S[2] of S[2] also known as D[2].

**Figure 38.5**   The S[13] tile has interesting geometry with potential for next-generation tiles.

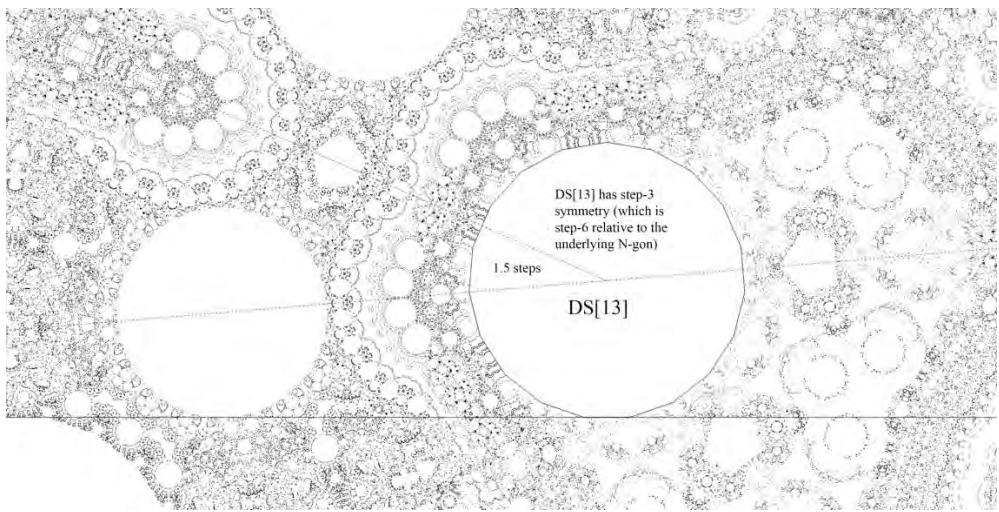

● **N = 39**

N = 39 has complexity 12 and is in the 8k+7 family so there will be a DS[3] and the 8k+7 Conjecture says that the predicted DS[3] will generate right and left-side DS[1]s.

**Figure 39.1** The symmetry of the combined S[1]-S[2] web\

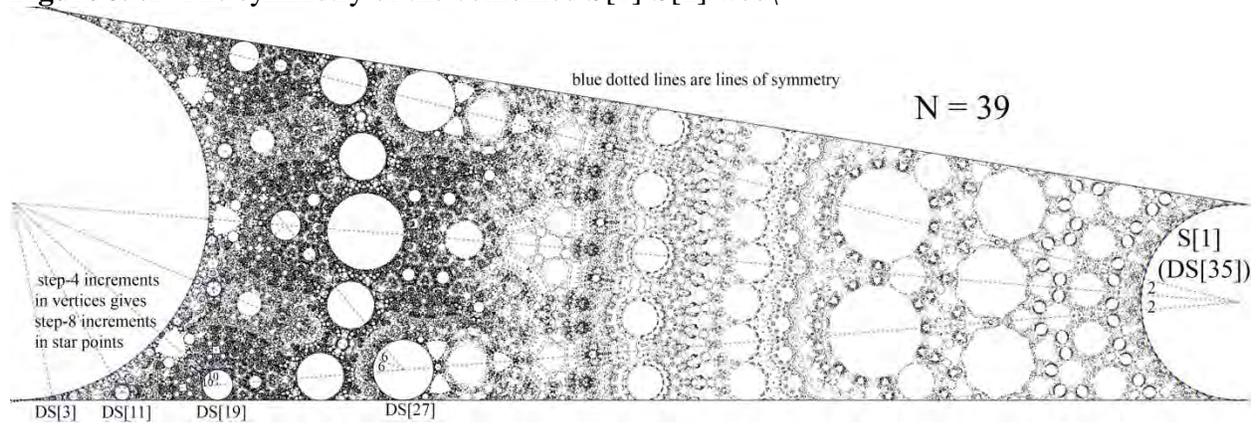

**Figure 39.2** Detail of the S[2] region. DS[3] showing the k′ steps increase mod 4.

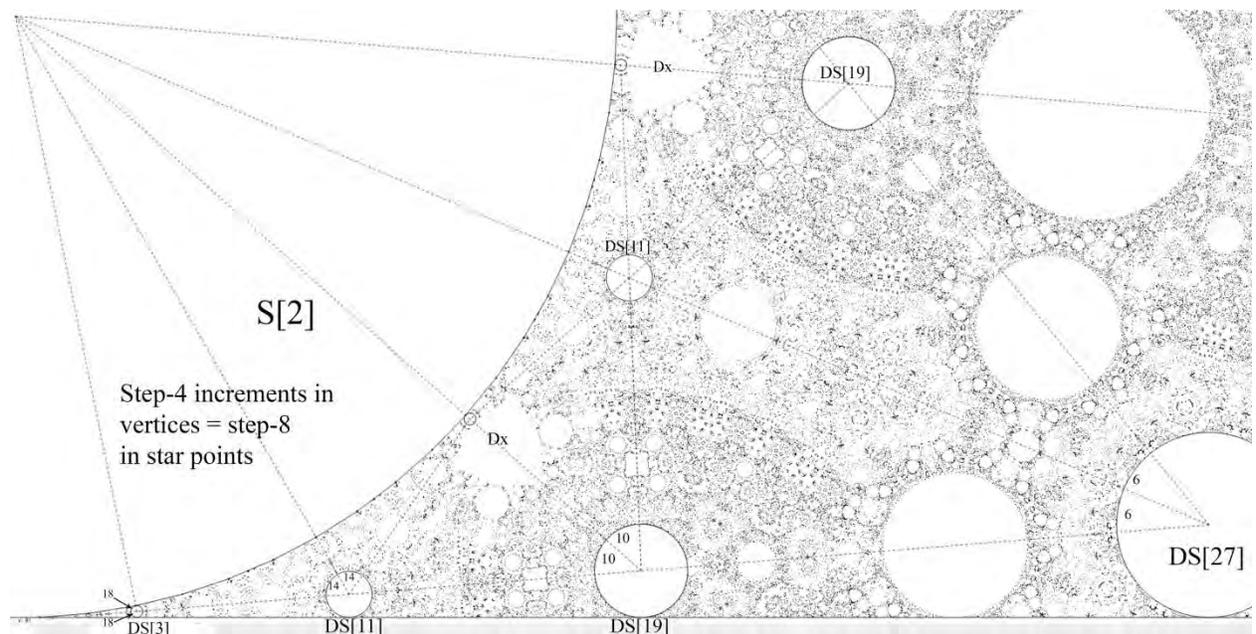

For N odd, S[2] is a 2N-gon so the DS[k] are also 2N-gons which evolve in the web with $k' = 2N/2 - k$ steps, DS[3] and DS[27] will be mutated since they have $k' = 36$ and 12 respectively so both have $2N/\gcd(k',2N) = 78/6 = 13$. This will imply that these two have identical mutations as the weave of two regular 13-gons as shown in Figure 39.3 below. (Unlike the N even case the S[3] of N does not share step sequences with DS[3] because $k' = N-2k$ for N odd, so the S[3] of N = 39 has $k' = 33$ with $78/\gcd(33,78) = 26$ which is barely noticeable.)

**Figure 39.4** DS[3] and DS[27] are shown here not to scale. They have the same mutation consisting of the weave of a blue regular 13-gon based on star[5] of the underlying DS[k] (in cyan) and a magenta regular 13-gon based on opposite-side star[1] of the underlying DS[k].

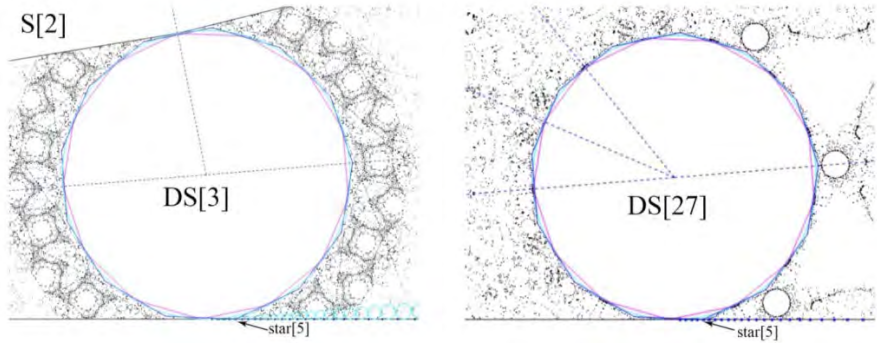

The 'base' of these mutations must span 6 star points since gcd(N,k′) = 6 and the Mutation Conjecture predicts that the minimum star point should be minimum of N-2-jk′ which here is 1 for both DS[3] and DS[27]. This means that both mutations will be identical and span star[1] to star[5] as shown above. These graphics do not reveal the interesting dynamics that occur below the horizontal 'base' line. Actually DS[3] shares it's star[5] point with a large Dx tile which seems to arise in the early web due to the mutation of DS[3].

**Figure 39.5** The DS[3] region showing the likely origin of the large Dx tiles as an extension of the star[5] edge of DS[3] shown in blue below.

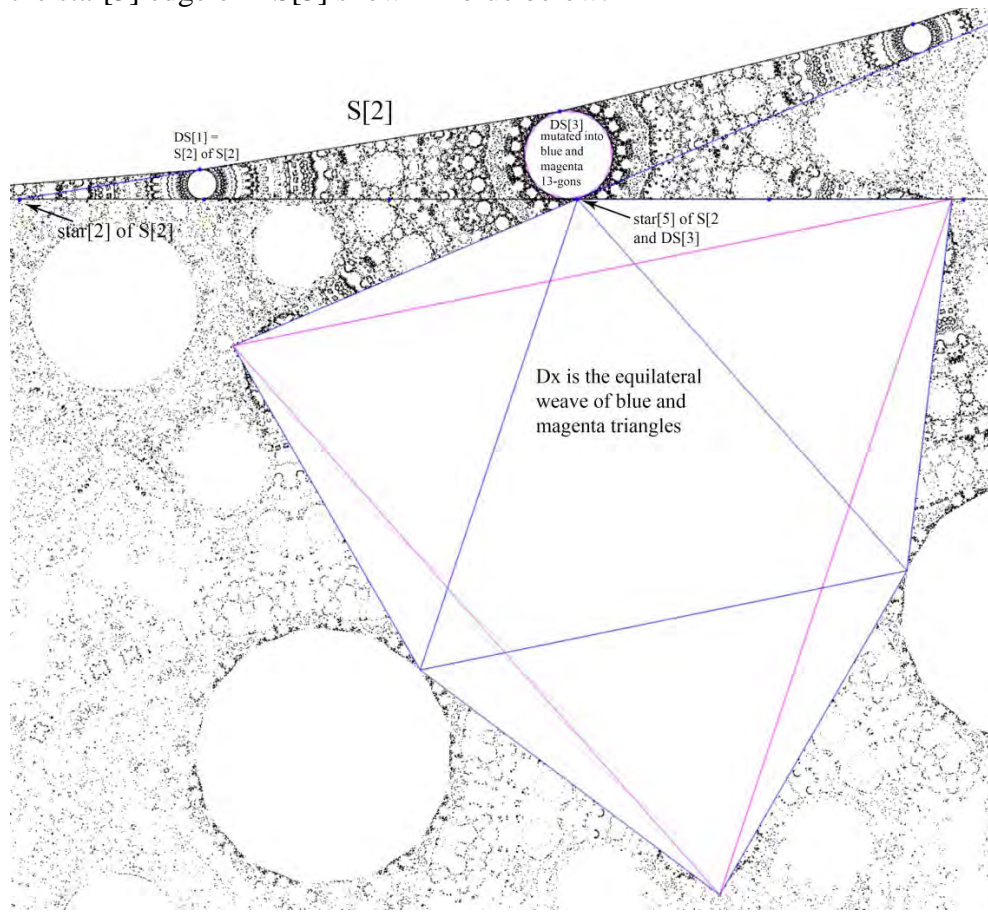

A clone of Dx occurs in N = 30 but it is based on the DS[1] – which for N-even is simply the S[1] tile of S[2]. In both cases the corresponding DS[3] tile was mutated . Here the only parameter needed to define Dx is edge length which is the offset $s_1$ between star[5] and star[3] of S[2] and this is also the offset of DS[1] and DS[3] as can be seen above. Since Dx is equilateral the 'top' edge is all that is needed if we assume correctly that it is step-6 just like DS[3]. Then working backwards we find that it is the weave or Riffle of the blue and magenta triangles. This of course means it is 'canonical' since $s_1/s_0$ will be in the scaling field $S_{39}$ generated by GenScale[39] where $s_0$ is the side of N.

**AlgebraicNumberPolynomial[ToNumberField[s1/s0,GenScale],x]** yields

$$\frac{1871}{26624} - \frac{545859x}{26624} + \frac{7923081x^2}{26624} - \frac{1868429x^3}{26624} - \frac{6409373x^4}{13312} - \frac{2266215x^5}{13312}$$

$$+ \frac{1182769x^6}{13312} + \frac{660331x^7}{13312} + \frac{4411x^8}{26624} - \frac{70847x^9}{26624} - \frac{6459x^{10}}{26624} + \frac{391x^{11}}{26624}$$

We would expect this equation to simplify by replacing $s_0$ with the side of DS[3] and indeed that is true as shown below for $s_2/s_1$ where $s_2$ is the side of DS[3]

$$\frac{11}{2048} + \frac{2567x}{2048} + \frac{20053x^2}{2048} - \frac{168639x^3}{2048} - \frac{3273x^4}{1024} + \frac{124931x^5}{1024} + \frac{86949x^6}{1024} + \frac{14353x^7}{1024} - \frac{8041x^8}{2048} - \frac{2941x^9}{2048} - \frac{175x^{10}}{2048} + \frac{13x^{11}}{2048}$$

For comparison we also tried a heuristic of iterating the nearby star[7] of S[2] under the Dc map to a depth of 400 million. As shown above, this star point is adjacent to a small tile that may share a local star point with Dx. If the star[7] point is indeed a vertex of this tile then iteration under Dc will eventually generate further vertices in the same way that vertices of N will generate N because the rotation angle w = $2\pi/N$ is fixed under Dc. So all the tiles generated in this fashion will be N-gons, but they share star points with any matching 2N-gon. See N = 47 below for multiple examples.

**Figure 39.6** Detail of DS[1] and DS[3]

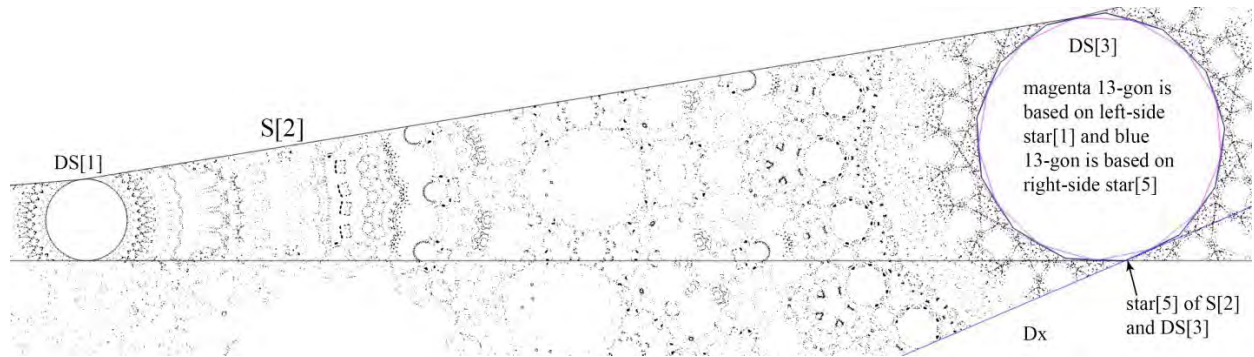

It is no surprise that this mutated DS[3] supports no First Family tiles, but it is still obvious that it evolved step-6 even after the mutations. DS[1] also appears to support no First Family tiles and this is consistent with other members of this family.

Note: This is just the third case done with the new 13[th] generation I9 computer running Mathematica 14 and the results have been promising with processing speed 4 times a last-generation I7. The detail above would have been very difficult to achieve without this new HP powerhouse and the improved memory management of Mathematica 14.

●N = 40

N= 40 has complexity 8 along with 32 and 48 – all of which are in the 8k family. This family has a DS[2] which can serve as a traditional D[2], but there are only matching M[k] for N = 16 and this seems to be sufficient to guarantees future generations. For N = 32 and N = 40 there is again no sign of an M[2] with little hope of an M[3]. However the DS[2] tile of N = 40 has a step-4 web and generates an S[17] which is just one step away from M[2] at S[18]. But this is really an S[17] in the First Family of D – which is virtual here. It is not unusual for virtual D ties to generate extended families. This occurs with DS[4] for the N = 20 subfamily of N =8k+4.

**Figure 40.1** The blue lines linking centers and star points are lines of symmetry

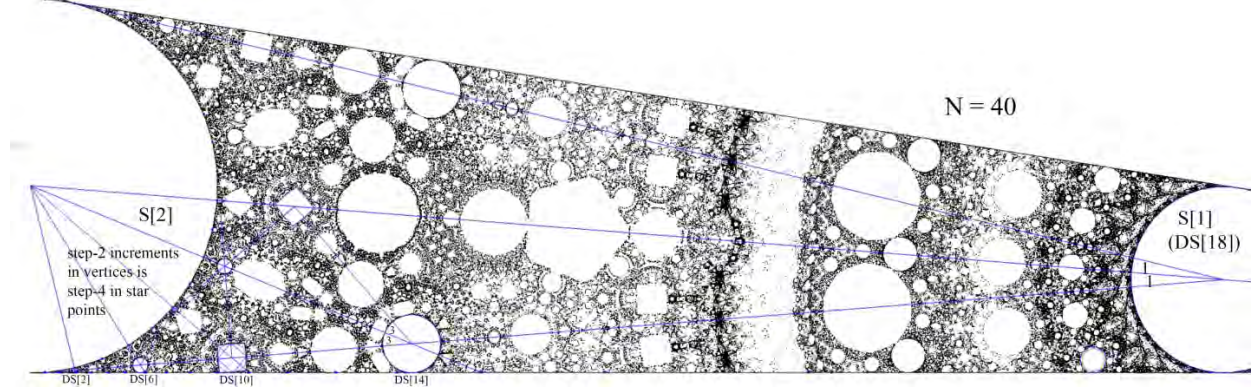

Here the DS[k] are simply the S[k] in the First Family of S[2] and the Rule of 4 says that existing DS[k] will occur mod-4 counting down from DS[N/2-4]. There is a strong resemblance here with the N-odd case because the DS[k] there only occur for k-odd and the Rule of 8 reduces to a Rule of 4 for odd DS[k]. See Figure 43.2. The real surprise is that these lines of symmetry reveal similar next-generation geometry across the whole spectrum of N-even and N-odd. Compare this geometry with N = 39 above and N = 45 to follow.

 Of course N = 40 is rife with mutations and there is no doubt that they play a role in the surrounding geometry. The graphic above does not show the large S[4] that shares its star[3] point with S[2]. Of course S[4] is mutated along with the neighboring S[5] and the resulting geometry will be described in Figure 40.3.

The basic formula for the web steps of the DS[k] is the same as the First Generation, namely $k' = N/2-k$. These same values can be read directly from the symmetry diagram, so S[1] at DS[18] is step 2 and likewise DS[14] will be step 6. Since gcd(6,40) = 2 it will not be mutated but DS[10] will consist of two squares since N/gcd(40,10) = 4. This is a canonical mutation with base spanning 10 star points – from left-side star[9] to right-side star[1]. DS[6] has gcd(14,40) = 2 and DS[2] has gcd(18,40) = 2 so they are not mutated. As shown in the enlargements below DS[2] actually has a step-4 web and this is compatible with the step-2 prediction.

**Figure 40.2.** This is just a tiny fragment of the unique geometry.

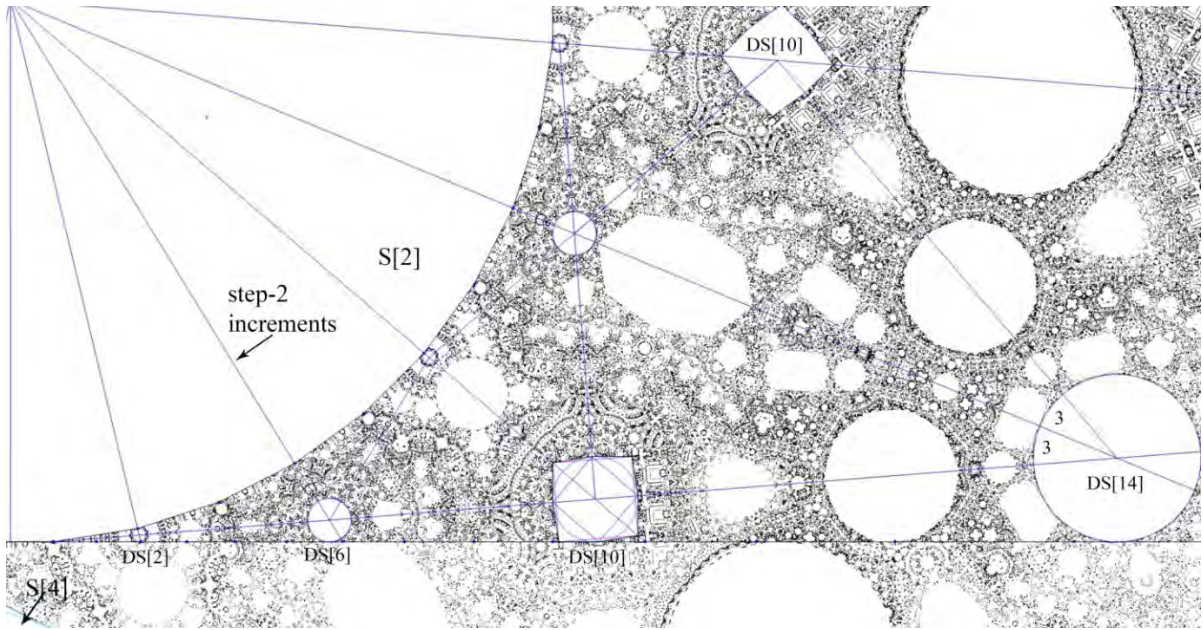

**Figure 40.2.** S[4] and S[5] illustrate the complementary 5 by 8 symmetry of the First Family

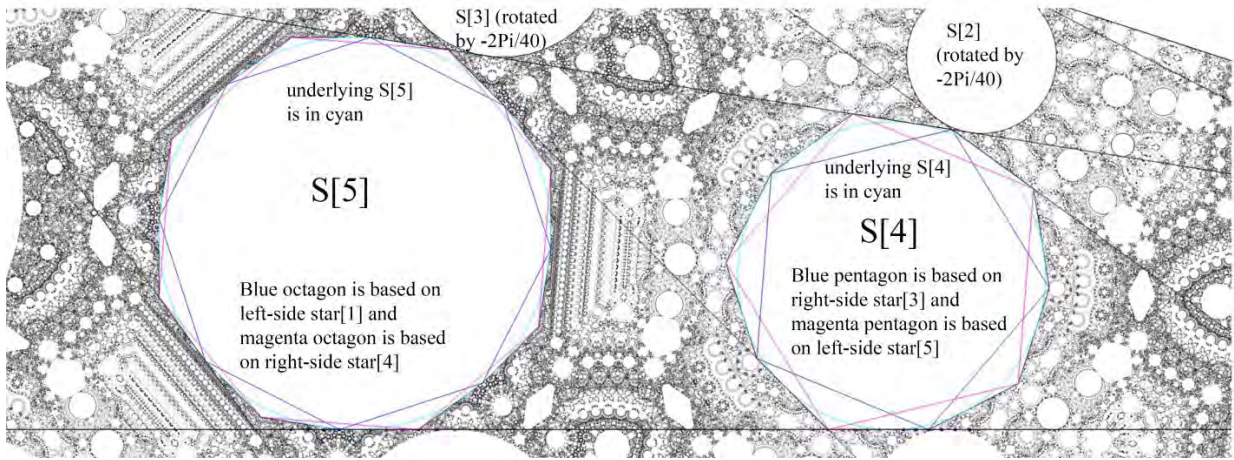

The linear structures around S[5] are found elsewhere to a smaller degree. These 'virtual edges' are more common with almost vertical edges. They 'form when $\tau^2$ 'return' maps have small offsets from the previous, so they are a similar to the 'period doubling' which is common for points inside 2N-gons.

For S[4] with k′ = 16, 40/gcd(16,40) = 40/8 = 5 so the mutation consist of two pentagons as shown here, and the base spans star[5] to star[3] . Like a good neighbor S[5] generates octagonal symmetry with k′ = 15 and 40/(gcd(15,40) = 8. The S[5] and S[4] 'towers' ending in S[1] and S[2] are no doubt strongly influenced by the mutations. To add to the mix, N = 40 has algebraic complexity 8.

These are Dc maps, but this map is only faithful to the τ-map above the x-axis, so a simple work-around to include the region below the x-axis, is to use a large crop area that includes a -2Pi rotation about the origin. Then just rotate back and crop again.

**Figure 40.2** Detail of the step-4 DS[2] geometry

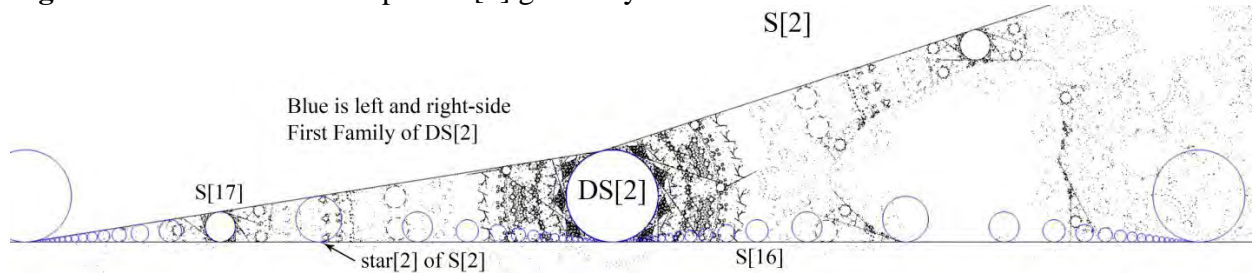

DS[2] has a virtual clone at S[19] which is the blue D tile shown on the left here. This D tile has a natural embedding at star[1] of S[2] and it has a DS[17] which survives the web and appears to generate its own step-4 web. The step-4 web of DS[2] has star[2] effective on the right but this implies that star[1] is effective on the left. There are no obvious First Family survivors but tiles like S[16] may hose surviving tiles.

**Figure 40.3** Detail of DS]2] showing the First Family

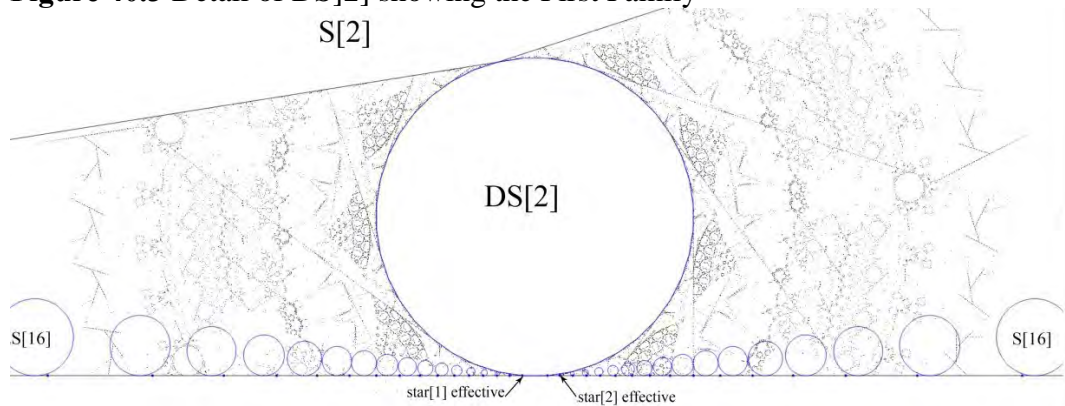

For all N, it appears that DS[1] and DS[2] will have step-2 and step-4 webs so for N-even any existing DS[1] would match up in a self-similar fashion with S[1] (and both of these will be 'mutated' into N/2 gons for N twice-odd, so the step-2 webs will essentially be step-1.) For N-odd DS[2] should have a self-similar relationship with S[1] but there are no predicted even DS[k]. The 8k family has an 'outreach' tile at S[N/2-3] of the D tile DS[2]. As shown above this S[17] has its own step-4 web and it is likely that it has some First Family survivors.

The twice-even family is the only one where S[1] has a true step-2 web and the Twice-Even S[1] Conjecture describes the possible 'step-2' families that could arise. Except for S[2] itself, there are no step-2 candidates here, but the large S[3] below always share its S[2] tile with a virtual S3x of S[1] - which has a scale[2] (D) relationship with S[3] of S[1]. This is a common feature in the 8k family.

**Figure 40.4** The S[1] region appears to be devoid of First Family survivors

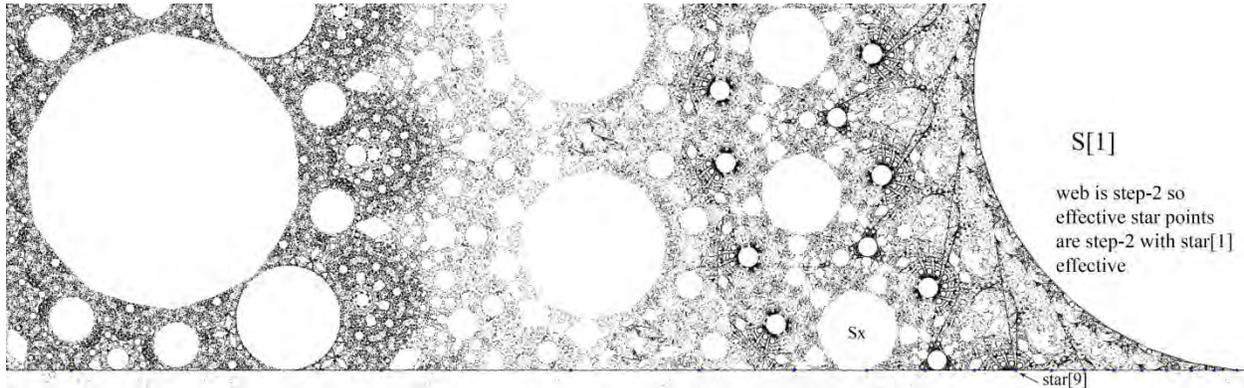

Typically these step-2 webs have little to offer in the way of surviving First Family members. As noted above the Twice-even S[1] Conjecture says that the normal S[k] in the First Family of S[1] may be replaced by Skx tiles which are conforming to S[1], but span two star points, so an S12x would span the S[12] and S[11] regions. This is not the case for the volunteer Sx but it does appear that it shares its star[13] point with the star[13] of S[1] and since it also appears to be conforming to S[1], it will also share its star[N/2-1] point with star[1] of S[1]. This allows us to find its parameters using the Two-Star Lemma.

Back in the First Family of N, the matching S[12] has k′ = N/2-k = 8, and gcd(40,8) = 5 so this tile should be the weave of two pentagons, and it is no surprise that this clone has the same mutation into green and magenta pentagons. This is not the first time that a mutation seems to facilitate the existence of a tile. N = 30 had an abundance of mutated volunteers

**Figure 40.4** An enlargement of the Sx region

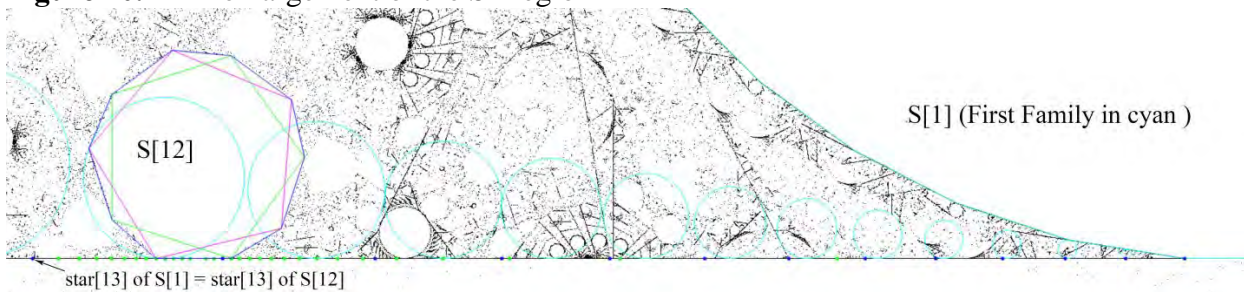

The characteristic polynomial for hSx/hS[3] is

$$-\frac{1}{128} + \frac{439x}{128} - \frac{3301x^2}{128} + \frac{12835x^3}{128} - \frac{12579x^4}{128} + \frac{3557x^5}{128} - \frac{183x^6}{128} + \frac{x^7}{128}$$

## ●N = 41

N = 41 has algebraic complexity 20 along with N = 100. It is in the 8k+1 family so there is a DS[5] and in all known cases there is also a 'volunteer' DS[2].

**Figure 41.1** The blue lines of symmetry

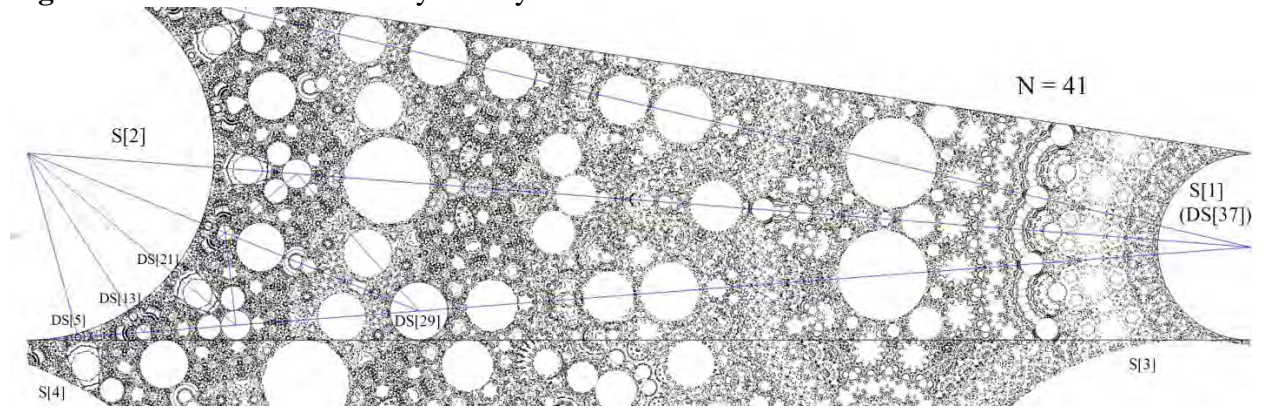

**Figure 41.2** The geometry local to S[2]

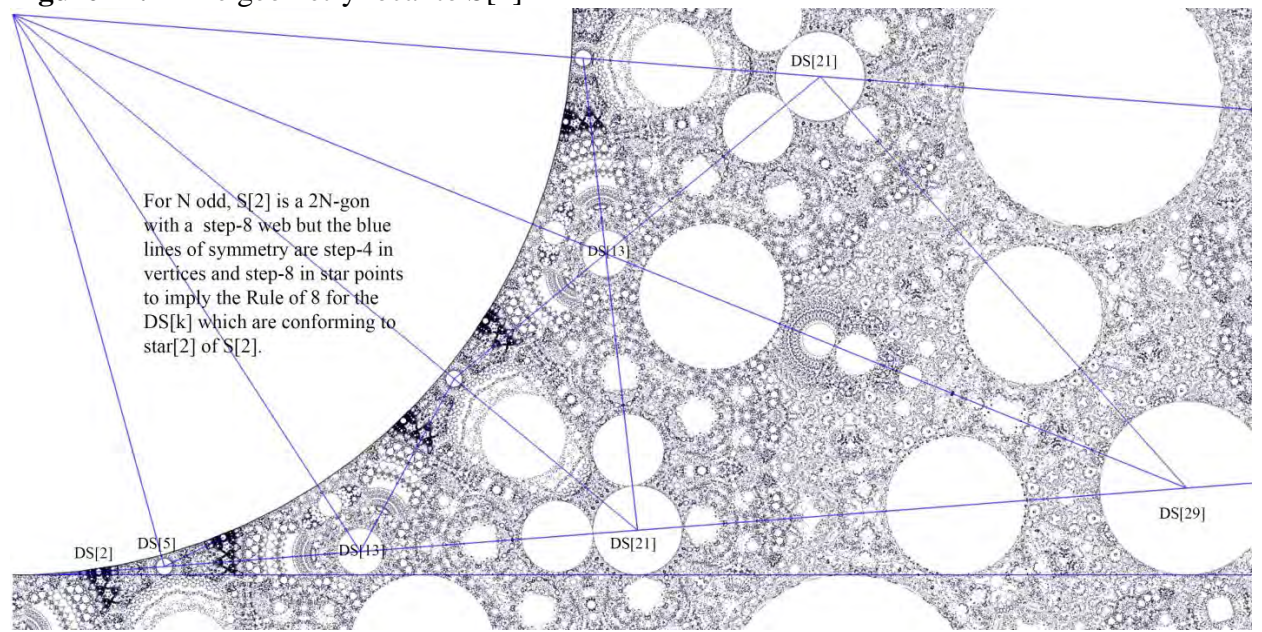

In the 8k+1 family the volunteer DS[2] has a step-4 web which has potential to support First Family members but in fact these tiles never seem to survive the limiting web. N = 41 will be no exception but DS[2] does have a very interesting Pac-Man local web as shown below.

**Figure 41.3** The geometry of DS[2] and DS[5]

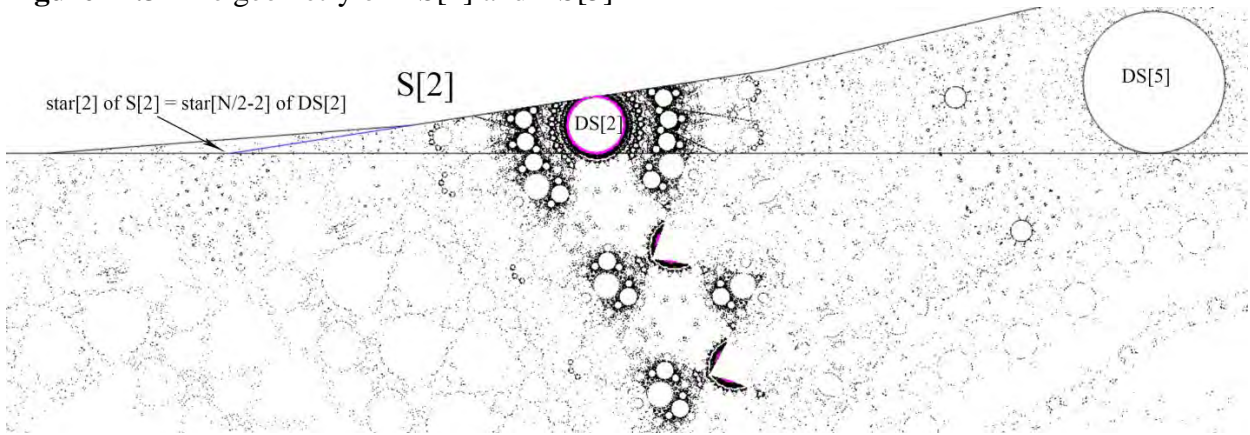

The region local to DS{2} is very dense with no surviving First Family members but the center of the virtual S[3] tile has a periodic orbit which is shown here in magenta. The orbit has period 5,536,968 = 67524*82. It is primarily local to DS[2] but it does visit the two Pac-Man satellites. Unlike periodic orbits for N = 17, 33 and N = 25, this one does not have 'period doubling', but DS[2] itself has period doubling with center period 461 = 21*41 and off-center period 1722.

**Figure 41.3** The local geometry of DS[2] showing orbits of cS[1] and cS[3]

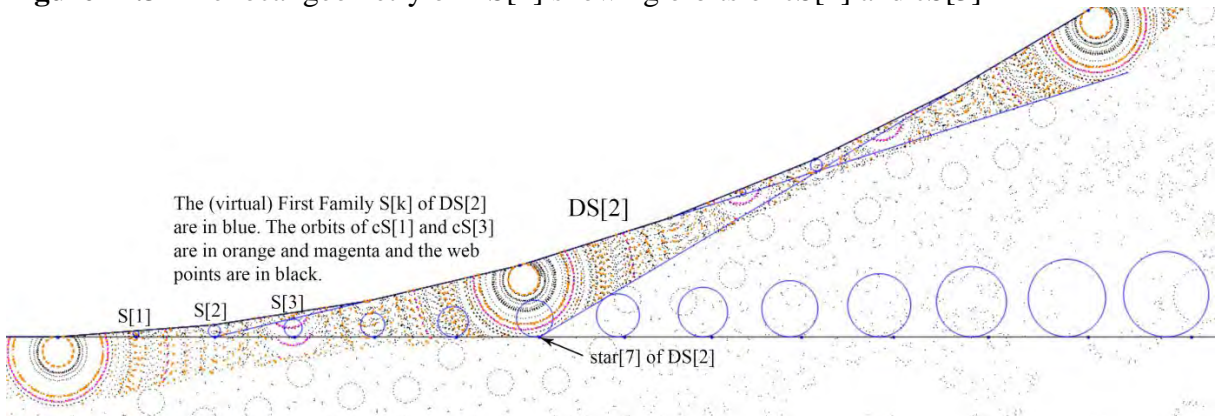

The orbit of other centers generates further rings but cS[3] is the only known periodic orbit. These rings can to get very dense, but in the limit they should have Lebesgue measure zero.

**Figure 41.5** The geometry local to S[1] is also step-4 as seen by the clusters below

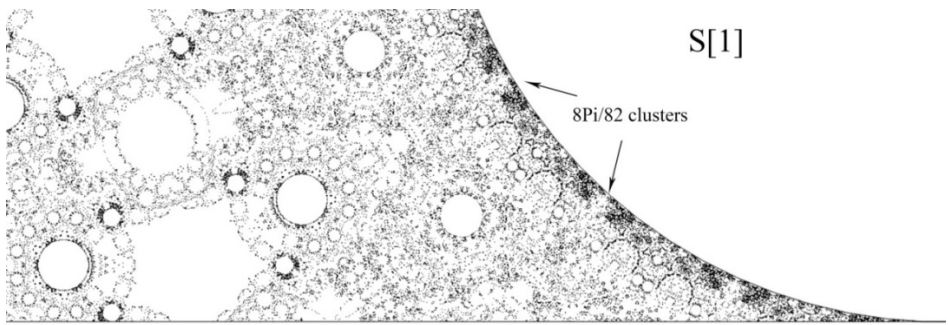

● **N = 42**

N = 42 and N = 21 have algebraic complexity 6 along with 13 ,26, 38 and 36.  N= 42 and 26 are both members of the 8k+2 family so they should have sequences of D[k] and M[k] converging to star[1] of S[2]. We noted in Part 1 that N = 18 will have mutations in S[3] and DS[3] and this will continue mod-24, so N = 42 is the next case with k′ = N/2-k = 18 for both S[3] and DS[3] with S[2] sitting-in for N. So in both cases 21/gcd(21,18) = 21/3 = 7, and they will be the weave of two heptagons with base length gcd(21,18) = 3 (relative to an N/2-gon).

**Figure 42.1**  The lines of symmetry have the same 2-step increments as N = 40 shown above. The web steps are also similar when the underlying DS[k] are used. This is illustrated with the blue underlying 42-gon of S[1].

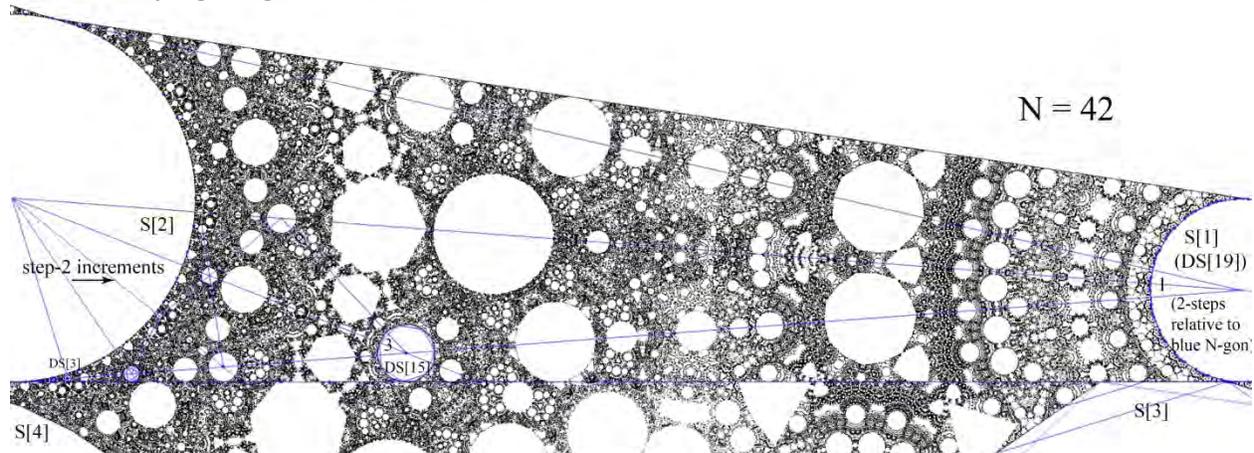

The mutation in S[3] carries over to DS[3], but since these mutations preserve centers they should have no effect on the D[k], M[k] convergence predicted by the 8k+2 Conjecture. DS[15] will also be mutated as described below. The web steps described here are the same as the twice-even case when the DS[k] are replaced with their underlying tiles, so it is possible to use the same k′ = N/2-k  predictions, with S[1] at step-2 relative to the underlying blue N-gon.

**Figure 42.2**  Detail of the S[2] region

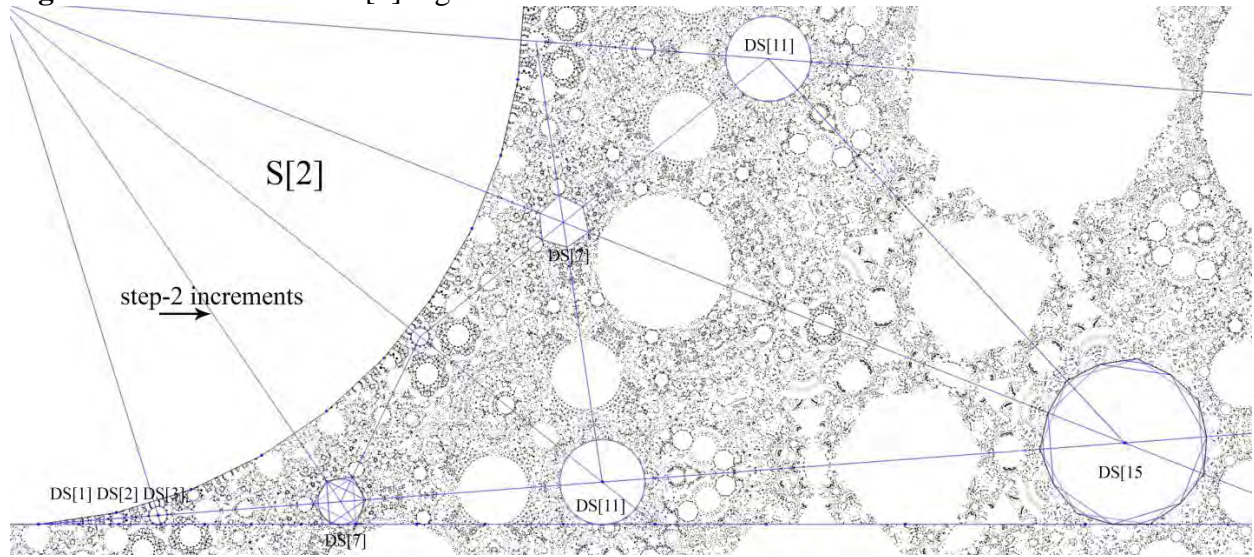

The Mutation Conjecture predicts the initial star point of the mutation using the minimum of N/2-1 -jk′ which is 20-18 = 2, but this prediction is always based on the underlying N-gon, so it corresponds to star[1] of the existing N/2-gon S[3] or DS[3]. This gives the correct mutation from star[1] to right-side star[2] as shown below. (Our First Family software automatically applies the gender-change mutations of the odd S[k] and DS[k] so the cyan polygons shown below are already mutated once-over.)

DS[15] will have the same mutation since k′ = N/2-15 = 6 and gcd(21,6) = 3. The minimum of 20-jk′ is still 2 but now it occurs with j = 3. These DS[k] mutations should be horizontal reflections of the S[k] cases because the web reverses polarity with each new generation. For DS[7] and S[7] , k′ = N/2-k = 14 so 21/gcd(21,14) = 3 and they will be the weave of two triangles with base of length 7. The Mutation Conjecture correctly predicts an initial star[3] because N/2-1-14 = 6 for the underlying N-gon (not shown here).

**Figure 42.3** Mutations in S[3], DS[15] (identical to DS[3]) , DS[7] and S[7] (not to scale). As expected the DS[k] mutations are horizontal reflections of the S[k] mutations.

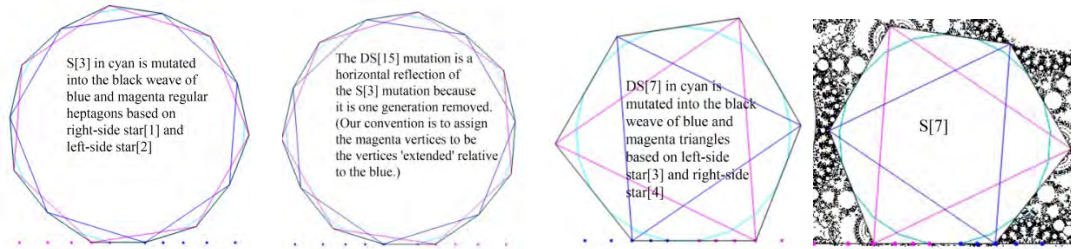

The enlargement below shows that the DS[3] mutation appears to have no adverse effect on the web evolution of DS[2] acting as D[2]. The combined 3$^{rd}$ generation M[2]-D[2] web is clearly step-4 and supports the same Rule of 4 DS[k] as M[1] and D[1] in the 2$^{nd}$ generation. Of course this D[2] web is inverted relative to S[2] as D[1] and this will continue for each generation. This should imply that 'even' and 'odd' generations will share greater similarity and it is possible that there may be a limiting form of self-similarity that emerges in the 8k+2 family. It is hard to imagine a doubly infinite family of distinct tiles emerging in this family

**Figure 42.4** Detail of the star[1] region of S[2] showing the early web evolution of the third generation presided over by D[2] and M[2] (a.k.a. D2 and M2).

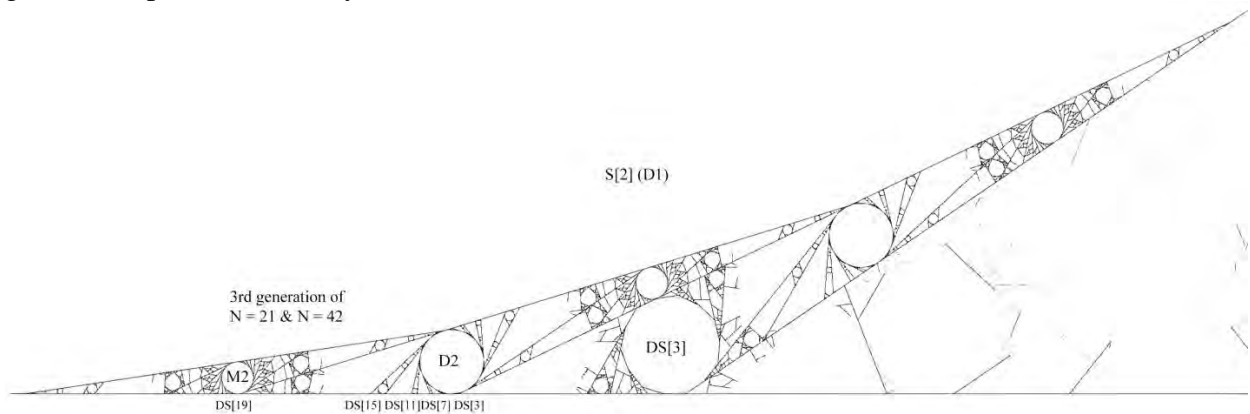

● N = 43

N = 43 has complexity 21 along with N = 49. N = 43 is in the 8k + 3 family so the first DS[k] predicted by the Edge Conjecture is DS[7]. This is the most 'extreme' S[2] tile gap of any N-gon and one issue is whether a family structure can arise here without any local DS[k]. Early versions of the local web showed promise of a possible DS[1] but in the plots below we show that this DS[1] is indeed virtual and it is not clear what extended family structure exists. The second issue is whether there are always volunteer conforming tiles between existing DS[k]. In some cases these volunteers could be odd DS[k] beyond those predicted by the Edge Conjecture. For N-odd there will always be three such candidates between predicted pairs

**Figure 43.1**  The level-21 web in magenta

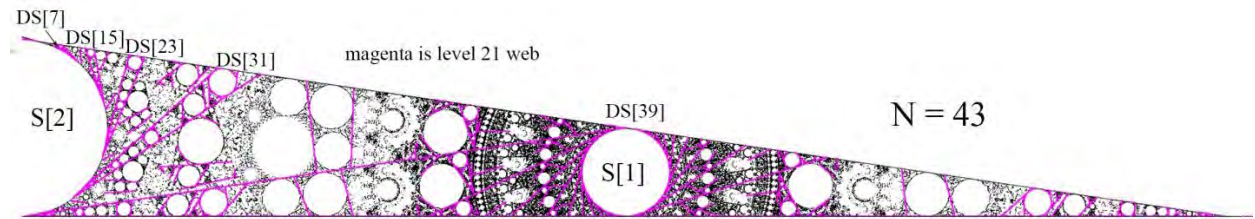

**Figure 43.2**  All possible DS[k] (with k odd). This would apply to any N-odd case.

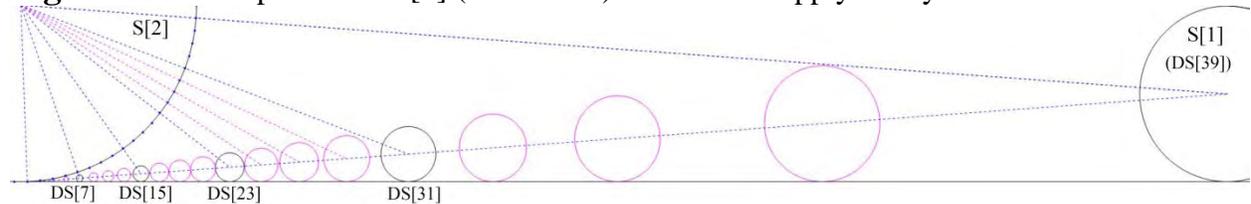

N =47 to follow has a volunteer adjacent to the predicted DS[19] which is a close match for the DS21] and below we can see that that N = 43 has a similar volunteer which is a close match with a DS[11] which would sit between the existing DS[7] and DS[15]. There are also near hits with DS[17] and DS[19].

**Figure 43.2**  The lines of symmetry

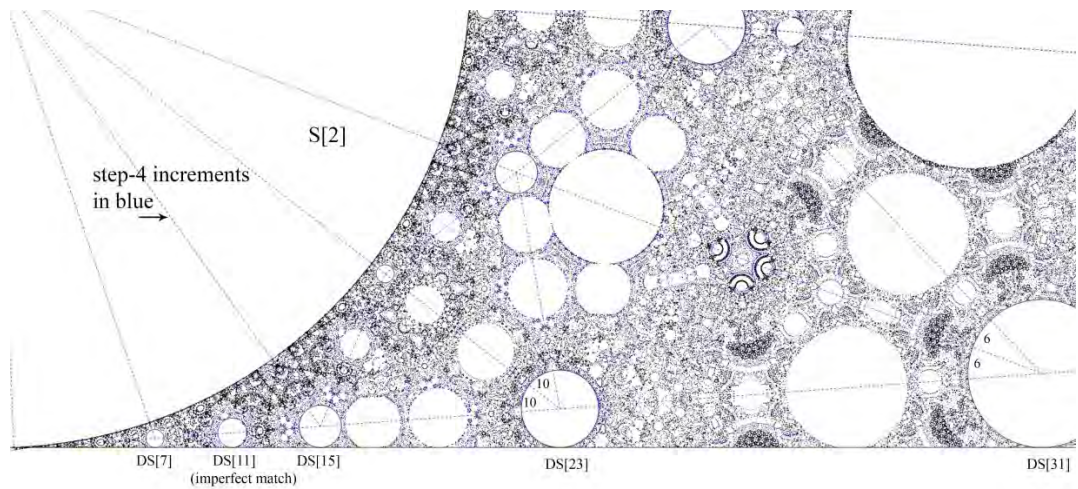

**Figure 43.3** The star[2] region of S[2] shows no obvious signs of extended family

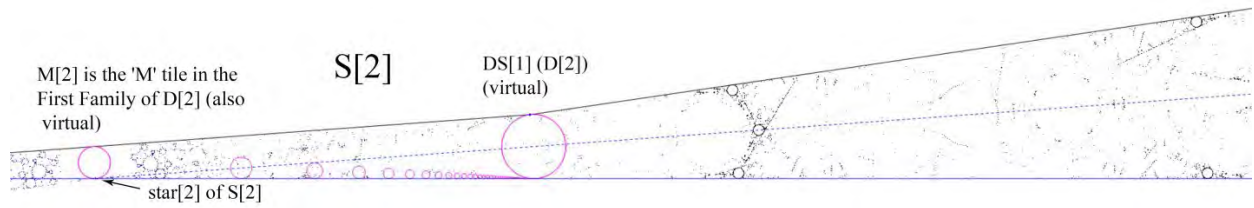

This region originally looked promising for the existence of a volunteer DS[1], but after extensive web iterations, DS[1] was shown to be virtual. The matching M tile in the First Family is the S[1] of S[2] which we call the M[2] of S[2]. It looks more promising that DS[1] – but it is also virtual. We also generated the First Family of the virtual M[2] but found no obvious hits with the existing web structure (which took 3 days of round the clock iteration with the Dc Map.) . This is not promising to the 8k +3 family – at least local to S[2].

Since the web local to S[1] is step-4, the 'effective' star points will include star[1] of S[1] as shown below. This is not always a promising alignment and there are no obvious First Family survivors. However every effective star point provides the potential for extended tile structure, and N = 43 clearly provides a diverse framework for 'alien life forms'.

**Figure 43.2** The step-4 web local to S[1] (which is a 2N-gon)

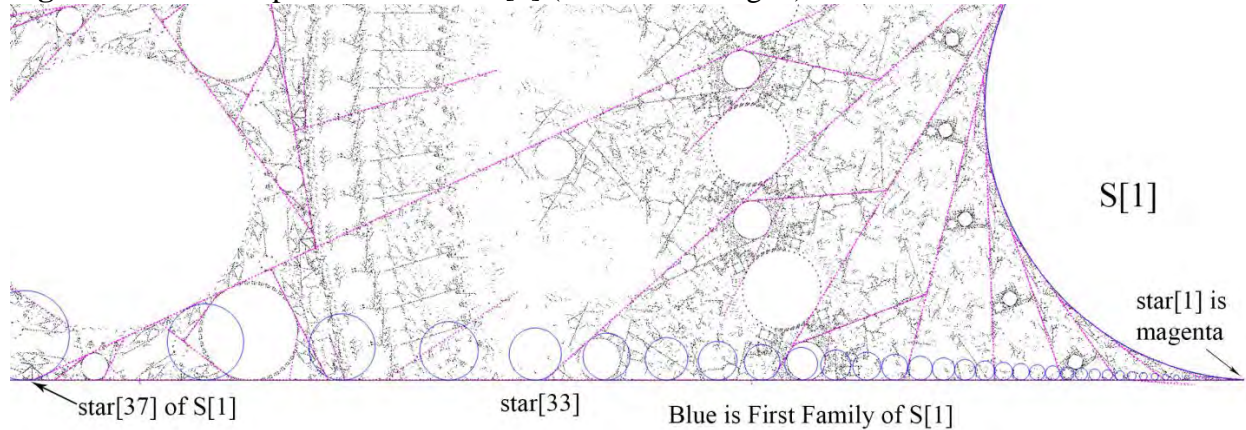

● N = 44

N = 44 is degree 10 and the 5<sup>th</sup> non-trivial member of the 8k+4 family and the 3<sup>rd</sup> member of the mod-16 subclass containing N = 28 and N = 12. Therefore the 8k + 4 Conjecture says that S[2] will be mutated into blue and magenta N/4-gons and that the DS[4] predicted by the Edge Conjecture will have a Px 'parent' which has a DS[4] as S[19]. We will use this fact to construct Px. As always the two star points that DS[4] and Sx share are GenStar and star[3] of DS[4] which match up with star[1] and star[19] of Px: Therefore $hPx = d/(Tan[19\pi/44]-Tan[\pi/44])$ where $d = StarDS4[[3]][[1]]- StarDS4[[19]][[1]] \approx 0.000656836$

**Figure 44.1** The blue lines of symmetry

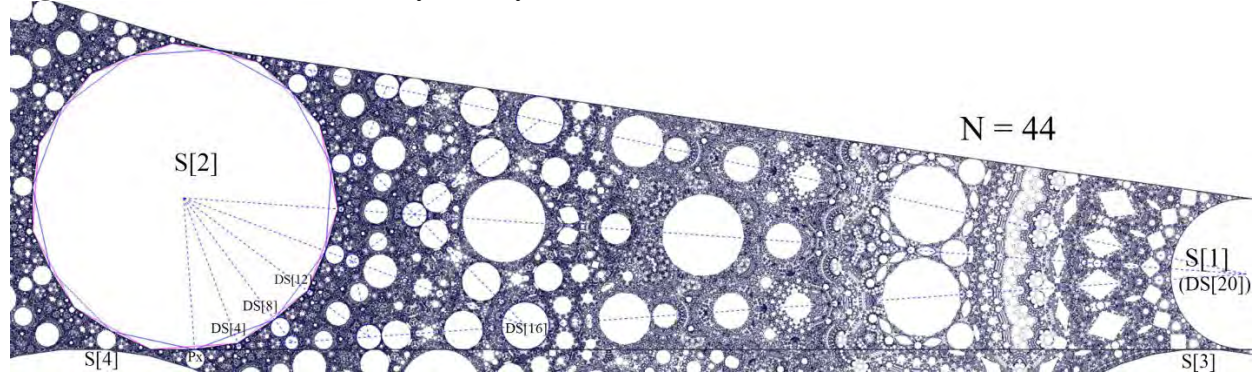

**Figure 44.1** The S[2] region

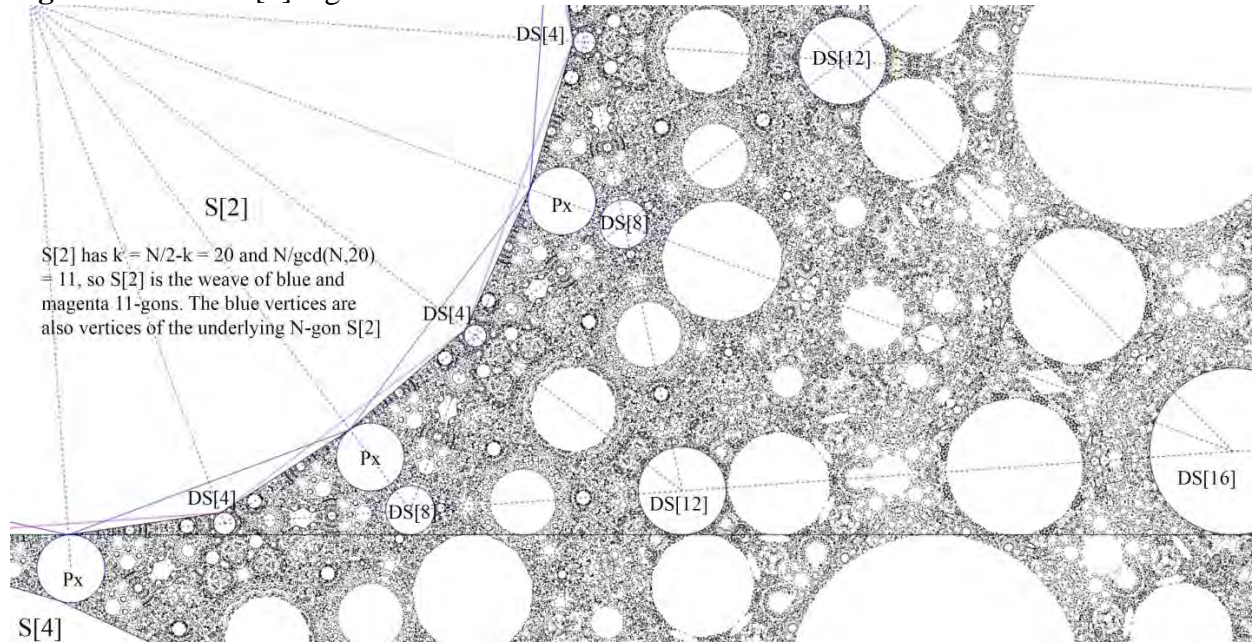

For N = 28, Px and the rotated S[4] almost shared a vertex and here the match is exact since Px[[1]] = S[4][[37]], but the characteristic polynomial of Px does not indicate a close algebraic relationship with any other tile except DS[4]. Since Px also shares a vertex with the blue hendecagon embedded in S[2] we checked to see if Px could possibly be in the First Family of this tile, but there was only a close match with S[11].

**AlgebraicNumberPolynomial[ToNumberField[hPx/hN,GenScale],x]** yields

$$\frac{57731}{131072} - \frac{12997321x}{131072} + \frac{87263363x^2}{32768} - \frac{552905261x^3}{32768} + \frac{2194226605x^4}{65536} - \frac{1542371279x^5}{65536} + \frac{215163863x^6}{32768} - \frac{22993969x^7}{32768} + \frac{3007163x^8}{131072} - \frac{13089x^9}{131072}$$

For the ratio hPx/hDS[4] the last term is just $\frac{6561x^9}{131072}$ because DS[4] is in the First Family of Px but relative to S[2] and S[4] this term is larger at $\frac{22321x^9}{131072}$ and $\frac{72427x^9}{1441792}$.

As always Px is a natural candidate for generating a 'next-generation' at star[1] of S[2] but as shown below the tiles that survive in the web have no clear relationship with either Px or DS[4].

**Figure 44.2** The 'next-generation' geometry of S[2] yields no obvious candidates.

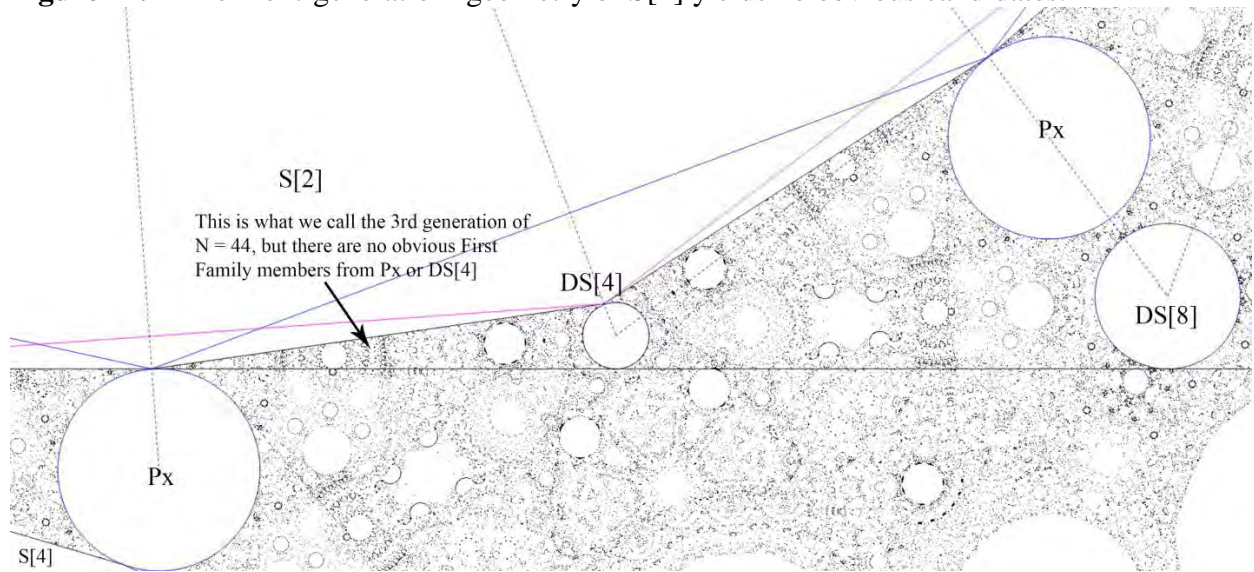

The Twice-even S[1] Conjecture says that the step-2 web of S[1] can support 'step-2 'tiles such as S3x and S9x shown below. S3x will always be an S[2] in the First Family of S[3], but S9x shown here is only virtual because the early web splits – most likely due to S[3] influence.

**Figure 44.3** Detail of the step-2 S[1] web

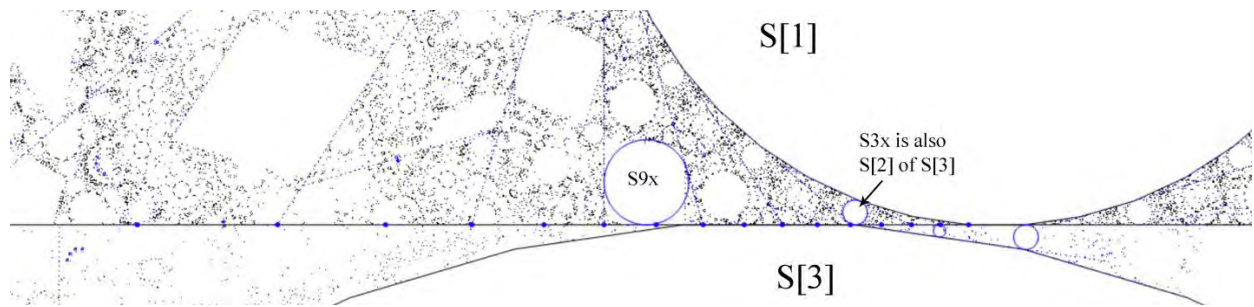

● **N = 45**

N = 45 has complexity 12 and is in the 8k+5 family. Therefore the S[2] step-2 family will include a DS[1] with a step-2 web.

**Figure 45.2** A symmetry diagram for the 6 DS[k]

The overall geometry appears to be similar for all odd N. The combined S[1],S[2] web will be step-8 and the Rule of 8 DS[k] will be 2N-gons which evolve in a simple k′ =2N/2-k cw fashion, so S[1] at DS[41] is step-4 as shown here, and DS[33] will have step-12 symmetry. DS[1] will be maximum step 44 which is 22 top and bottom. Since Mod(2N,N-1) is always 2, DS[1] will have a simple step-2 web as shown below. There should be mutations here because 2N/gcd(2N,k) is 10, 18 and 30 for DS[9], DS[25] and DS[33]. The DS[9] case is shown below.

**Figure 45.3** The mutation of DS[9] into the weave of blue and magenta decagons

For DS[9], k′ = N-k = 36 just like DS[3] of N = 39. Here gcd(N,k′) is 9 so the mutation is based on two regular decagons. As with N = 39, the expected surviving star points should match the minimum of N-2 - jk′ which is 7. The base should span 9 star points (with included edge of the underlying cyan DS[k]), so the right-side star point is star[2].

There are also mutations among the volunteers and DS[25] and DS[17] appear to support a non-regular hexagon in a fashion similar to DS[19] and DS[11] for N = 39. It is no surprise that the shared symmetry between the 8k+1,8k+3, 8k+5 and 8k+7 classes translates into similar geometry. Compare the DS[33] region above with DS[27] of N = 39 and the connections between DS[25] and DS[17] with the matching DS[19] and DS[11] for N = 39. In addition the 'volunteer' DS[1] of N = 39, has a step-2 web just like the 'real' DS[1] here.

**Figure 45.3** The DS[1] (D[2]) region showing the step-2 web which has the potential to foster Fist Family tiles.

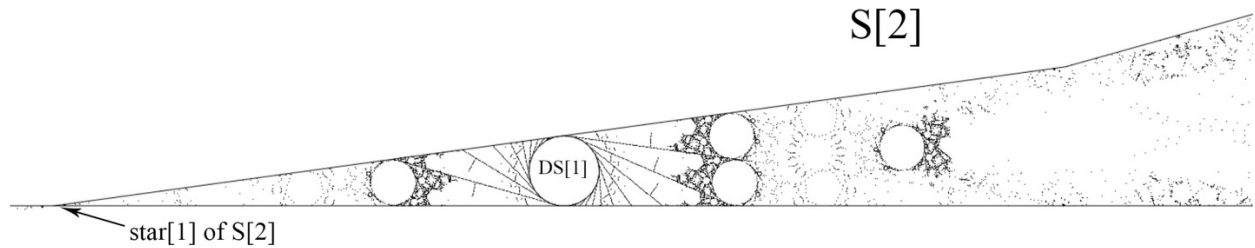

The large volunteer on the right may be a failed DS[3]. Below we show the First Family of D[2] which includes a prominent 'M' tile at star[2]. This virtual M tile is an M[2] relative to S[2]. It has a promising local family structure which is on the border of the star-2 region defined by S[2]. The local geometry at star[2] of S[2] is a literal reflection of the geometry at star[1] of N, and in the Dc world this point has coordinates {-1,0} while star[1] of N is {0,0}.

When N is odd S[1] and S[3] will always share their star[4] points because these points are also star points of S[2] and N.

**Figure 45.4** The web local to S[1] showing the shared geometry with S[3]

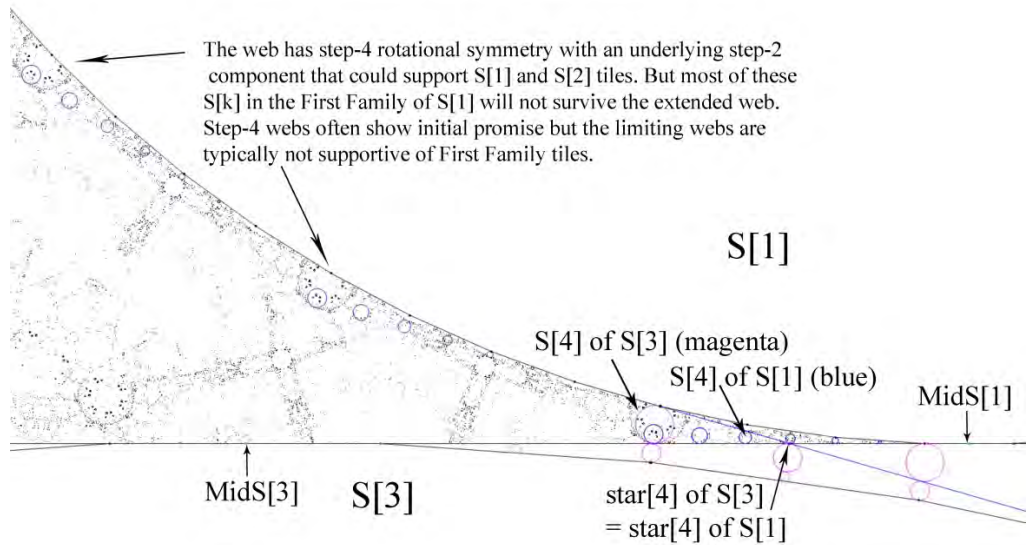

- **N = 46**

N = 46 has complexity 11 and is a member of the 8k+6 family along with N = 14, 22, 30 and 38 This family has a prominent DS[1] which could serve as an M[2] but typically there is no matching D[2]. As indicated earlier N = 46 and N = 38 have very different embedded N/2-gons and since DS[1] is also an N/2-gon, this may have an effect on the geometry of DS[1]. It appears that the DS[1] of S[2] always has a step-1 web in a manner similar to S[1], but the issue is whether this web can support any First Family tiles or any potential self-similarity.

**Figure 46.1** The early web

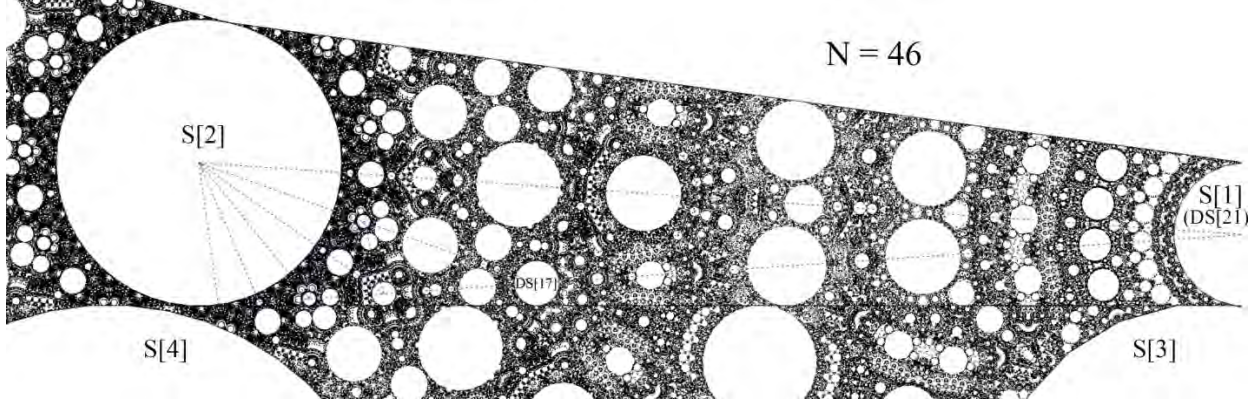

**Figure 46.2** The region local to S[2]. The large volunteer below DS[5] is similar to the Mx of N = 38

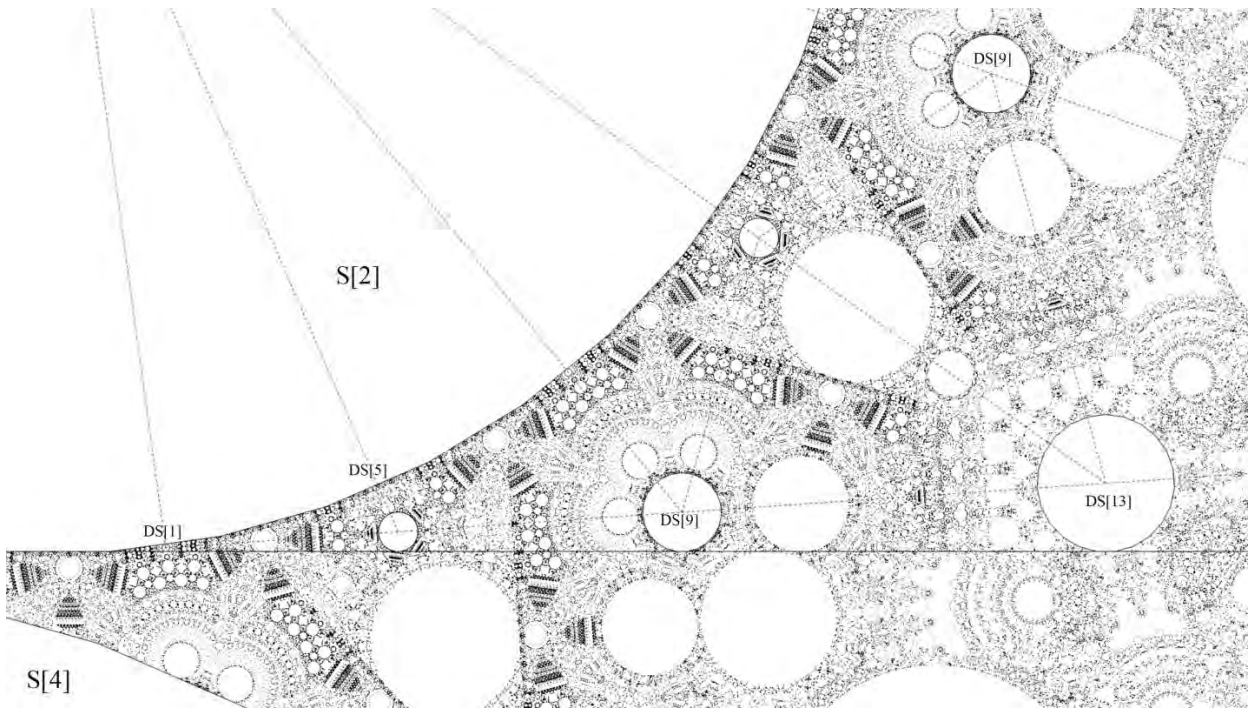

**Figure 46.3** The region local to DS[1] showing the step-1 web

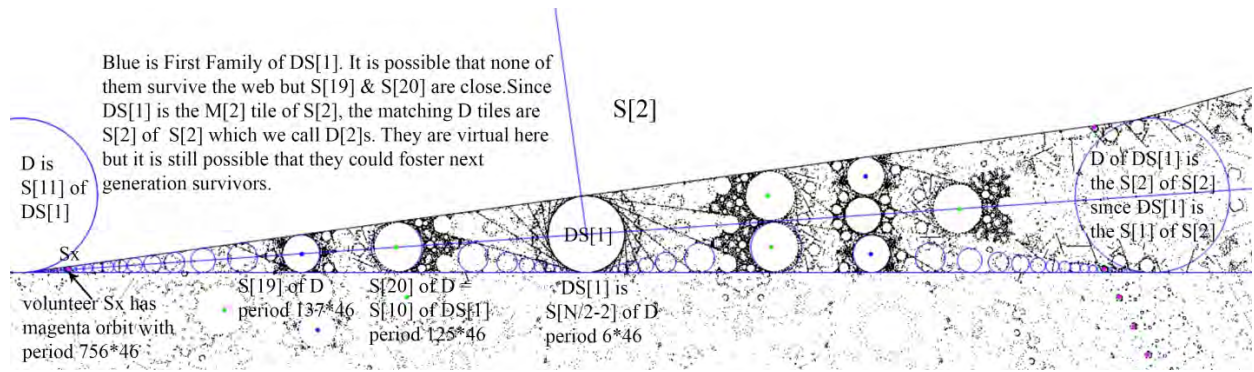

DS[1] always has a step-1 web in the 8k + 6 family but such webs do not mimic the step-1 web, of N because DS[1] is a satellite of S[2] which is itself a satellite of N. Relative to N, all the points in S[2] have the same predictable step-2 orbit with period N/gcd[2,N] = 23 but any point outside S[2] will have a more complex orbit.

The points inside DS[1] can be regarded as 'satellites' of S[2] because their orbits will repeat after 6*46 iterations, but points outside of DS[1] are not well-behaved 'satellites of satellites' and in fact the orbit of the center of the S[2] tile of DS[1] has no known period and could be 'dense' in the large-scale invariant region of N that extends out to S[11] of N. The first 6 million points in this orbit are shown below. It should be clear that even though this center is exact, we have no idea what the actual tile looks like and it could be just a single non-periodic point.

**Figure 46.3** The orbit of the center of the S[2] tile of DS[1]

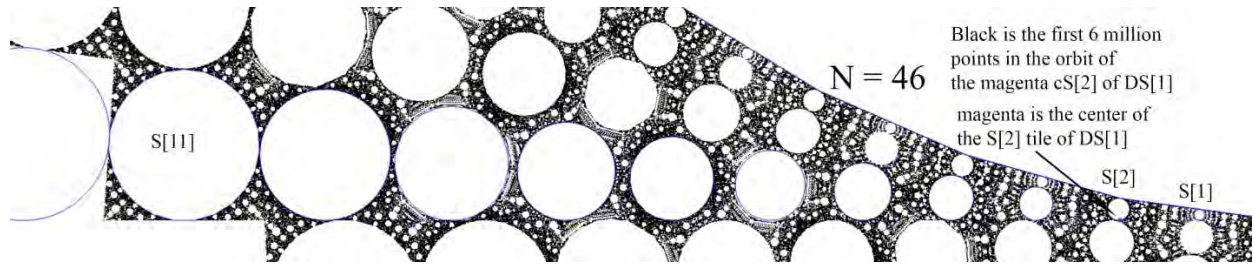

This behavior is not unusual for points around DS[1] for any member of the 8k+6 family, but as shown above the S[10] tile of DS[1] is actually a close match with the web and the period is just 126*46. This is the only tile that is in the individual First Families of both D and DS[1], but by default we use their combined families for analysis. The neighboring S[19] of D is an even closer match with the web and it has period 137*46. For tiles closer to D the best we can do is the volunteer Sx which is about half-way between S[6] and S[7] with period 756*46. The S[2] of D has an orbit that probably would rival the orbit of the S[2] of DS[1] shown above. The critical star[1] of S[2] is congruent to star[1] of N and like all star points the orbital periods will approach infinity with the temporal scaling being the main issue. Even in the well-behaved 8k+2 family we would expect a spectrum of sales since the geometry is typically multi-fractal.

Notice in Fig. 46.3 above that the S[10] of DS[1] and the S[19] of D are both displaced horizontally by the web. These small horizontal shifts are fairly common and the plot below shows how this can occur through the expected outward displacements of their star points.

**Figure 46.4** The geometry of the DS[1] family showing the virtual D at S[11] and the neighboring S[19] and S[20] tiles of D. Relative to DS[1], this S[20] of D is an S[10] and this is the only tile that is shared by the families of D and DS[1]. In the limiting web these two both experience a slight inward (right-side) displacement and we can see here how this occurs by tracking the early magenta and red orbits of p1 and p2 ,

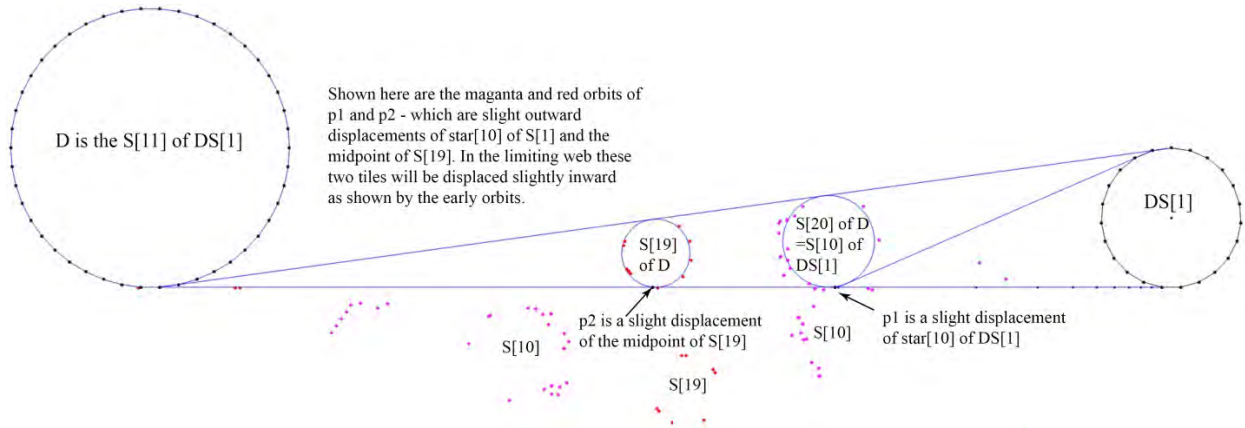

If DS[1] were at the origin as an N/2-gon, D and these other two tiles would be a perfect fit for the web. In this case the points on the horizontal 'base' line will experience an outward displacement of sDS[1] because this is a clockwise web obtained by iterating the trailing edges of DS[1] under $\tau$ and as the reference point for $\tau$ changes by one step from star[1] of DS[1] to star[2] the corresponding shear will be sDS[1] outwards (and the corresponding rotation for an S[k] would be $k' = (N/2-2k)$ times $\pi/(N/2)$ which would be $k' = 1$ and 3 for D and S[10] and they would both be N-gons – but of course S[19] of D would be an N/.2 gon).

Here the initial horizontal base line is unchanged when DS[1] is embedded in N = 46, but the actual shear for $\tau$ will be sN, so the blue star points of DS[1] will have an effective displacement outwards which is simulated here by the displacements of p1 and p2. (The actual displacements used here are very small at p1 = star[10]-{.000002,0} and p2 = MidDS19-{.0000032,0}.) These two points were then iterated about 10 million times under $\tau$ to generate the magenta and red points shown above.

If a similar iteration was done for points closer to DS[1] , the early points might rotate randomly around DS[1] like they would do for a true S[k], but soon they would diverge in a manner similar to the virtual S[2] of DS[1] as shown above in Fig. 46.3

- **N = 47**

N = 47 has complexity 23, which is the largest of any of the N-gons in this survey. It is a member of the 8k+7 family so there is a DS[3] and matching 'volunteer' DS[1]s. The 8k+7 Conjecture predicts that these DS[1] will always exist as S[N-3] tiles in the First Family of the 2N-gon DS[3]. If we regard DS[3] as a parent the matching D and M tiles will be S[N-1] and S[N-2] and this S[N-3] will be the unique tile in the First Families of both D and M.

**Figure 47.1** The geometry of the star[2] region for N odd

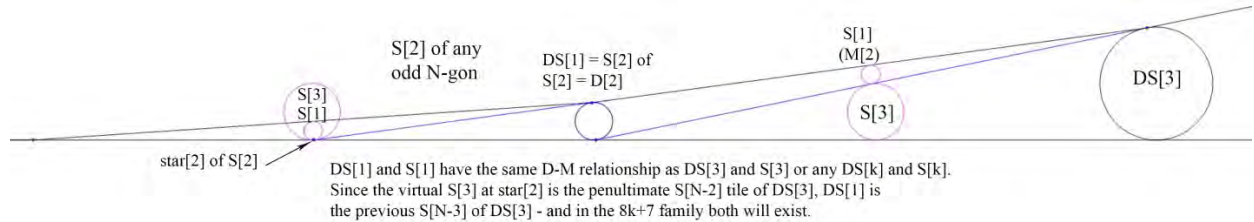

Every DS[k] of S[2] is a 'D' tile relative to the matching S[k] of S[2] but these S[k] never seem to exist . Here DS[3] is an S[N-2] 'D' tile relative to the virtual S[3] of S[2] and this implies that so this is an N-odd version of the D-M relationships for N even. S[2] here is a 2N-gon and so together they form a D-M unit similar to S[2] and S[1]

**Figure 47.2** The web local to S[2] and S[1] showing lines of symmetry in blue

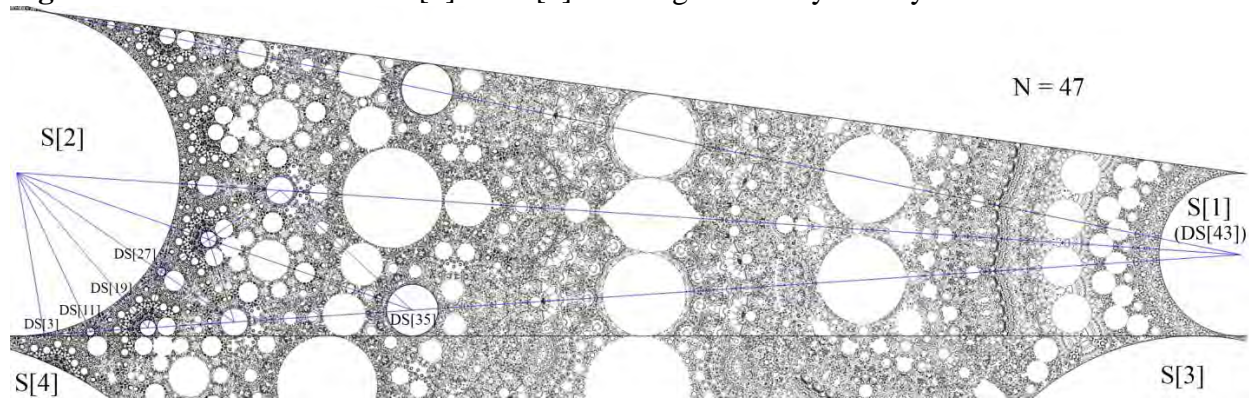

For N odd, the combined S[1]-S[2] web local to S[2] is step-8 so the DS[k] survivors exist at 16*Pi/94 rotations about the center of S[2]. This means that the rotated DS[11] will be aligned with DS[27] and the rotated DS[19] will align with DS[35]. This DS[19] supports a large colony of tiles and has a neighbor which is a slightly displaced DS[21].

**Figure 47.3** The web local to S[2]

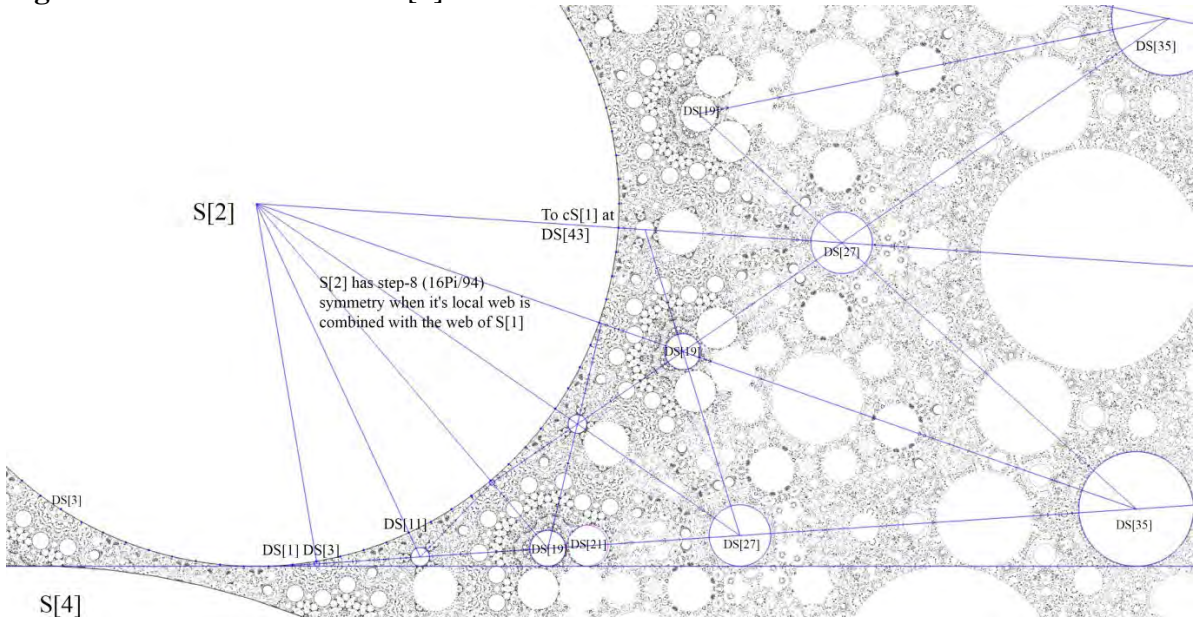

Each DS[k] has an internal step-k′ symmetry which increases mod-4 from S[1] at step-2. Therefore DS[35] will be step-6 as shown by the blue lines of symmetry – which are 6 steps each. This matches the DS[27] from N = 39 which was mutated into two 13-gons since 78/gcd(12,78) = 78/6 = 13. Likewise DS[27] matches the DS[19] from N = 39 at step-10 , This is still step-10 since k′ = N/2-k = 28 and Mod[94,28] = 10.

This volunteer DS[1] is also the S[2] of S[2] which we call D[2]. It has a step-2 web that probably does not support S[k] First Family tiles. The step-6 web of DS[3] counts down from star[N-3], so the smallest effective star points are periodic {0,2,4} starting with N = 15. Therefore N = 47 should have star[2] of DS[3] effective just like N = 23 earlier. This could provide support of First Family members, but these tiles do not survive so we turn to DS[3]

**Figure 47.4** The DS[3] region showing the step-6 web development and satellites

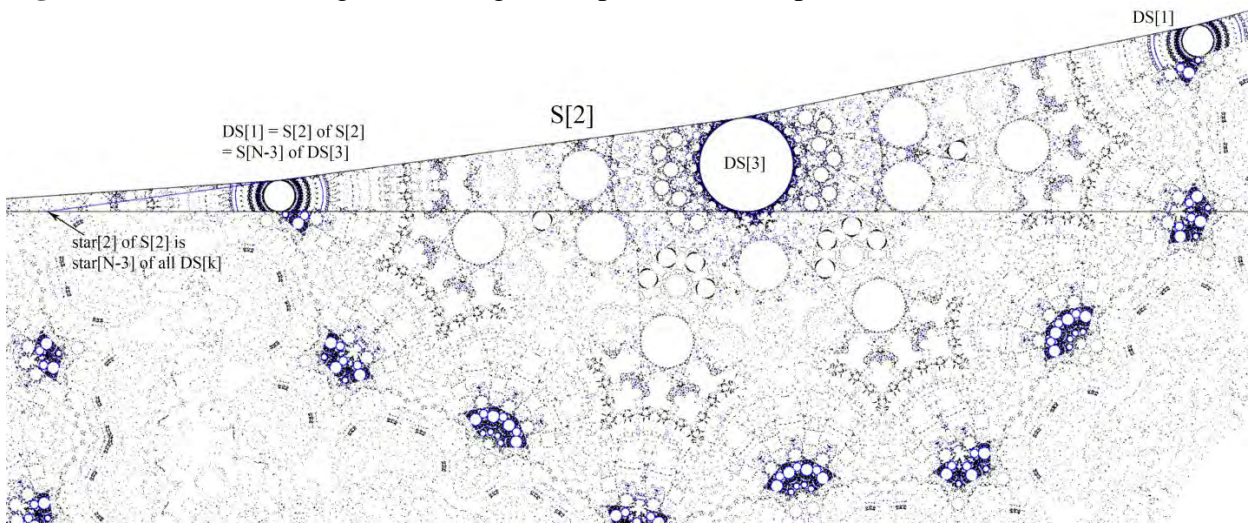

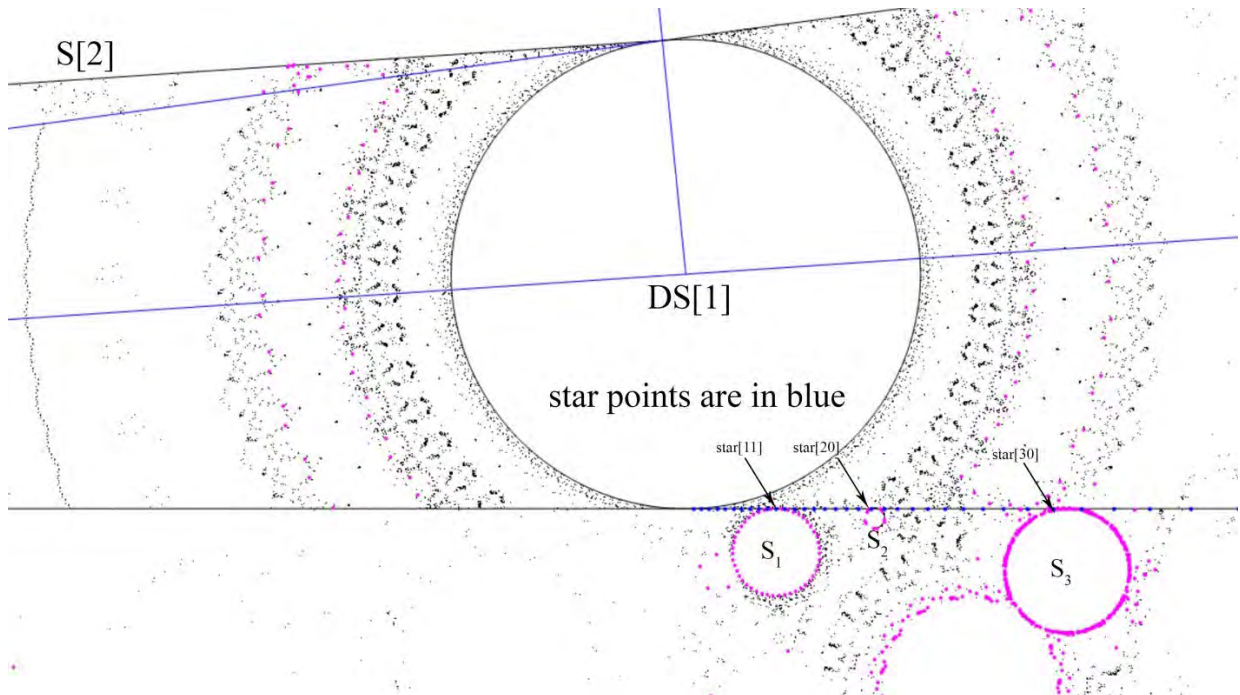

**Figure 47.2** The web local to DS[3] showing potential survivors – which actually do not exist

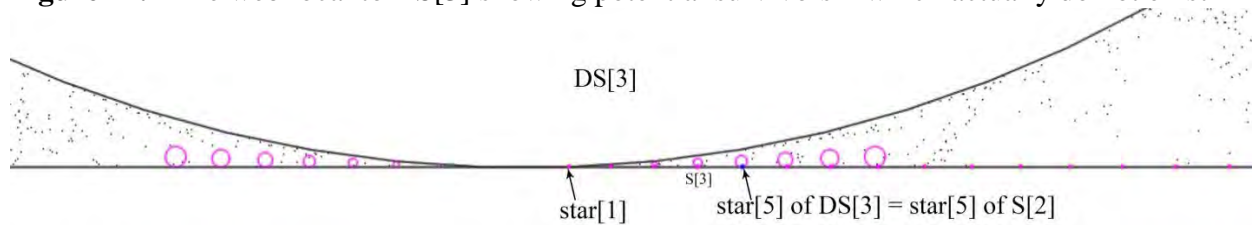

In the 8k+7 family DS[3] always shares its star[5] point with S[2] as shown above, but this is never synchronized with the step-6 web. However as noted above, the right-side star[2] of DS[3] should be effective for all N of the form 23 + 24k and this creates a potential geometry similar to N = 23 where S[1], S[2] and S[4] survive in the 2N-gon family, but in fact there are no known survivors.

S[1] shown below is also a 2N-gon but it has a step-4 web. S[1] shares the same penultimate N-2 star point as all of the DS[k], so star[1] will be effective just like DS[1]. This is not necessarily conducive to a local family structure but there may be a surviving S[2] shown below in magenta.

**Figure 47.3** The web local to S[1] showing the (ideal) First Family of S[1]

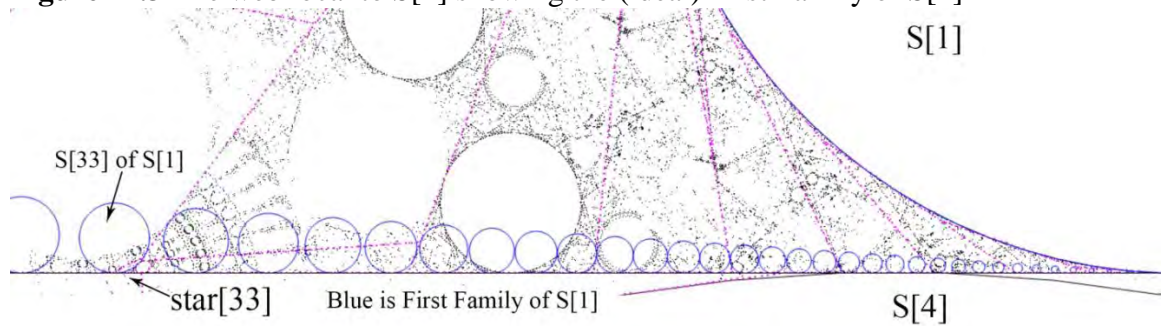

● N = 48

N = 48 has complexity 8 along with N = 40. They are both in the 8k family so they have a DS[2] which can serve as a second generation D[2] but there is little evidence that its step-4 web can support any traditional First Family tiles of DS[2].

**Figure 48.1** The web local to S[1] and S[2] showing mutations in DS[18], S[3] and S[4]

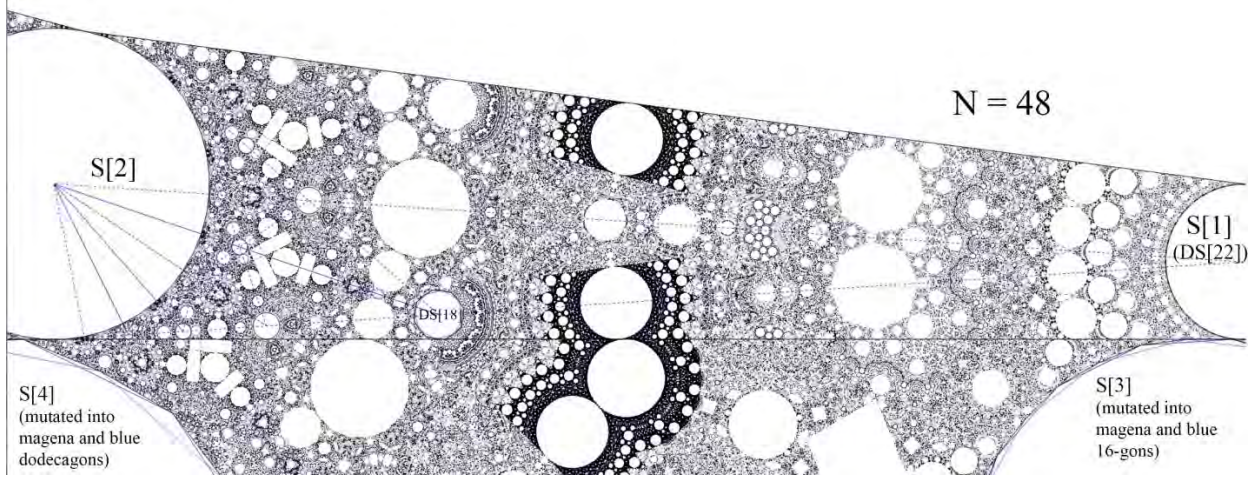

By the Rule of Four the surviving DS[k] include DS[2], DS[6], DS[10], DS[14] , DS[18] as well as S[1] at DS[22]. When N is even these DS[k] are simply the S[k] tiles in the First Family of S[2]. With S[2] playing the part of a surrogate N, the step sequences of the DS[k] are the same as the S[k] of  N, namely $k' = N/2-k$. This means that mutations in the S[k] of the First Family of N = 48  should be matched by mutations in the DS[k]. The most important S[k] mutations here are S[3] and S[4] of N with $k' = 21$ and 20 respectively. For S[3], N/gcd(N,21) = 16 so S[3] will be the barely perceptual weave of two 16-gons as shown above. The case for S[4] is more extreme with N/gcd(N,20) = 12 and there is no doubt that this will impact the web local to S[2].

**Figure 48.2**  The web local to S[2] showing the step-8 symmetry which translates to a step-4 DS[k] evolution.

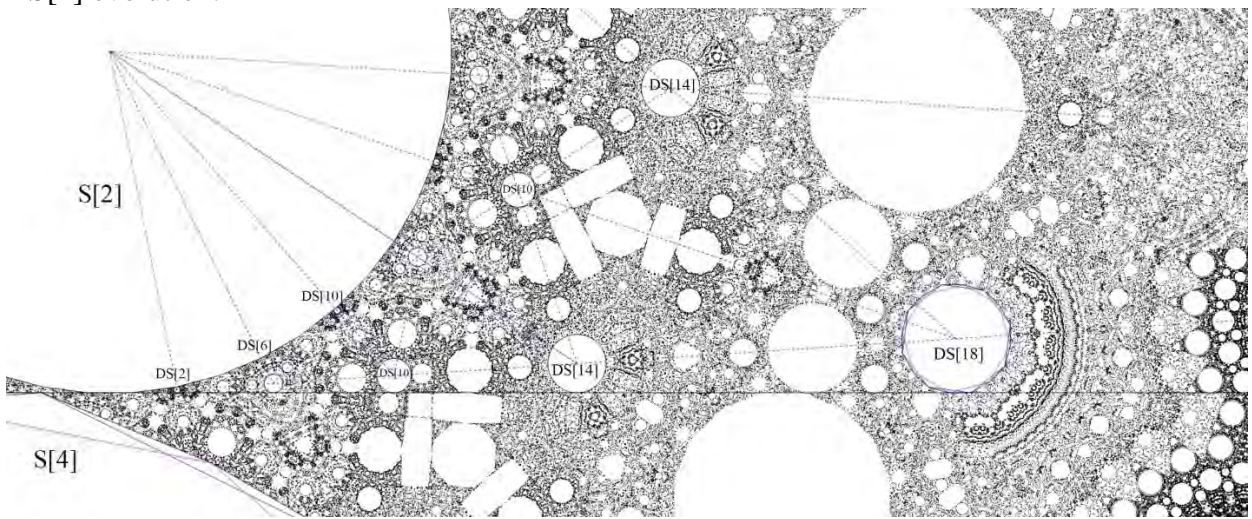

Since the mutation condition for N twice-even is $gcd(N,k') > 2$, the DS[2], DS[6], DS[10] and DS[14] tiles escape, but for DS[18], $N/gcd(N,6) = 8$ so it is weave of two octagons and this should match the mutation of the S[18] tile of N. For N twice-even this penultimate DS[N/2-6] will always have step 3 symmetry with respect to the underlying cyan N-gon and this symmetry should be preserved by mutations. There appear to be no survivors of the First Family of DS[18]

**Figure 48.3** The S[18] tile of S[2] shows clearly the contributions from S[1] (primarily left-side) and S[2] (right-side). By contrast the large 'volunteer' shows little sign of left-right bias.

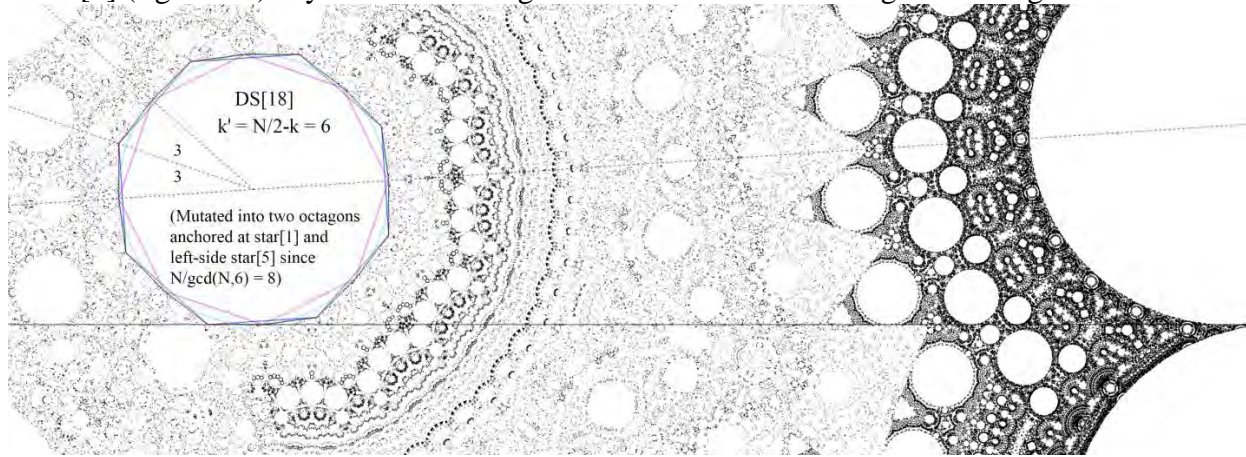

This DS[18] lies on the boundary of the 'sphere of influence' of S[2]. For N twice-even, S[1] is still in the First Family of S[2], but the First Family of S[1] does not include S[2], and the local web of S[1] is even more 'autonomous' than the twice-odd case. In all the combined S[1]- S[2] webs we have seen, the DS[k] feel a divided influence from S[1] and S[2] and it would be expected that a 'boundary' tile like DS[18] will have the most extreme right-left dichotomy.

**Figure 48.4** The 3$^{rd}$ generation D[2] tile has no known surviving First Family members but the large Dx tiles are based on the virtual D[3] of D[2]. These Dx are N/2-gons with their own families. Since D[2] evolves cw step4, the Dx tiles exist at $-8\pi/48$ rotations about cD[2] but they map to each other in a very non-trivial fashion, beginning with just the bottom 3 as a triad.

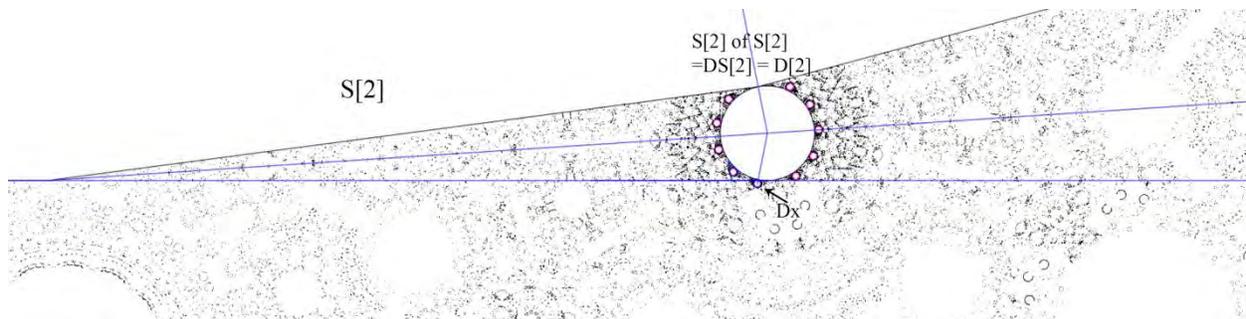

The period of the center of Dx is 860*48 and by comparison, the period of D[2] is 6*48 so these orbits clearly wandered well beyond the edges of S[2]. Part of this orbit is shown below in magenta – where the black background points are part of 400 million Dc iterations of a single point in the 'halo' surrounding Dx – namely star[15] of D[3]. This orbit may be non-periodic yet it is linked to the 41,280 Dx clones.

**Figure 48.5** The magenta periodic orbit of Dx superimposed on the black orbit of star[15] of D[3] – which is in the 'halo' around Dx. This orbit may be non-periodic.

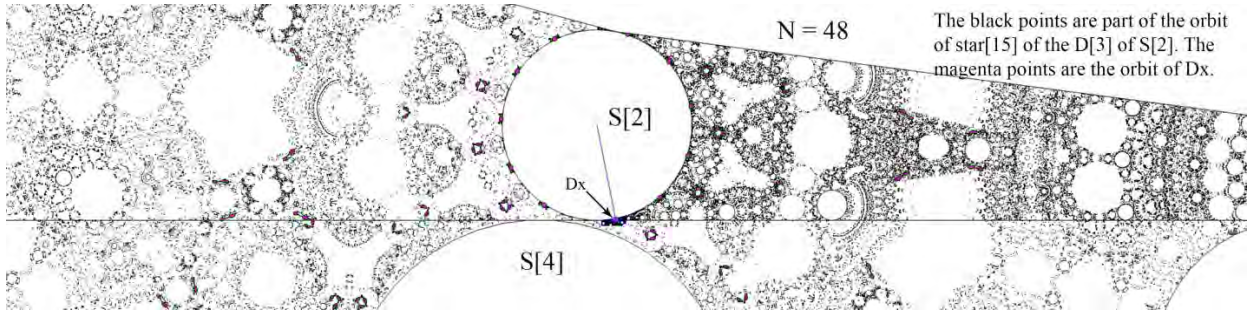

**Figure 48.6** Detail of the 'volunteer Dx which shares its star[2] with D[2] and the virtual D[3]. Note the blue extended edges which become part of the local web of Dx. This web in turn fosters surrounding tiles and this geometry is repeated throughout the invariant region local to N. N = 39 earlier has a Dx which appears to be formed in a similar fashion from 'left-over' snippets of extended edges.

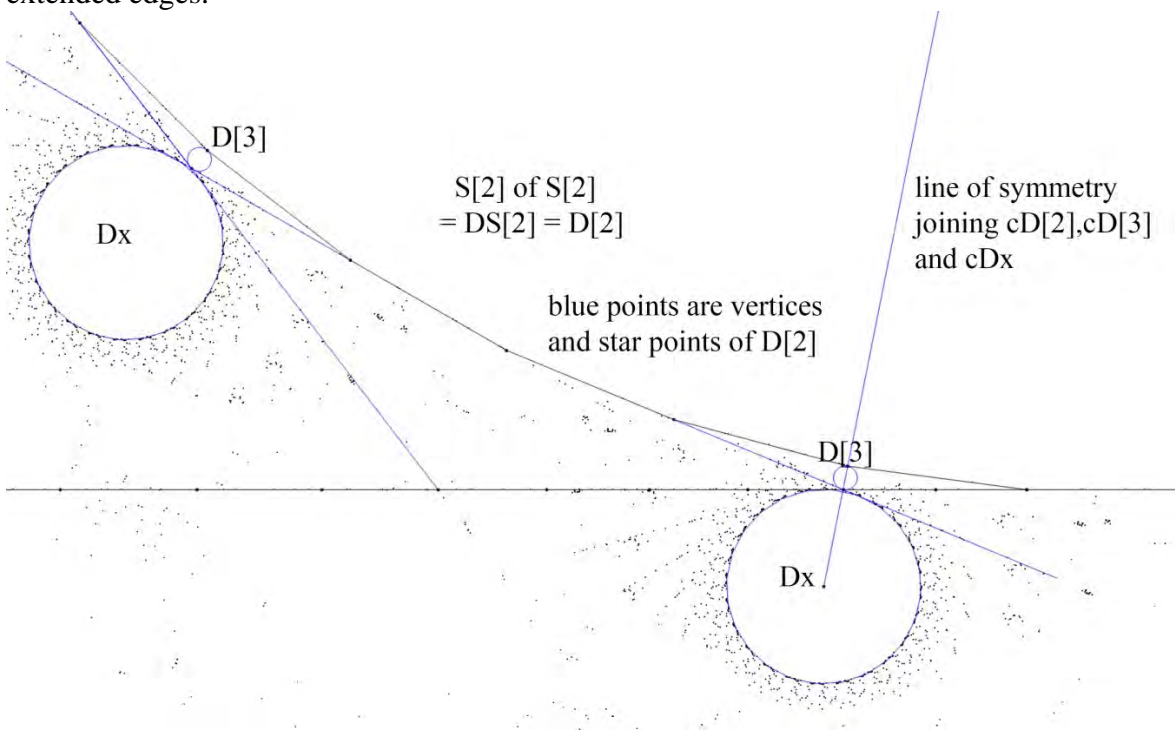

Using a modified Two-Star Lemma we have formulated a plausible regular 24-gon for Dx using star[3] and star[5] of D[2]. It is shown in blue above. Note the arc of 'triads' which exist at each vertex of Dx. Both S[2] and D[2] have similar triads.

● **N = 49**

N = 49 has complexity 21 along with N = 43. Algebraically these two are among the most extreme cases that we will see here. In terms of edge geometry, N = 49 is in the 8k+1 dynamical family so the Edge Conjecture and the Rule of 8 predict the existence of S[1] at DS[45] as well as DS[37], DS[29], DS[21], DS[13] and DS[5] as the main members of the modified First Family of S[2]. The 8k+1 Conjecture predicts that there will be a 'volunteer' DS[2].

**Figure 49.1** The web symmetry local to S[2]

It is fitting that DS[21] is mutated into blue and magenta heptagons as shown below. Since $k' = 49-21 = 28$ and $2N/(\gcd(2N,28)) = 98/14$, DS[21] will be the weave of heptagons with a base spanned by 14 star points. The Generalized First Family Theorem of [H5] predicts that in the local step $k'$ web, one surviving 'effective' star point will be the absolute minimum of $N-2-jk'$ which here is $|47-56| = 9$, so the matching opposite-side star point will be star[5]. This yields the blue and magenta heptagons, which must be compatible with the lines of symmetry.

**Figure 49.2** Enlargement showing the mutation in DS[21] and 'towers' of (virtual) S[k]

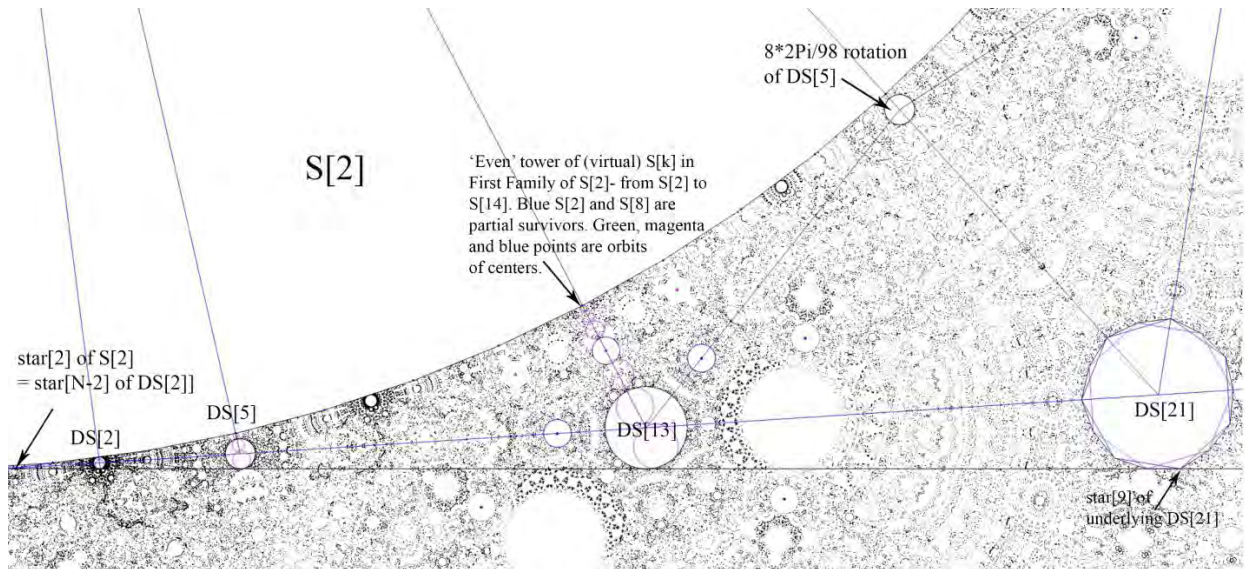

In Section 2 of Part 1, we discussed the 'even' and 'odd' 'towers' of S[k] that arise from the rotational (step-1) symmetry of N. Modified towers can also arise from multi-step symmetry. S[2] here has partial step-4 symmetry and full step-8 symmetry, so it is possible that some S[k] in the First Family of S[2] could survive the web. To test that hypothesis we have constructed virtual S[k] 'towers' above, based on the step-8 DS[k].

Looking at the DS[13] tower above, the 'base' is the (unrotated) magenta S[14] tile of S[2]. The First Family Theorem says that relative to S[2], $hS[k] = \text{Tan}[\pi/98] \cdot \text{Tan}[k\pi/98]$ which we abbreviate $s_1 \cdot s_k$, so relative to S[2], $hS[14]+hS[12] = s_1 \cdot s_{14} + s_1 \cdot s_{12}$ which is well approximated by $2hS[13]$. But S[13] and DS[13] have an M-D relationship in size which means that $hDS[13] = hS[13]/\text{scale}[2]$ where $\text{scale}[2] = s_1/s_2$ approaches ½ with N (For N = 98 scale[2] ≈ .49968.)

This explains the close fit in the three cases above. Note that the S[2] and S[8] tiles in this tower may be partial survivors and in this fashion such towers may help to explain the geometry.

All the magenta tiles in these towers are either S[k] of S[2] or k*2Pi/98 rotations of such S[k], so their parameters are known and we can iterate their centers. These S[k] will inherit the period doubling of S[2] so all of-center points will have orbits which visit S[k] twice in a symmetric fashion relative to the center so they will have even periods while the center has odd period.

To avoid confusion we refer to the S[k] of S[2] as Skx  Here cS2x has period 91826 = 916*98 so this is not a true center -which is half-way from cS2x to $\tau^{45913}$ of cS2x. It is no coincidence that this true (green) center is on the line of symmetry joining cDS[13] and cS[2]. For cS8x the corresponding period is 15974 = 163*98 and the offset is comparable to that for S2x. Therefore the true (magenta) center is known and we can use the web it estimate its radius. This blue candidate we call Sy and of course its center is also on the line of web symmetry.

At first glance it seems that the two rotations of Sy about cDS[13] would also be web survivors, but in fact they have slight offsets and distinct orbits with period 11074 = 113*98 for off-center points and 5537 for the (blue) center. We call these two Sz and Szz. It is entirely possible that Sy and the Sz may be mutated 2N-gons. There is a hint of this in the web.

**Figure 49.3**  The DS[2]-DS[5] region showing their towers

The basic geometry here is identical to all members of the 8k+1 family but the issue is whether any members of the First Family of DS[2] or M[2] survive the web. The situation with DS[2] is typically not very promising but there may be some hope of survivors in the First Family of the virtual M[2]. Here a slight displacement of the S[23] of M[2] survives with period
30674 = 313*98 and another survivor at about S[16] with period 208446 =2127*98.
It is not surprising that neither of these have period doubling because M[2] is an N-gon.

For DS[2] the web steps are k′ = N-2 and Mod[2N,N-2] is always 4 so this DS[2] will have a step-4 web in a manner similar to S[1]. Like all next-generation DS[k], DS[2] will have a cw web orientation with left side star[N-2] equal to star[2] of S[2] as shown above. These step-4 webs of DS[2] are deceptive because they initially show promise of support for First Family tiles, but the limiting webs that develop are very complex and have no clear support for such tiles. As N grows, it takes many hours of round-the clock iteration to get '3$^{rd}$ generation' detail only to discover that we know very little about the resulting geometry. Here there seems to be a 'odd' 'tower' with the virtual S[1] and S[3] tiles of DS[2] and a volunteer which we call Sx.

**Figure 49.4** The geometry local to DS[2] showing 'volunteer' tiles Sx and Sxx

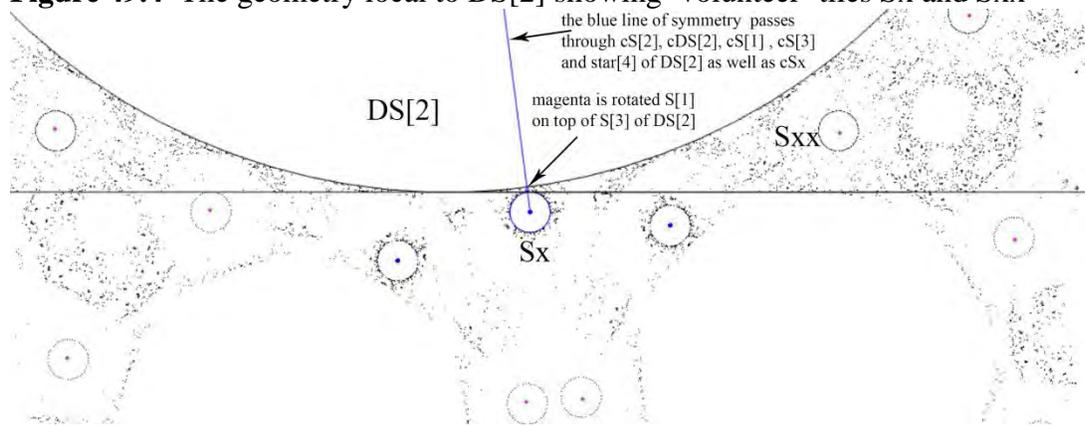

The magenta virtual S[1] and S[2] of DS[2] are barely visible at this resolution. Sx may have a vertex at star[2] of S]3] which is star[4] of DS[2] but it is difficult to verify this with the current web resolution. Because Sx has period doubling its center is known and it must be a 2N-gon, which does not normally form an 'odd' tower. The period of cSx is 11445 so all other points in Sx (and the two matching tiles) will have period 22890. This is an unusual period because it is divisible by 7 and not 49. This seems to imply that it is a 'resonant' orbit which visited only one of the two identical heptagons embedded in N = 49. The Sxx have period 150038 = 98*1531 with no doubling so they should be N-gons. We do not know their exact centers. They map to each other at step-4 but they are not strictly satellites of DS[2] since their orbits 'wander'. The period of cDS[2] is 1225 = 25*49 and it must have doubling.

**Figure. 49.5** The web of S[1] is step-4 with star[3] effective as with DS[2], but once again there appear to be no First Family survivors.

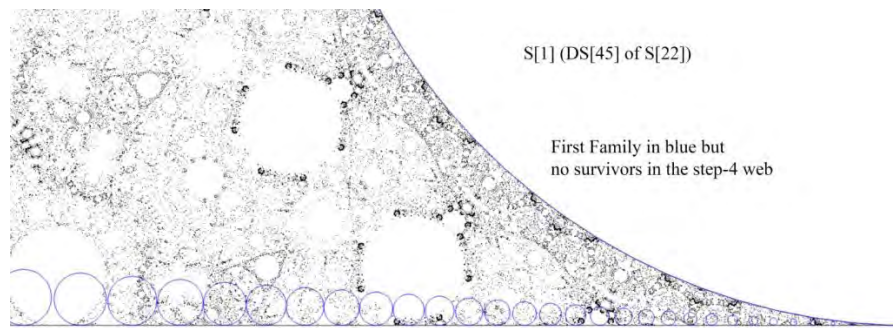

● N = 50

N = 50 and the matching N = 25 have complexity 10. Since N = 50 is in the 8k+2 family, there should be sequences of M[k] and D[k] tiles converging to star[1] of N and star[1] of S[2]. There is little doubt that these sequences exist here, but as N grows it becomes increasingly difficult to generate the local webs and track the generations converging to star[1]. The crude plot of the 4th generation in Figure 50.5 below took a week of round the clock calculations.

**Figure 50.1** - The second generation on the edges of N = 50.

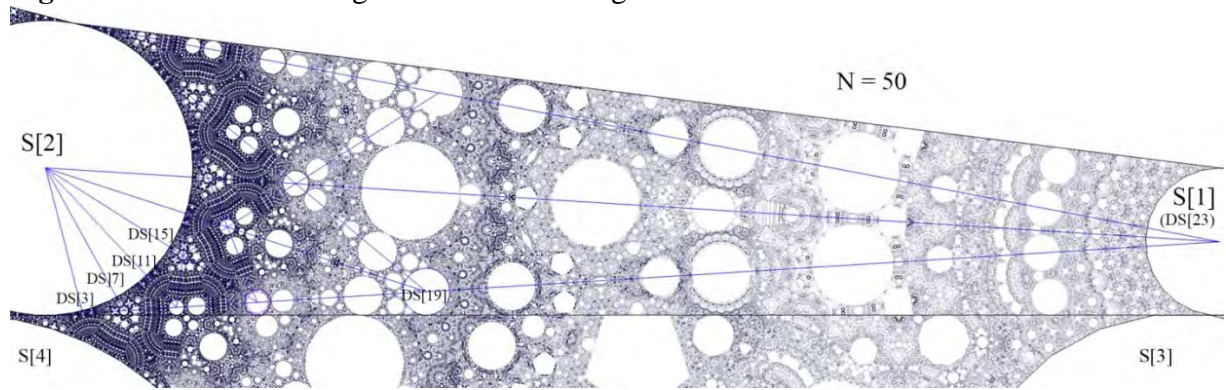

**Figure 50.2** The step-4 symmetry local to S2

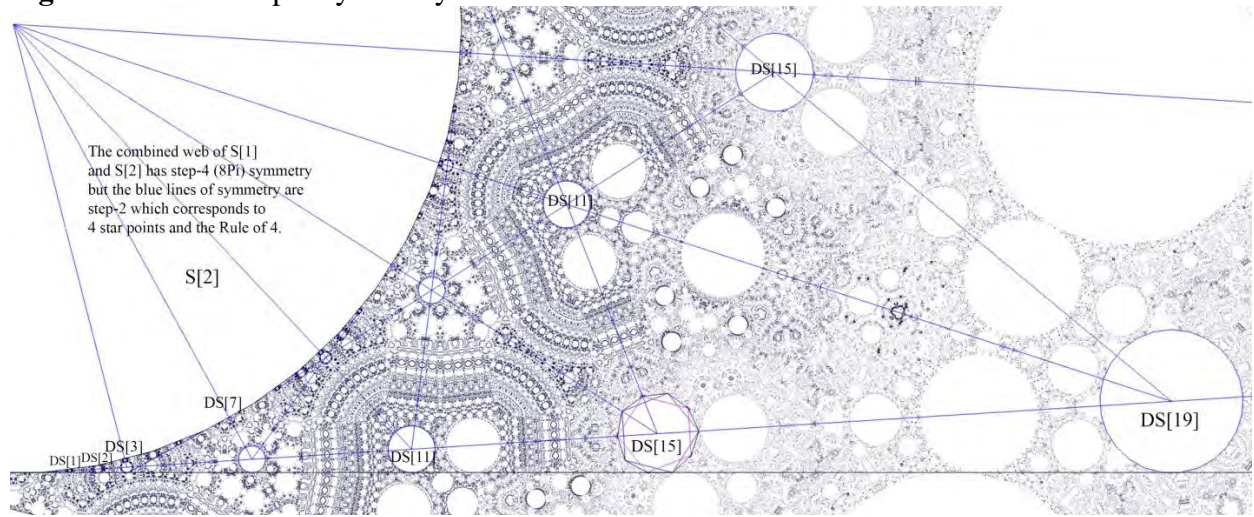

The Edge Conjecture predicts DS[k] at step-4 counting down from S[1] at DS[23] and these match the blue lines of symmetry above. It is not unusual for predicted DS[k] to have matching 'binaries' as seen here with DS[11] and DS[19], but these volunteers are seldom one of the 3 possible DS[k] that could co-exist between pairs of DS[k] so their parameters are not obvious. DS[15] is mutated because $k' = N/2-k = 25-15 = 10$. Since DS[15] is an N/2-gon the mutation condition is $\gcd(25,10) > 2$, and here $25/(\gcd(25,10) = 5$ s0 DS[15] will be the weave to blue and magenta pentagons as shown above.

A salient feature of all 8k+2 N-gons is the fact that the predicted DS[3] will so-exist with right and left side DS[2]s which can serve as D[2]s. Each DS[3] 'cluster' will also contain 3 matching

M[2]s so each D[2] will have right and left side M[2]s (with one M[2] shared). There is ample evidence that this two-family web evolution would not occur without DS[3] providing overall 3-edge structure to form a semi-invariant region on the edges of S[2]. The shared M[2] will always share a vertex with DS[3].

The 8k+2 Conjecture implies that there should be a total of 2k of these 'clusters' where k is 6 here. These 12 clusters yields 24 D[2]s and there is another D[2] on the line of symmetry as shown below, for a total of N/2 D[2]s for each D[1]. Therefore the initial conditions for the $D_k$ periods of the D[k] will be $D_1 = N/2$ (which we call n). Since D[1] = S[2] always has period n, $D_2$ will be $n^2$. These will be the initial conditions for the second-order difference equation as derived in Section 2 of Part 1. This part is easy but the critical step is obtaining a correct count for the D[3] surrounding each D[2].

To first approximation the number of D[3] for each D[2] would be just another power of n and jndeed each D[2] does support n D[3]s on its edges for a total of $n^3$ D[3]s. But the D[2] also have potential right and left side D[3] at their GenStar points. For D[1] these tiles were in different dynamical families, but for D[2] they are part of the count. For the D[2]s in clusters, the GenStar D[3]s will be shared so each D[2] will only contribute one D[3], but the lone D[2] on the line of symmetry will contribute 2 D[3] for a total of n+1 'outlier' D[3]s for each D[1], This says that the total number of D[3] will be $n^3 + n(n+1)$ so for N = 50 the count will be $25^2 \cdot 25 + 25 \cdot 26$ = 16275. This is what we call $D_3$. Assuming this continues recursuvely the formula would be

$D_k = nD_{k-1} + (n+1)D_{k-2}$ where n = N/2 and $D_1 = n$ and $D_2 = n^2$ with solution $D[n,k] = -\dfrac{n\left((-1)^k - (1+n)^k\right)}{2+n}$

**Figure 50.3** Each D[2] on the edges of D[1] contributes just one D[3] outlier except for the lone D[2] on the blue line of symmetry. This contributes 2 D[3], so the total number of 'outliers is n+1 and this is true for the n D[1]s, for a total of n(n+1) D[3] outliers.

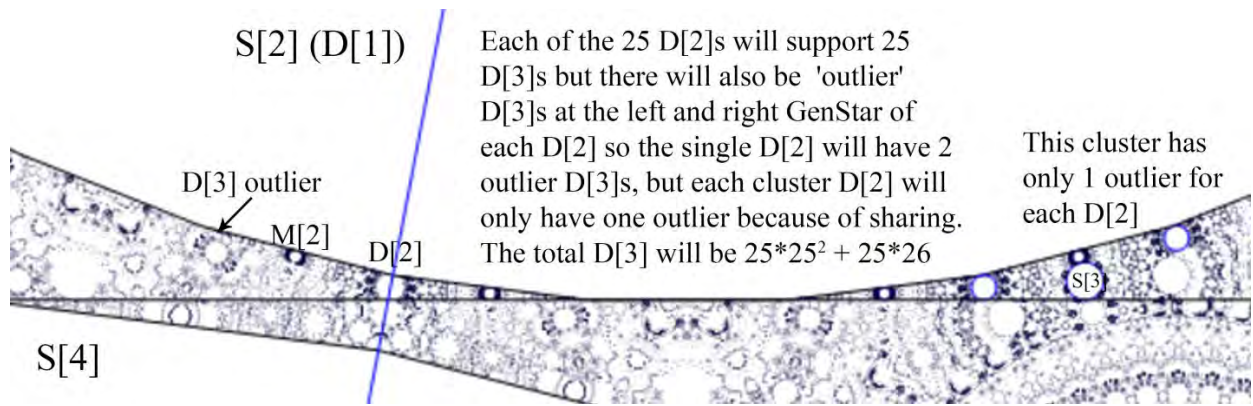

S[2] (D[1])

D[3] outlier
M[2]  D[2]

S[4]

Each of the 25 D[2]s will support 25 D[3]s but there will also be 'outlier' D[3]s at the left and right GenStar of each D[2] so the single D[2] will have 2 outlier D[3]s, but each cluster D[2] will only have one outlier because of sharing. The total D[3] will be $25*25^2 + 25*26$

This cluster has only 1 outlier for each D[2]

S[3]

The first few D[k] periods are **Table**[D[25,k],{k,1,5}] = {25, 625, 16275, 423125, 11001275}

**Figure 50.4** Detail of the 3$^{rd}$ generation presided over by D[2] and M[2]. This generation is clearly not self-similar to any of the first 2 generations.

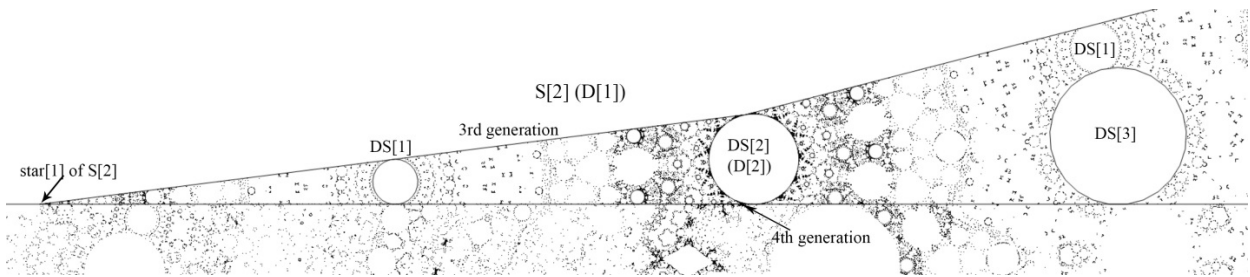

**Figure 50.5** The local edge geometry of each D[2] is self similar to the edge geometry of D[1] so D[2] can recursively play the part of a D[1].

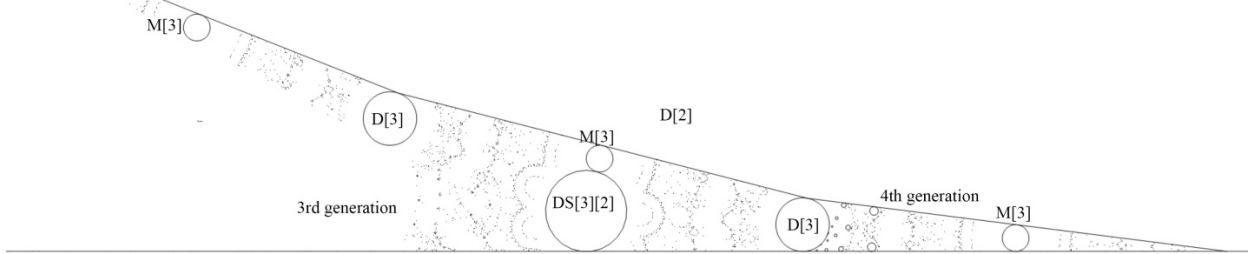

**Figure 50.6** A comparison of the 4$^{th}$ generation and the 2$^{nd}$ generation from Figure 50.1 shows that they are destined to be very different even though both D[1] and D[3] have step-4 webs that evolve in a ccw fashion and these webs support DS[k] at mod-4.

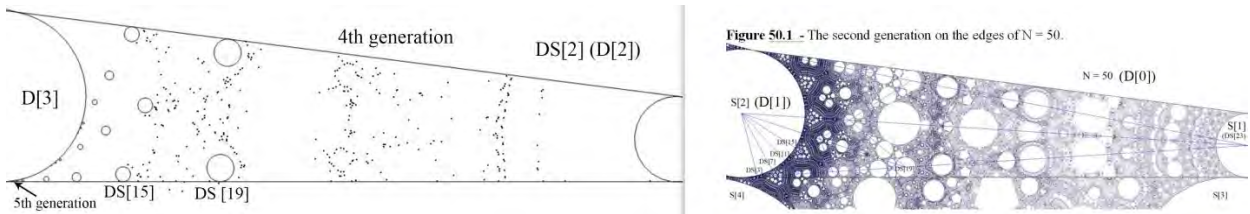

This means that we have only a vague notion of what to expect for the 5$^{th}$ generation geometry but in all cases the periods of the D[k] and the matching DS[k] are known. Here for N = 50 the periods are D[k] = -25[(-1)$^k$ - 26$^k$]/27 so the first few are 25, 625, 16275, 423125, 11001275. This implies an expected return 'time' of about 11 million iterations for the 5$^{th}$ generation so it will be a daunting task to generate even fuzzy images.

If these families presided over by the D[k] and M[k] are distinct there will be a doubly infinite series of families in the 8k+2 family. They will be embedded in a multi-fractal environment with increasing algebraic complexity so there could be some spectacular small-scale geometry. It is likely that for all regular N-gons, this 'edge' geometry will give insight into the global geometry of the web W. In cosmological terms every N-gon can be regarded as defining a distinct mode-B oscillation of the CBM . These oscillations are still poorly understood but it is likely that their 'footprint' has evolved over time through interaction with the expanding universe.

**Links**

(i) The author's web site at DynamicsOfPolygons.org is devoted to the outer billiards map and related maps from the perspective of a non-professional. The items below are available there.

(ii) A Mathematica notebook called FirstFamily.nb will generate the First Family and related star polygons for any regular polygon. It is also a full-fledged outer billiards notebook which works for all regular polygons. This notebook includes the Digital Filter map and the Dual Center map. The default height is 1 to make it compatible with the Digital Filter map. This notebook is not necessary to implement the Digital Filter or Dual Center maps, but it may be useful to have a copy of the matching First Family to be used as reference.

(iii) Outer Billiards notebooks for all convex polygons (height 1 convention for regular cases). There are four cases: Nodd, NTwiceOdd, NTwiceEven and Nonregular.

(iv) For someone willing to download the free Mathematica CDF reader there are many 'manipulates' that are available at the Wolfram Demonstrations site - including an outer billiards manipulate of the author and two other manipulates based on the author's results in [H2].